\documentstyle[12pt]{article} 
\makeatletter
\def\mineappendix{
        \setcounter{section}{1}
        \setcounter{subsection}{0}
        \def\thesection{\Alph{section}}
        \def\sectionap{\@startsection  {section}{1}{\z@}
                        {-3.5ex plus-1ex minus-.2ex} {0ex plus.2ex}
                        {\reset@font\Large\bf  Appendix:  \, }
                        } 
        }
\makeatother
\def\Proclaim #1. #2\par{\bigbreak\noindent{\sc#1.\enspace}{\it#2}\par}
\newtheorem{lem}{Lemma}[section]
 
\newtheorem{thm}[lem]{Theorem} 
\newtheorem{pro}[lem]{Proposition}

\newtheorem{conj}[lem]{Conjecture}

\title{Virasoro Constraints For Quantum Cohomology}
\author{Xiaobo Liu, Gang Tian}
\date{}

\begin{document}
\maketitle

\input amssym.def

\newcommand{\R}{\mbox{${\Bbb R}$}}
\newcommand{\om}{\mbox{$\omega$}}
\newcommand{\la}{\mbox{$\lambda$}}
\newcommand{\ffi}{\mbox{$\varphi$}}
\newcommand{\ar}{\mbox{$\alpha$}}

In \cite{EHX2}, Eguchi, Hori and Xiong, proposed a conjecture that
the partition function of topological sigma model coupled to gravity
is annihilated by infinitely many differential operators which form 
half branch of the Virasoro algebra. A similar conjecture was also proposed
by S. Katz \cite{Ka} (See also \cite{EJX}). Assuming this conjecture is
true, they were able to reproduce certain instanton numbers of some projective
spaces known before (cf. the above cited references and \cite{EX}
for details). This conjecture is also referred to as the Virasoro conjecture
by some authors. The main purpose of this paper is to give a proof of 
this conjecture for the genus zero part. 

The theory of topological sigma model coupled to gravity
has been extensively studied recently by both mathematicians and physicists.
This theory is built on the intersection theory of moduli spaces
of stable maps from Riemann surfaces to a fixed manifold $V^{2d}$,which is
 a smooth projective variety (or more  generally, a  symplectic
manifold). 
To each cohomology class of $V$ (denoted by ${\cal O}$) and a non-negative
integer $n$, there is associated a quantum field theory operator, denoted
by $\tau_{n}({\cal O})$. When $n = 0$, the corresponding operator is
simply denoted by ${\cal O}$ and is called a {\it primary field}. 
For $n > 0$, $\tau_{n}({\cal O})$ is called the n-th {\it (gravitational)
descendent} of ${\cal O}$.
The so called $k$-point
genus-$g$ {\it correlators} in topological field theory can be defined
via the {\it Gromov-Witten invariants} as follows:
\begin{eqnarray*}
&& \left< \tau_{n_{1}}({\cal O}_{1})\tau_{n_{2}}({\cal O}_{2}) \cdots
	\tau_{n_{k}}({\cal O}_{k}) \right>_{g}  \\
& := &	\sum_{A \in H_{2}(V, {\Bbb Z})} q^{A} 
	\int_{\left[\overline{\cal M}_{g, k}(V, A) \right]^{\rm virt}}
	c_{1}(E_{1})^{n_{1}} \cup {\rm ev}_{1}^{*}({\cal O}_{1}) \cup
		\cdots \cup
	c_{1}(E_{k})^{n_{k}} \cup {\rm ev}_{k}^{*}({\cal O}_{k}),
\end{eqnarray*}
where $q^{A}$ belongs to the Novikov ring (i.e. the multiplicative ring
spanned by monomials $q^{A} = q_{1}^{a_{1}} \cdots q_{r}^{a_{r}}$
over the ring of rational numbers, where
$\{q_{1}, \cdots, q_{r}\}$ is a fixed basis of $H_{2}(V, {\Bbb Z})$ and
$A = \sum_{i=1}^{r} a_{i}q_{i}$),
$\left[\overline{\cal M}_{g, k}(V, A) \right]^{\rm virt}$
 is the virtual
moduli space of degree $A$ stable maps from $k$-marked genus-$g$ curves 
to $V$ (cf. \cite{LT1}),  
$c_{1}(E_{i})$ is the first
Chern class of the tautological line bundle $E_{i}$ over
$\left[\overline{\cal M}_{g, k}(V, A) \right]^{\rm virt}$
whose fiber over each stable map
is defined by the cotangent space of the underlying curve at the
$i$-th marked point, 
 and ${\rm ev}_{i}$ is the evaluation
map from $\left[\overline{\cal M}_{g, k}(V, A) \right]^{\rm virt}$
to $V$ defined by evaluating each stable map
at the $i$-th marked point. 
We also refer to \cite{RT2} for more discussions in the case of 
semi-positive symplectic manifolds, which include 
all Fano Manifolds and Calabi-Yau manifolds as special cases.

All genus-$g$ correlators can be assembled into a generating function,
called the genus-$g$ {\it free energy function}, in the following way:
\[ F_{g}(T) := \left< \exp \sum_{n, \alpha}t^{\alpha}_{n} 
	\tau_{n}({\cal O}_{\alpha}) \right>_{g} 
	= \sum_{\{k_{n, \alpha}\}} \left(\prod_{n, \alpha} 
	  \frac{(t^{\alpha}_{n})^{k_{n, \alpha}}}{k_{n, \alpha}!}\right)
	  \left< \prod_{n, \alpha} 
		(\tau_{n}({\cal O}_{\alpha}))^{k_{n, \alpha}}
	  \right>_{g}, \]
where ${\cal O}_{1}, \cdots, {\cal O}_{N}$ form a basis of 
$H^{*}(V, {\Bbb Q})$, $\alpha$ ranges from $1$ to $N$, $n$ ranges over
the set of all non-negative integers ${\Bbb Z}_{+}$, and $\{k_{n, \alpha}\}$
ranges over the set of all collections of non-negative integers, almost
all (except finite number) of them are zero, labeled by $n$ and $\alpha$,
and $T = \{t^{\alpha}_{n} \mid n \in {\Bbb Z}_{+}, \alpha = 1, \cdots, N\}$
is an infinite set of parameters. The space of all parameters $T$
is called the {\it big phase space}. Its subspace 
$\{T \mid t^{\alpha}_{n} = 0 \, \, \, {\rm for \, \, \, all} \, \, \,
	n > 0\}$ 
is called the {\it small phase space}. The genus zero free energy $F_{0}$
restricted to the small phase space is the potential function
of  the {\it Quantum cohomology}
of $V$, whose third derivatives define the quantum ring structure
on $H^{*}(V, {\Bbb Q})$.  The generating function of 
all free energy functions, i.e.
\[ Z(T; \lambda) := \exp \sum_{g \geq 0} \lambda^{2g-2} F_{g}(T), \]
is called the {\it partition function} and $\lambda$ is called the
{\it genus expansion parameter}. 

It is widely expected that the partition function $Z$ has many interesting 
properties. For example, it always satisfies the (generalized)
string equation and dilaton equation (cf. \cite{W2}, \cite{DW}). 
When $V$ is just a point, it was 
conjectured by
Witten \cite{W2} and proved by Kontsevich \cite{Kon} that $Z$ is a 
$\tau$-function of the KdV hierarchy. On the other hand, it is well known
that the $\tau$-function of the KdV hierarchy which satisfies the string 
equation is annihilated by a sequence of differential operators, which
form half branch of the Virasoro algebra (cf. \cite{DVV}, \cite{FKN}, and
\cite{KS}). For general manifold $V$, it is not clear at this stage what
kind of integrable system might govern $Z$. However it seems very promising
that an analogue of the Virasoro constraints could still exist.
In \cite{EHX2}, Eguchi, Hori and Xiong constructed a sequence of
linear differential operators, denoted by $L_{n}$ with $n \in {\Bbb Z}$, 
on the big phase space 
(see section~\ref{sec:VirOp} for the precise form of these operators).  
They  checked that these operators define a representation of the Virasoro
algebra with the central charge equal to the Euler characteristic number
of $V$, i.e. the commutators of these operators satisfy the following relation
\begin{equation} \label{eqn:vir}
	[L_{m}, L_{n}] = (m-n)L_{m+n} + \delta_{m, -n} 
			\frac{m(m^{2}-1)}{12} \cdot \chi(V)
\end{equation}
if the following condition is satisfied:
\begin{equation} \label{eqn:chernclass}
	\frac{1}{4} \sum_{\alpha = 1}^{N} b_{\alpha}(1-b_{\alpha})
	= \frac{1}{24} \left( \frac{3-d}{2} \chi(V) -
				\int_{V} c_{1}(V) \wedge c_{d-1}(V) \right),
\end{equation}
where $d$ is equal to half of the (real) dimension of $V$,
$b_{\alpha}= \frac{1}{2}{\rm dim}({\cal O}_{\alpha}) - \frac{1}{2}(d-1)$, 
and $c_{i}(V)$ is the $i$-th Chern class of $V$.
Condition (\ref{eqn:chernclass}) is needed in order that 
$[L_{-1}, L_{1}] = L_{0}$. 
The following conjecture was proposed for Fano manifolds with only
even dimensional cohomology classes 
(See also \cite{EJX} for a more general conjecture)
\begin{conj}[Eguchi-Hori-Xiong, Katz] \label{conj:vir}
	$L_{n} Z \equiv 0$ for all $n \geq -1$.
\end{conj}
We will call equation $L_{n} Z = 0$ the $L_{n}$ constraint. 
The $L_{-1}$ constraint is the string equation. 
The $L_{0}$ constraint is a combination of the selection rule, the divisor
equation and the dilaton equation. All these equations hold for general
manifold $V$ (cf. \cite{RT2} \cite{W2}, as well as \cite{G2}).
Moreover, due to the Virasoro type relation (\ref{eqn:vir}), 
if $L_{1}$ and $L_{2}$
constraints are true, then $L_{n}$ constraint is true for all $n > 0$. 

If we write $(L_{n}Z)/Z$ as a Laurent series in $\lambda$, where $\lambda$
is the genus expansion parameter, then each $L_{n}$ constraint gives
a sequence of differential equations for the free energy functions $F_{g}$,
corresponding to the coefficients of different powers of $\lambda$.
Notice that these differential equations are no longer linear when $n > 0$
since they contain some quadratic terms.
The coefficient of $\lambda^{-2}$ gives a differential equation which only
involves genus-$0$ free energy $F_{0}$. We call this equation the
{\bf genus-$0$ $L_{n}$ constraint}. If this equation holds, we also say 
that $F_{0}$ satisfies the $L_{n}$ constraint. 
The main result of this paper can be stated as
\begin{thm} \label{thm:main}
	If $V$ has only even dimensional cohomology classes (or if we
only consider even dimensional cohomology classes in the topological
sigma model), then the genus-$0$ free energy $F_{0}$ satisfies the 
$L_{1}$ and $L_{2}$ constraints.
\end{thm}
{\bf Remark}: {\it

	(1)In this theorem, we do not assume that $V$ is Fano. In fact,
	   we even do not assume that $V$ is algebraic. All what are
	   needed in the proof of this theorem are the string equation, the
	   dilaton equation,
	   the genus-0 topological recursion relation, and Hori's $L_{0}$
	   constraint, which in turn follows from the selection rule and
	   the divisor equation (see section~\ref{sec:specialVF} and
	   \ref{sec:TRR} for precise forms of these equations). 
		Therefore this theorem
	   should be true for all manifolds where these equations hold, e.g.
	   smooth projective varieties and semi-positive symplectic manifolds.

	(2)In this theorem, we also do not assume condition 
	   (\ref{eqn:chernclass}), which is needed to guarantee that
	   $[L_{-1}, L_{1}] = L_{0}$. The reason behind this is that the
	   constant term in the $L_{0}$ operator does not affect the genus-$0$
	   constraints. As it was pointed out in \cite{Bor},
	   if $V$ has only even dimensional cohomology classes, then
	   condition (\ref{eqn:chernclass}) is equivalent to $h^{p, q}(V)=0$
	   for $p \neq q$, where $h^{p,q}(V)$ is the Hodge number of $V$.  

	(3)As we mentioned above, as long as
	   the Virasoro relation (\ref{eqn:vir}) holds for $m, n >0$,
	   this theorem implies
	   the genus-$0$ $L_{n}$ constraint for all $n>0$.
	   Consequently,  the genus-$0$ part of 
	   Conjecture \ref{conj:vir} is true.

	(4)The assumption that $V$ has only even dimensional cohomology
	   classes is not essential. The general case may be treated
	   by the same method. However, in this paper, we only
	   consider  this case  for simplicity.
}

This theorem is a combination of Proposition~\ref{pro:g0L1} and
Proposition~\ref{pro:g0L2}, which will be proved in section~\ref{sec:g0L1}
and section~\ref{sec:g0L2} respectively.
The main idea of our proof can be described as follows:
If we consider
the first derivative part of Eguchi, Hori, and Xiong's $L_{n}$ operator
as a vector field, denoted by ${\cal L}_{n}$, on the big phase space, then
we have the equation:
\begin{eqnarray}
        && \sum_{\sigma, \rho = 1}^{N}
        \left<\left< {\cal L}_{n} \left({\cal L}_{0} - (n+1){\cal D}\right)
                {\cal O}_{\sigma} \right>\right>_{0} \eta^{\sigma \rho}
        \left<\left< {\cal O}_{\rho} \tau_{k}({\cal O}_{\mu})
                \tau_{l}({\cal O}_{\nu}) 
                \right>\right>_{0} \nonumber \\
        &=& \sum_{\sigma, \rho = 1}^{N}
        \left<\left< {\cal L}_{n} \tau_{k}({\cal O}_{\mu})
                {\cal O}_{\sigma} \right>\right>_{0} \eta^{\sigma \rho}
        \left<\left< {\cal O}_{\rho} \left({\cal L}_{0} - (n+1){\cal D}\right)
                \tau_{l}({\cal O}_{\nu}) 
                \right>\right>_{0},   \label{eqn:sWDVV}
\end{eqnarray}
where ${\cal D}$ is the dilaton vector field defined in 
section~\ref{sec:specialVF},  $\left<\left< \cdot \cdot \cdot
 \right>\right>_{0}$  
is the 3-point genus-$0$ correlation function which is a symmetric tensor
on the  big phase space defined by third derivatives of $F_{0}$, and
$(\eta^{\sigma\rho})$
is the inverse matrix of the intersection form on $H^{*}(V, {\Bbb Q})$. 
This is a simple corollary of the {\it generalized WDVV equation}:
\[ \sum_{\sigma, \rho = 1}^{N}
	\frac{\partial^{3} F_{0}}{\partial t^{\alpha}_{m}
	\partial t^{\beta}_{n}
	\partial t^{\sigma}_{0}} \eta^{\sigma \rho}
	\frac{\partial^{3} F_{0}}{\partial t^{\rho}_{0}
	\partial t^{\mu}_{k}
	\partial t^{\nu}_{l}} 
  = \sum_{\sigma, \rho = 1}^{N}
	\frac{\partial^{3} F_{0}}{\partial t^{\alpha}_{m}
	\partial t^{\mu}_{k}
	\partial t^{\sigma}_{0}} \eta^{\sigma \rho}
	\frac{\partial^{3} F_{0}}{\partial t^{\rho}_{0}
	\partial t^{\beta}_{n}
	\partial t^{\nu}_{l}},
\] 
which is satisfied by the genus-$0$ free energy function (cf. \cite{W2}).  
If $F_{0}$ satisfies the $L_{n}$ constraint, then we can compute
both sides of equation (\ref{eqn:sWDVV})
by using the genus-$0$ $L_{0}$ and $L_{n}$ constraints and
the dilaton equation.
The result for the left hand side is an expression which contains
infinitely many terms, while the result for the right hand side
only contains finitely many terms. Our crucial observation is that 
the difference of these
two expressions is the second derivative,
i.e. $\frac{\partial}{\partial t^{\nu}_{l}} 
        \frac{\partial}{\partial t^{\mu}_{k}}$, 
of a function which does not
depend on $\tau_{k}({\cal O}_{\mu})$ and 
$\tau_{l}({\cal O}_{\nu})$. Moreover, up to some
linear terms, this function is just
the coefficient of $\lambda^{-2}$ in the Laurent expansion of
$(L_{n+1}Z)/Z$. Vanishing of this function is the genus-$0$ $L_{n+1}$
constraint. This observation provides us with a general strategy
for proving the genus-$0$ Virasoro constraints, which will be described
in more detail in section~\ref{sec:VirOp}.  Such a strategy
could be easily adapted to prove many other constraints, as it will be
demonstrated in section~\ref{sec:tilde}.

We would like to mention that in \cite{EHX2}, a heuristic argument
for deriving the genus-$0$ $L_{1}$ constraint for ${\Bbb C}P^{n}$ was given.
This argument is based on a recursion formula, called fundamental
recursion relation, which was discovered in \cite{EHX1} (see also equation
(\ref{eqn:FRR})). However 
it seems that there is a serious gap in this argument.
What was really proved in \cite{EHX2} is the following:
\[
	\frac{\partial}{\partial t^{\alpha}_{0}}
	\left( b_{\alpha}^{2} \widetilde{\Psi}_{0, 1}
		+ b_{\alpha} \widetilde{\Psi}_{0, 1}
		+ \Psi_{0, 1} \right) = 0,
\]
for every $\alpha = 1, \ldots, N$, where $\Psi_{0, 1}$ and
$\widetilde{\Psi}_{0, 1}$ are two functions on the big phase
space which do not depend on $\alpha$.
Integrating this equation with respect to $t^{\alpha}_{0}$
 and assuming that the integration constant
is zero, one obtains
\[
	b_{\alpha}^{2} \widetilde{\Psi}_{0, 1}
		+ b_{\alpha} \widetilde{\Psi}_{0, 1}
		+ \Psi_{0, 1} = 0.
\]
If ${\rm dim} H^{*}(V, {\Bbb Q}) \geq 3$, this equation would imply
$\Psi_{0, 1} = 0$, which is equivalent to the genus-$0$ $L_{1}$
constraint, and $\widetilde{\Psi}_{0, 1} = 0$, which is a new constraint
called the $\widetilde{L}_{1}$ constraint (which  will also be proved in
section~\ref{sec:tilde} of this paper). 
However, in this procedure, it is not clear why the integration 
constant, which still depends on infinitely many other parameters,
should be zero. It seems that to prove the vanishing of the 
integration constant
is as difficult as to prove the $L_{1}$ constraint itself.

This paper is organized as follows. In section~\ref{sec:Corr}, we first define
the basic notations used in this paper. We then review some well known
facts about correlation functions and derive some simple but very useful 
applications of these facts. Virasoro operators
of Eguchi, Hori, and Xiong are introduced in section~\ref{sec:VirOp}.
We then give the precise interpretation of Conjecture~\ref{conj:vir}
for free energy functions. At the end of section~\ref{sec:VirOp}, we
describe a general strategy for using the generalized WDVV equation to prove
the genus-$0$ part of Conjecture~\ref{conj:vir}. This strategy is carried
out for $L_{1}$ and $L_{2}$ constraints in section~\ref{sec:g0L1}
and section~\ref{sec:g0L2} respectively.
In section~\ref{sec:tilde}, we prove two other genus-$0$ constraints,
called $\widetilde{L}_{1}$ and $\widetilde{L}_{2}$ constraints, 
which were also conjectured in \cite{EHX2}. We will discuss higher
genus cases in a forthcoming paper.

The first author is partially supported by an NSF postdoctoral fellowship.

\section{Relations among different correlation functions} 
\label{sec:Corr}

In this section, we review some well known formulas
for correlation functions and derive some of their immediate consequences.
We will always identify quantum field theory operators 
$\tau_{m}({\cal O}_{\alpha})$
with the tangent vector fields $\frac{\partial}{\partial t^{\alpha}_{m}}$
and view the {\it genus-$g$ correlation functions}, denoted
by 
$\left<\left< \, \right>\right>_{g}$, 
as symmetric tensors on the big phase
space defined by
\[ \left<\left< \tau_{m_{1}}({\cal O}_{\alpha_{1}})
	\tau_{m_{2}}({\cal O}_{\alpha_{2}}) 
	\cdots \tau_{m_{k}}({\cal O}_{\alpha_{k}}) 
	 \right>\right>_{g} := 
	\frac{\partial^{k}}{\partial t^{\alpha_{1}}_{m_{1}}
	\partial t^{\alpha_{2}}_{m_{k}} \cdots
	\partial t^{\alpha_{k}}_{m_{k}}} F_{g},
\]
where $F_{g}$ is the genus-$g$ free energy function.

\subsection{Convention of notations}
\label{sec:notation}
We will use the following convention of notations
throughout the paper unless otherwise stated.
We will use $d$ to denote one half of the real dimension of $V$.
$N$ is the dimension of the space of cohomology classes $H^{*}(V, {\Bbb Q})$.
Lower case Greek letters, e.g. $\alpha$, $\beta$, $\gamma$, ..., etc., will
be used to index the cohomology classes. The range of these indices is from
$1$ to $N$. Lower case English letters, e.g. $i$, $j$, $k$, $m$, $n$,
..., etc., will
be used to index the level of gravitational descendents. Their range is the
set of all non-negative integers, i.e. ${\Bbb Z}_{+}$.
All summations are over the entire ranges of the indices unless otherwise
indicated.
We fix a basis ${\cal O}_{1}$, ${\cal O}_{2}$, ..., ${\cal O}_{N}$ of 
$H^{*}(V, {\Bbb Q})$ which is arranged in such an order  that the 
dimension of ${\cal O}_{\alpha}$ is non-decreasing with respect to $\alpha$.
In particular, ${\cal O}_{1}$ is equal to the identity element
of the ordinary cohomology ring. Gravitational descendents are denoted by
$\tau_{m}({\cal O}_{\alpha})$ whose corresponding parameters are 
$t^{\alpha}_{m}$, where $m \in {\Bbb Z}_{+}$ and $\alpha = 1, \ldots, N$.
$\tau_{0}({\cal O}_{\alpha})$ is always identified with ${\cal O}_{\alpha}$.
We also consider $\tau_{m}({\cal O}_{\alpha})$ with $m<0$ as a zero operator.
Let $\eta_{\alpha \beta} = \int_{V} {\cal O}_{\alpha} \cup {\cal O}_{\beta}$
be the intersection form on $H^{*}(V, {\Bbb Q})$. 
We will use $\eta = (\eta_{\alpha \beta})$ and 
$\eta^{-1} = (\eta^{\alpha \beta})$ to lower and raise indices. Let
${\cal C} = ({\cal C}_{\alpha}^{\beta})$ be the matrix of multiplication
by the first Chern class $c_{1}(V)$ in the ordinary cohomology ring, i.e.
\begin{equation} \label{eqn:chernmatrix}
 c_{1}(V) \cup {\cal O}_{\alpha} = \sum_{\beta} {\cal C}_{\alpha}^{\beta}
		{\cal O}_{\beta}. 
\end{equation}
Since we are dealing with even dimensional cohomology classes only,
both $\eta$ and ${\cal C} \eta$ are symmetric matrices, where the entries
of ${\cal C}\eta$ are given by
${\cal C}_{\alpha \beta} = \int_{V}c_{1}(V) \cup {\cal O}_{\alpha} 
	\cup {\cal O}_{\beta}$.
Let $q_{\alpha} = (1/2){\mathrm dim}({\cal O}_{\alpha})$ and
\begin{equation} \label{eqn:balpha}
	 b_{\alpha} = q_{\alpha} - \frac{1}{2}(d-1).
\end{equation}
The following simple observations will be used throughout the calculations
without mentioning: If $\eta^{\alpha \beta} \neq 0$ or
$\eta_{\alpha \beta} \neq 0$, then $b_{\alpha} = 1-b_{\beta}$.
${\cal C}_{\alpha}^{\beta} \neq 0$ implies $b_{\beta} =
1 + b_{\alpha}$, and ${\cal C}_{\alpha\beta} \neq 0$ implies 
$b_{\beta} = - b_{\alpha}$.

Instead of coordinates 
$\{t^{\alpha}_{m} \mid m\in {\Bbb Z}_{+}, \, \alpha = 1, \ldots, N\}$, 
it is very convenient to use the following shifted
coordinates on the big phase space
\begin{equation} \label{eqn:ttilde}	
	\tilde{t}^{\alpha}_{m} = 
	t^{\alpha}_{m} - \delta_{m, 1} \delta_{\alpha, 1}
	= \left\{ \begin{array}{ll}
			t^{\alpha}_{m}-1, & \textrm{ if }m=\alpha=1, \\
			t^{\alpha}_{m}, & \textrm{ otherwise.}
		  \end{array}
		\right. 
\end{equation}
Notice that these two coordinate systems have different origins.

\subsection{Some special vector fields on the big phase space}
\label{sec:specialVF}
The first vector field, which will be used extensively later, is the following
\[ {\cal S} := - \sum_{m, \alpha} \tilde{t}^{\alpha}_{m} 
	\frac{\partial}{\partial t^{\alpha}_{m-1}} \]
We call this vector field the {\bf string vector field}.
The famous {\it string equation} (cf. \cite{RT2} and \cite{W2})
can be expressed as
\[ \left<\left< {\cal S} \right>\right>_{g} = {\cal S} F_{g} =
	\frac{1}{2} \delta_{g, 0} \sum_{\alpha, \beta} \eta_{\alpha \beta} 
		t^{\alpha}_{0} t^{\beta}_{0}. \]
This equation is equivalent to Eguchi, Hori, and Xiong's $L_{-1}$ 
constraint. 
Using this equation and the fact that 
$[{\cal S}, \frac{\partial}{\partial t^{\alpha}_{m}}] 
	= \frac{\partial}{\partial t^{\alpha}_{m-1}}$, 
we can show the following
\begin{lem} \label{lem:StringCorr}
\begin{eqnarray*}
&(1)&  \left<\left< {\cal S} \right>\right>_{0} =
	\frac{1}{2}  \sum_{\alpha, \beta} \eta_{\alpha \beta} 
		t^{\alpha}_{0} t^{\beta}_{0}.    \\
&(2)&	\left<\left< {\cal S} \tau_{m}({\cal O}_{\alpha}) \right>\right>_{0} =
		\left<\left< \tau_{m-1}({\cal O}_{\alpha}) \right>\right>_{0}
		+ \delta_{m, 0} \sum_{\beta} \eta_{\alpha \beta}
		t^{\beta}_{0}.   \\
&(3)& \left<\left< {\cal S} \tau_{m}({\cal O}_{\alpha}) 
		\tau_{n}({\cal O}_{\beta}) \right>\right>_{0} =
		\left<\left< \tau_{m}({\cal O}_{\alpha}) 
			\tau_{n-1}({\cal O}_{\beta})\right>\right>_{0}
		+ \left<\left< \tau_{m-1}({\cal O}_{\alpha}) 
			\tau_{n}({\cal O}_{\beta})\right>\right>_{0}
		+ \delta_{m, 0} \delta_{n, 0} \eta_{\alpha \beta}. 
\end{eqnarray*}
\end{lem}

Another special vector field is
\[ {\cal D} := - \sum_{m, \alpha} \tilde{t}^{\alpha}_{m}
	\frac{\partial}{\partial t^{\alpha}_{m}}. \]
We call ${\cal D}$ the {\bf Dilaton vector field}. Notice that
some authors call $\tau_{1}({\cal O}_{1})$ the dilaton operator
which is different from ${\cal D}$. The name for ${\cal D}$ is justified
by the so called {\it dilaton equation}, which can be expressed as
\[ \left<\left< {\cal D} \right>\right>_{g} = {\cal D} F_{g} =
	(2g-2) F_{g} + \frac{1}{24} \, \chi(V)\delta_{g, 1}. \]
Using this equation and the fact that 
$[{\cal D}, \frac{\partial}{\partial t^{\alpha}_{m}}] 
	= \frac{\partial}{\partial t^{\alpha}_{m}}$, 
we can show the following
\begin{lem}  \label{lem:DilatonCorr}
\begin{eqnarray*}
&(1)& \left<\left< {\cal D} \right>\right>_{0} = -2 F_{0}.  
			\textrm{ \hspace{240pt}} \\
&(2)&\left<\left< {\cal D} \tau_{m}({\cal O}_{\alpha}) \right>\right>_{0} =
		-\left<\left< \tau_{m}({\cal O}_{\alpha}) \right>\right>_{0}.
			   \\
&(3)& \left<\left< {\cal D} \tau_{m}({\cal O}_{\alpha}) 
		\tau_{n}({\cal O}_{\beta}) \right>\right>_{0} \equiv 0.
\end{eqnarray*}
\end{lem}

Probably the most important vector field in deriving the Virasoro constraints
is 
\[ {\cal X} := - \sum_{m, \alpha} \left(m + b_{\alpha} - \frac{3-d}{2}
			\right)\tilde{t}^{\alpha}_{m}
	\frac{\partial}{\partial t^{\alpha}_{m}} 
	- \sum_{m, \alpha, \beta} 
	{\cal C}_{\alpha}^{\beta}\tilde{t}^{\alpha}_{m}
	\frac{\partial}{\partial t^{\beta}_{m-1}}, 
	\]
where ${\cal C}$ is the matrix of multiplication by the first Chern class
defined by (\ref{eqn:chernmatrix}) and $b_{\alpha}$ is defined by
(\ref{eqn:balpha}). When restricted to the small phase space,
${\cal X}$ is the Euler vector field of the Frobenius manifold defined
by the restriction of the genus-$0$ free energy $F_{0}$ (cf. \cite{Du}).
Therefore we also call ${\cal X}$ itself the {\bf Euler vector field}.
It seems that the significance of this vector field on the big phase space
was first noticed in \cite{EHX1} where it is called the perturbed
first Chern class. As noted in \cite{EHX1}, the divisor equation for
the first Chern class $c_{1}(V)$ together with the selection rule implies
the following
\begin{lem}
\[ \left<\left< {\cal X} \right>\right>_{g} = {\cal X} F_{g} =
	(3-d)(1-g) F_{g} 
	+ \frac{1}{2} \delta_{g, 0}\sum_{\alpha, \beta} {\cal C}_{\alpha \beta}
	t^{\alpha}_{0}t^{\beta}_{0}
	- \frac{1}{24} \delta_{g, 1} \int_{V} c_{1}(V) \cup c_{d-1}(V). \]
\end{lem}
Adopting the language of Frobenius manifolds, we call this equation
the {\it quasi-homogeneity equation}.
Using this equation and the fact that 
\[ [{\cal X}, \, \frac{\partial}{\partial t^{\alpha}_{m}} ] 
	= \left(m + b_{\alpha} -\frac{3-d}{2}\right)
		\frac{\partial}{\partial t^{\alpha}_{m}}
	+ \sum_{\beta} {\cal C}_{\alpha}^{\beta}
 		\frac{\partial}{\partial t^{\beta}_{m-1}},\]
we can show the following
\begin{lem} \label{lem:EulerCorr}
\begin{eqnarray*}
&(1)&  \left<\left< {\cal X} \right>\right>_{0} = (3-d) F_{0}
		+ \frac{1}{2} \sum_{\alpha, \beta} {\cal C}_{\alpha \beta}
	t^{\alpha}_{0}t^{\beta}_{0}.    \\
&(2)&	\left<\left< {\cal X} \tau_{m}({\cal O}_{\alpha}) \right>\right>_{0} =
		\left(m + b_{\alpha} + \frac{3-d}{2} \right)
		\left<\left< \tau_{m}({\cal O}_{\alpha}) \right>\right>_{0}
		+ \sum_{\beta} {\cal C}_{\alpha}^{\beta}
	\left<\left< \tau_{m-1}({\cal O}_{\beta}) \right>\right>_{0}  \\
   &&\textrm{ \hspace{90pt}}	+ \delta_{m, 0} \sum_{\beta} {\cal C}_{\alpha \beta} t^{\beta}_{0}.
			   \\
&(3)& \left<\left< {\cal X} \tau_{m}({\cal O}_{\alpha}) 
		\tau_{n}({\cal O}_{\beta}) \right>\right>_{0} =
	\delta_{m, 0} \delta_{n,0} {\cal C}_{\alpha \beta} +
	(m+n+b_{\alpha}+b_{\beta}) \left<\left<\tau_{m}({\cal O}_{\alpha}) 
		\tau_{n}({\cal O}_{\beta}) \right>\right>_{0} \\
 &  & \textrm{ \hspace{110pt}}	+ \sum_{\gamma} {\cal C}_{\alpha}^{\gamma}
	\left<\left<\tau_{m-1}({\cal O}_{\gamma}) 
		\tau_{n}({\cal O}_{\beta}) \right>\right>_{0} 
	+ \sum_{\gamma} {\cal C}_{\beta}^{\gamma}
	\left<\left<\tau_{m}({\cal O}_{\alpha}) 
		\tau_{n-1}({\cal O}_{\gamma}) \right>\right>_{0}.
\end{eqnarray*} 
\end{lem}

Let ${\cal L}_{0} := - {\cal X} - \frac{3-d}{2}{\cal D}$. Then the dilaton
equation and the quasi-homogeneity equation imply
\begin{equation} \label{eqn:Hori}
 \left<\left< {\cal L}_{0} \right>\right>_{g} = {\cal L}_{0} F_{g} =
	 -\frac{1}{2} \delta_{g, 0}\sum_{\alpha, \beta} {\cal C}_{\alpha \beta}
	t^{\alpha}_{0}t^{\beta}_{0}
	- \frac{1}{24} \delta_{g, 1} \left(
	\frac{3-d}{2} \chi(V) - \int_{V} c_{1}(V) \cup c_{d-1}(V)
		\right). 
\end{equation}
This equation was first discovered by Hori \cite{H}. 
It is equivalent to Eguchi, Hori and Xiong's 
$L_{0}$ constraint for the partition function.

\subsection{Genus-$0$ topological recursion relation and its 
applications}
\label{sec:TRR}

Topological recursion relations make it possible to express many
correlation functions involving gravitational descendents by those only
involve primary fields. Such relations have been proven to exist
in genus-$0$ (cf \cite{RT2} and \cite{W2}) and genus $1$ and $2$ (cf.
\cite{G1}, \cite{G2} and \cite{BP}). 
In this paper we only consider the genus-$0$ case.
Genus-$0$ {\it topological recursion relation} has the following
form:
\[ \left<\left< \tau_{m}({\cal O}_{\alpha}) 
		\tau_{n}({\cal O}_{\beta})
		\tau_{k}({\cal O}_{\gamma}) \right>\right>_{0} =
	\sum_{\sigma} \left<\left< \tau_{m-1}({\cal O}_{\alpha}) 
			{\cal O}_{\sigma} \right>\right>_{0}
		\left<\left< {\cal O}^{\sigma} 
		\tau_{n}({\cal O}_{\beta}) 
		\tau_{k}({\cal O}_{\gamma}) \right>\right>_{0}, \]
for $m > 0$. In this formula, we used the convention that the indices
of primary fields are raised by $\eta^{-1}$. Therefore ${\cal O}^{\sigma}$
should be understood as $\sum_{\rho}\eta^{\sigma \rho}{\cal O}_{\rho}$. 
As noted by Witten (\cite{W2}), this recursion relation
implies the {\it generalized WDVV equation}:
\begin{eqnarray*} \sum_{\sigma} \left<\left< \tau_{m}({\cal O}_{\alpha}) 
		\tau_{n}({\cal O}_{\beta})
		{\cal O}_{\sigma} \right>\right>_{0}
		\left<\left< {\cal O}^{\sigma} 
		\tau_{k}({\cal O}_{\mu}) 
		\tau_{l}({\cal O}_{\nu}) \right>\right>_{0}  
	\textrm{ \hspace{100pt}}   \\
	\textrm{ \hspace{10pt}}
	 = \sum_{\sigma} \left<\left< \tau_{m}({\cal O}_{\alpha}) 
			\tau_{k}({\cal O}_{\mu}) 
			{\cal O}_{\sigma} \right>\right>_{0}
		\left<\left< {\cal O}^{\sigma} 
		\tau_{n}({\cal O}_{\beta}) 
		\tau_{l}({\cal O}_{\nu}) \right>\right>_{0}. 
		\textrm{ \hspace{50pt}}
\end{eqnarray*}
When restricted to the space of primary fields, this equation implies
the associativity of the algebra defined by the third
derivatives of $F_{0}$ and $\eta^{-1}$. However when gravitational
descendents are involved, the exact algebraic structure hidden in 
this equation seems not very clear. As we will see later in this paper,
the genus zero Virasoro constraints are actually disguised in this equation.

Genus-$0$ topological recursion relation is a recursion formula for
3-point functions. It can be used to derive recursion formulas for
2-point functions when combined with other equations. For example,
applying the topological recursion relation to the 3-point functions
in Lemma~\ref{lem:EulerCorr} (3),
we get
\begin{eqnarray} 
	&&\sum_{\mu, \nu}
	\left\{ \delta_{m, 0}\delta_{\mu, \alpha} +
		\left<\left< \tau_{m-1}({\cal O}_{\alpha}) 
		{\cal O}^{\mu} \right>\right>_{0}
		\right\}
	\left<\left< {\cal O}_{\mu} 
		{\cal X} {\cal O}_{\nu} \right>\right>_{0}
	\left\{ \left<\left< {\cal O}^{\nu} 
		\tau_{n-1}({\cal O}_{\beta})\right>\right>_{0}
		+\delta_{n, 0} \delta_{\nu, \beta} \right\}   
		\nonumber  \\
	&=&  \delta_{m, 0} \delta_{n,0} {\cal C}_{\alpha \beta} +
	(m+n+b_{\alpha}+b_{\beta})\left<\left< \tau_{m}({\cal O}_{\alpha}) 
		\tau_{n}({\cal O}_{\beta}) \right>\right>_{0}  
		\nonumber \\
	&& + \sum_{\sigma}{\cal C}_{\alpha}^{\sigma}
		\left<\left< \tau_{m-1}({\cal O}_{\sigma}) 
		\tau_{n}({\cal O}_{\beta}) \right>\right>_{0}
	+\sum_{\sigma}{\cal C}_{\beta}^{\sigma}
		\left<\left< \tau_{m}({\cal O}_{\alpha}) 
		\tau_{n-1}({\cal O}_{\sigma}) \right>\right>_{0}.
	\label{eqn:FRR}
\end{eqnarray}
Notice that, by Lemma~\ref{lem:EulerCorr} (3),
\begin{equation}
\left<\left< {\cal O}_{\mu} 
		{\cal X} {\cal O}_{\nu} \right>\right>_{0} =
	{\cal C}_{\mu \nu} + (b_{\mu} + b_{\nu})
	\left<\left< {\cal O}_{\mu} 
		 {\cal O}_{\nu} \right>\right>_{0},
\end{equation}
which only involves primary fields.
Therefore (\ref{eqn:FRR}) is really a recursion relation if
$m+n+b_{\alpha}+b_{\beta} \neq 0$. This
recursion relation was first noticed in \cite{EHX1}, where it was
called the fundamental recursion relation. It was also used in \cite{EHX2}
to give a heuristic argument (with some serious gaps) 
to the genus-$0$ Virasoro constraints for ${\Bbb C}P^{n}$.

Applying the topological recursion relation to the 3-point functions in
Lemma~\ref{lem:StringCorr} (3) and notice that
\[ \left<\left< {\cal O}_{\mu} 
		{\cal S} {\cal O}_{\nu} \right>\right>_{0} = \eta_{\mu \nu}, \]
we get another recursion formula:
\begin{eqnarray} \label{eqn:StringRec}
&&	\left<\left< \tau_{m}({\cal O}_{\alpha}) 
		\tau_{n-1}({\cal O}_{\beta}) \right>\right>_{0} +
	\left<\left< \tau_{m-1}({\cal O}_{\alpha}) 
		\tau_{n}({\cal O}_{\beta}) \right>\right>_{0} 
	\nonumber \\
&=&	\delta_{m,0} \left<\left< {\cal O}_{\alpha} 
		\tau_{n-1}({\cal O}_{\beta}) \right>\right>_{0} +
	\delta_{n, 0} \left<\left< \tau_{m-1}({\cal O}_{\alpha}) 
		{\cal O}_{\beta} \right>\right>_{0}   
	\nonumber \\
&&	+ \sum_{\sigma}\left<\left< \tau_{m-1}({\cal O}_{\alpha}) 
		{\cal O}_{\sigma} \right>\right>_{0}
	\left<\left< {\cal O}^{\sigma} 
		\tau_{n-1}({\cal O}_{\beta}) \right>\right>_{0}.
\end{eqnarray}
In this paper, this formula will mainly be used to shift the
level of descendents from one primary field to another. It's also
interesting to observe that sometimes it is very effective to use
this formula to reduce the level of descendents. For example,
for $m = n > 0$ and $\alpha = \beta$, this formula takes the following
simple form:
\[ \left<\left< \tau_{m}({\cal O}_{\alpha}) 
		\tau_{m-1}({\cal O}_{\alpha}) \right>\right>_{0} 
	= \frac{1}{2} 
	\sum_{\sigma}\left<\left< \tau_{m-1}({\cal O}_{\alpha}) 
		{\cal O}_{\sigma} \right>\right>_{0}
	\left<\left< {\cal O}^{\sigma} 
		\tau_{m-1}({\cal O}_{\alpha}) \right>\right>_{0}.
\]

\section{Virasoro operators}
\label{sec:VirOp}

In this section, we first give the constructions of Virasoro operators
by Eguchi, Hori, and Xiong. We then describe the relationship between
these operators and the generalized WDVV equation. This provide us
with a general strategy to prove the genus-$0$ part of the Virasoro
constraints. We will use
the normalizations in \cite{EJX} which are more consistent with \cite{RT2}
and \cite{W2}.

Define
\begin{equation}
	L_{-1} := \sum_{m, \alpha} \tilde{t}^{\alpha}_{m} 
			\frac{\partial}{\partial t^{\alpha}_{m-1}}
			+ \frac{1}{2 \lambda^{2}} \sum_{\alpha, \beta}
			\eta_{\alpha \beta} t^{\alpha}_{0} t^{\beta}_{0}, 
			\hspace{130pt}
\end{equation}
\begin{eqnarray}
	L_{0}& := &\sum_{m, \alpha} (m+b_{\alpha})\tilde{t}^{\alpha}_{m} 
			\frac{\partial}{\partial t^{\alpha}_{m}}
			+ \sum_{m, \alpha, \beta} {\cal C}_{\alpha}^{\beta}
				\tilde{t}^{\alpha}_{m} 
			\frac{\partial}{\partial t^{\beta}_{m-1}}
			+ \frac{1}{2 \lambda^{2}} \sum_{\alpha, \beta}
			{\cal C}_{\alpha \beta} t^{\alpha}_{0} t^{\beta}_{0}
			\nonumber \\
	   & &	+ \frac{1}{24}\left( \frac{3-d}{2} \chi(V) -
			\int_{V} c_{1}(V) \cup c_{d-1}(V) \right),
\end{eqnarray}
and for $n \geq 1$,
\begin{eqnarray}
	L_{n}& := &\sum_{m, \alpha, \beta} \sum_{j=0}^{m+n} 
			A^{(j)}_{\alpha}(m, n)({\cal C}^{j})_{\alpha}^{\beta}
			\tilde{t}^{\alpha}_{m} 
			\frac{\partial}{\partial t^{\beta}_{m+n-j}} 
			\hspace{100pt} 
			\nonumber \\
		&&	+ \frac{\lambda^{2}}{2}
			\sum_{\alpha, \beta, \gamma} \sum_{j=0}^{n-1}
			\sum_{k=0}^{n-j-1} B^{(j)}_{\alpha}(k, n)
			({\cal C}^{j})_{\alpha}^{\beta}\eta^{\alpha \gamma}
			\frac{\partial}{\partial t^{\gamma}_{k}}
			\frac{\partial}{\partial t^{\beta}_{n-k-1-j}}
			\nonumber \\
		&&	+ \frac{1}{2 \lambda^{2}} \sum_{\alpha, \beta}
			({\cal C}^{n+1})_{\alpha \beta}
			t^{\alpha}_{0} t^{\beta}_{0},	
\end{eqnarray}
where ${\cal C}^{j}$ is the $j$-th power of the matrix ${\cal C}$,
$({\cal C}^{n+1})_{\alpha \beta}$ are entries of the matrix
${\cal C}^{n+1}\eta$, $A^{(j)}_{\alpha}(m, n)$ and $B^{(j)}_{\alpha}(m, n)$ 
are constants defined in terms of Gamma function by
\[ A^{(j)}_{\alpha}(m, n) := \frac{\Gamma(b_{\alpha} + m + n + 1)}
			{\Gamma(b_{\alpha}+m)}
			\sum_{m \leq l_{1} < l_{2} < \cdots < l_{j} \leq m+n}
			\left( \prod_{i=1}^{j} \frac{1}{b_{\alpha} + l_{i}}
			\right),   \]
and
\[ B^{(j)}_{\alpha}(m, n) := \frac{\Gamma(m+2-b_{\alpha})\Gamma(n-m+b_{\alpha})}
			{\Gamma(1-b_{\alpha}) \Gamma(b_{\alpha})}
			\sum_{-m-1 \leq l_{1} < l_{2} < \cdots < l_{j} 
					\leq n-m-1}
			\left( \prod_{i=1}^{j} \frac{1}{b_{\alpha} + l_{i}}
			\right).   \]
When $j = 0$, the last factors in $A^{(j)}_{\alpha}(m, n)$ and
$B^{(j)}_{\alpha}(m, n)$ should be understood as equal to 1.
Any term which contains $t^{\alpha}_{m}$ with $m < 0$ should be understood as 
zero. Eguchi, Hori, and Xiong also construct $L_{-n}$ for $n>0$. However, the
significance of these operators is not clear and we do no deal with them
in this paper.

It is well known that $L_{n} Z(T; \lambda) \equiv 0$ for $n = -1$ or $0$,
where 
$T=\{t^{\alpha}_{m} \mid m \in {\Bbb Z}_{+}, \, \, \alpha = 1, \ldots, N\}$
and  $Z(T; \lambda)$ is the partition 
function defined in the introduction.
The first equation (i.e. for $n = -1$) is the string equation. 
The second equation (i.e. for $n =0$) is equivalent to (\ref{eqn:Hori}). 
The analogous equations for $n \geq 1$
is the content of Conjecture \ref{conj:vir}.
Let $\Psi_{g, n}(T)$ be the coefficient of $\lambda^{2g-2}$ in the Laurent
expansion of $(L_{n}Z(T; \lambda))/Z(T; \lambda)$.
In other words, $\Psi_{g, n}$ is defined
by 
\begin{equation} \label{eqn:genuscons}
	 L_{n}Z(T; \lambda) = 
		\left\{\sum_{g \geq 0} \Psi_{g, n} \lambda^{2g-2} \right\}
			Z(T; \lambda). 
\end{equation}
We call the equation $L_{n}Z = 0$ the {\it $L_{n}$-constraint} for the 
partition function. It is equivalent to $\Psi_{g, n} = 0$ for all $g$.
The equation $\Psi_{g, n} = 0$ will be called 
{\bf genus-$g$ $L_{n}$-constraint}. 
For $n = -1$ or $0$, this is a first order linear differential equation 
for the genus-$g$ free energy $F_{g}$. When $n \geq 1$, it is a second order
non-linear differential equation involving 
all free energy functions $F_{g^{'}}$ with $0 \leq g^{'} \leq g$.
The genus-$0$ constraints are special in the sense that only 
$F_{0}$ is involved in these equations.
It is straightforward to check the following fact:
\begin{lem}
 Suppose that the $L_{n}$ operators satisfy the Virasoro relation
\[ [L_{m}, L_{n}] = (m-n) L_{m+n} \, \, \, \, \textrm{  for } m, n \geq 1. \] 
Given $m, n \geq 1$ and $m \neq n$, if 
$\Psi_{g^{'}, m} = \Psi_{g^{'}, n} \equiv 0$ for all $g^{'}$ satisfying
$0 \leq g^{'} \leq g$, then $\Psi_{g, m+n} \equiv 0$.
\end{lem}

In this paper, we are only interested to the genus-$0$ constraints
$\Psi_{0, n} = 0$. 
We first observe that to prove the genus-$0$ $L_{n}$ constraints, it suffices
to show that all second derivatives of $\Psi_{0, n}$ vanish.
In fact, Lemma~\ref{lem:DilatonCorr} (2) and (3) at the origin
trivially imply the following:
\[ \left. \frac{\partial^{2} }{\partial t_{1}^{1} 
		\partial t_{k}^{\mu}} \Psi_{0, n} 
	\right|_{T=0} =
        - \left. \frac{\partial}{\partial t_{k}^{\mu}} \Psi_{0, n} 
		\right|_{T=0} 
	\, \, \, \, \, {\rm and} \, \, \, \, \,
	\left. \frac{\partial }{\partial t_{1}^{1}} \Psi_{0, n} \right|_{T=0} 
		=  \left. - 2 \Psi_{0, n} \right|_{T=0}. 
\]
(Same formulas also hold for $\tilde{\Psi}_{0,n}$ defined in 
section~\ref{sec:tilde}.) Therefore once we know that all secend derivatives
of $\Psi_{0, n}$ are zero, $\Psi_{0, n}$ and all of its first derivatives
have to vanish at the origin. Consequently $\Psi_{0, n}$ is constantly equal
to zero.

It is also interesting to observe that
 all the vector fields introduced in section
\ref{sec:specialVF}, i.e. ${\cal S}$, ${\cal D}$, 
and ${\cal X}$,  vanish at a very special point
\[\widetilde{T}_{0} = \{ \tilde{t}^{\alpha}_{m} = 0 \mid m \in {\Bbb Z}_{+},
\, \alpha = 1, \ldots, N\}. \] 
It follows from Lemma~\ref{lem:DilatonCorr} (2) that
all 1-point genus-$0$ correlation functions vanish at this point, i.e.,
\begin{equation}
	\left<\left< \tau_{m}({\cal O}_{\alpha}) \right>\right>_{0} 
	\mid_{\widetilde{T}_{0}} = 0
\end{equation}
for all $m$ and $\alpha$. Consequently, $\Psi_{0, n}$ and
all of its first partial derivatives vanish at $\widetilde{T}_{0}$
since each term of these functions either contains $\tilde{t}^{\alpha}_{m}$
for some $\alpha$ and $m$, or contains a 1-point genus-$0$ 
correlation function. However, there is a little problem with this argument
since the genus-0 energy function is just a formal power series at the origin and it may not
converge at $\widetilde{T}_{0}$ (we would like to thank Getzler for pointing
out this to us). 
Although one might expect that such a nice function should converge,
rigorously speaking, we need to use the arguments in the last paragraph,
which are simply obtained by applying Lemma~\ref{lem:DilatonCorr} 
at another point.

In the rest of this paper, we will show that all second derivatives of 
$\Psi_{0, n}$ vanish by using the generalized WDVV equation as described
in the following strategy.
Write the first derivative
part of the operator $L_{n}$ as a vector field
${\cal L}_{n}$ on the big phase space. 
We already saw two of these vector fields in 
section~\ref{sec:specialVF}, i.e., ${\cal L}_{-1} = -{\cal S}$
and ${\cal L}_{0} = -{\cal X} - \frac{3-d}{2}{\cal D}$.
For any two operators $\tau_{k}({\cal O}_{\mu})$ and 
$\tau_{l}({\cal O}_{\nu})$, the generalized WDVV equation implies
\begin{eqnarray*}
	&& \sum_{\alpha}
	\left<\left< {\cal L}_{n} \left({\cal L}_{0} - (n+1){\cal D}\right)
		{\cal O}_{\alpha} \right>\right>_{0} 
	\left<\left< {\cal O}^{\alpha} \tau_{k}({\cal O}_{\mu})
		\tau_{l}({\cal O}_{\nu}) 
		\right>\right>_{0} \nonumber \\
	&=& \sum_{\alpha}
	\left<\left< {\cal L}_{n} \tau_{k}({\cal O}_{\mu})
		{\cal O}_{\alpha} \right>\right>_{0} 
	\left<\left< {\cal O}^{\alpha} 
		\left({\cal L}_{0} - (n+1){\cal D}\right)
		\tau_{l}({\cal O}_{\nu}) 
		\right>\right>_{0}.
\end{eqnarray*}
Compute
both sides of this equation by using the genus-$0$ $L_{n}$ constraint
(which is assumed to be true).
It can be shown that the difference of the 
resulting expressions is equal to
 $\frac{\partial^{2}}
{\partial t^{\nu}_{l}\partial t^{\mu}_{k}}\Psi_{0, n+1}$.
Therefore the generalized WDVV equation implies
that all second derivatives of $\Psi_{0, n+1}$ are zero. As noted above,
this proves the genus-$0$ $L_{n+1}$ constraint. 
Although the computation involved in this process is a little tedious,
it is in fact quite straightforward. The only subtleties here, if there is
any,  are
when and where to use the recursion 
formula~(\ref{eqn:StringRec}) and Lemma~\ref{lem:EulerCorr}.
In the rest of the paper, we carry out this strategy for the $L_{1}$
and $L_{2}$ constraints in full details. 
 Due to the existence of the Virasoro type
relations between $L_{n}$ operators, this implies all the genus-$0$
Virasoro constraints.

\section{$L_{1}$ constraint for genus zero free energy function}
\label{sec:g0L1}

As explained in Section \ref{sec:VirOp}, the genus-$0$ $L_{1}$ constraint
is equivalent to the equation $\Psi_{0, 1} = 0$, where
\begin{eqnarray}
\Psi_{0, 1}& = & \sum_{m, \alpha} (m+b_{\alpha})(m+b_{\alpha}+1) 
	\tilde{t}^{\alpha}_{m} \left<\left< \tau_{m+1}({\cal O}_{\alpha})
		\right>\right>_{0}  \nonumber \\
	&& + \sum_{m, \alpha, \beta} (2m+2b_{\alpha}+1) 
		{\cal C}_{\alpha}^{\beta}
	\tilde{t}^{\alpha}_{m} \left<\left< \tau_{m}({\cal O}_{\beta})
		\right>\right>_{0}  \nonumber \\
	&& + \sum_{m, \alpha, \beta} 
		({\cal C}^{2})_{\alpha}^{\beta}
	\tilde{t}^{\alpha}_{m} \left<\left< \tau_{m-1}({\cal O}_{\beta})
		\right>\right>_{0}  \nonumber \\
	&& + \frac{1}{2}\sum_{\alpha} b_{\alpha}(1-b_{\alpha}) 
	    \left<\left< {\cal O}_{\alpha}\right>\right>_{0} 
		\left<\left< {\cal O}^{\alpha}\right>\right>_{0} 
		 \nonumber \\
	&& + \frac{1}{2}\sum_{\alpha, \beta} 
		({\cal C}^{2})_{\alpha \beta} t^{\alpha}_{0} t^{\beta}_{0}.
\end{eqnarray}

As noted at the end of section~\ref{sec:VirOp}, 
to prove $\Psi_{0, 1} = 0$, it suffices to show that
all second partial derivatives of $\Psi_{0, 1}$ are equal to zero.
We will see that this fact actually follows from the 
generalized WDVV equation. According to the general strategy described
at the end of section~\ref{sec:VirOp}, we should compute 3-point correlation
functions involving two vector fields 
\[ 
	{\cal L}_{0}  =  -{\cal X} - \frac{3-d}{2} {\cal D} 
	  =   \sum_{m, \alpha} (m + b_{\alpha})\tilde{t}^{\alpha}_{m}
        \frac{\partial}{\partial t^{\alpha}_{m}} 
        + \sum_{m, \alpha, \beta} 
        {\cal C}_{\alpha}^{\beta}\tilde{t}^{\alpha}_{m}
        \frac{\partial}{\partial t^{\beta}_{m-1}}
\]
and 
\[ {\cal L}_{0} - {\cal D}  =  
	\sum_{m, \alpha} (m + b_{\alpha}+1)\tilde{t}^{\alpha}_{m}
        \frac{\partial}{\partial t^{\alpha}_{m}} 
        + \sum_{m, \alpha, \beta} 
        {\cal C}_{\alpha}^{\beta}\tilde{t}^{\alpha}_{m}
        \frac{\partial}{\partial t^{\beta}_{m-1}}.
\]
We first compute the following 3-point correlation function
\begin{lem} \label{lem:XXCorr}
\begin{eqnarray*}
	\left<\left< {\cal L}_{0} ({\cal L}_{0} - {\cal D})
		\tau_{m}({\cal O}_{\alpha}) \right>\right>_{0} 
	&=& - \sum_{n, \sigma} (n+b_{\sigma})(n+b_{\sigma}+1) 
	\tilde{t}^{\sigma}_{n} \left<\left< \tau_{n}({\cal O}_{\sigma})
		\tau_{m}({\cal O}_{\alpha})
		\right>\right>_{0}   \\
	&& - \sum_{n, \sigma, \rho} (2n+2b_{\sigma}+1) 
		{\cal C}_{\sigma}^{\rho}
	\tilde{t}^{\sigma}_{n} \left<\left< \tau_{n-1}({\cal O}_{\rho})
		\tau_{m}({\cal O}_{\alpha})
		\right>\right>_{0}   \\
	&& - \sum_{n, \sigma, \rho} 
		({\cal C}^{2})_{\sigma}^{\rho}
	\tilde{t}^{\sigma}_{n} \left<\left< \tau_{n-2}({\cal O}_{\rho})
		\tau_{m}({\cal O}_{\alpha})
		\right>\right>_{0}   \\
	&& + (m+ b_{\alpha})(m+b_{\alpha}-1) 
	    \left<\left< \tau_{m}({\cal O}_{\alpha})\right>\right>_{0} \\
	&& + \sum_{\sigma} (b_{\alpha} + b_{\sigma} + 2m -2)
		{\cal C}_{\alpha}^{\sigma}
		\left<\left< \tau_{m-1}({\cal O}_{\sigma})\right>\right>_{0} 
		 \\
	&& + \sum_{\sigma} ({\cal C}^{2})_{\alpha}^{\sigma}
 		\left<\left< \tau_{m-2}({\cal O}_{\sigma})\right>\right>_{0} 
		 \\
	&& + \delta_{m, 0} \left\{
		 \sum_{\sigma}(2b_{\alpha} -1) {\cal C}_{\alpha\sigma}
			t^{\sigma}_{0} 
		- \sum_{\sigma}({\cal C}^{2})_{\alpha \sigma} 
		\tilde{t}^{\sigma}_{1}
			\right\} \\
	&& +\delta_{m, 1} \sum_{\sigma}
		({\cal C}^{2})_{\alpha \sigma} t^{\sigma}_{0}.
\end{eqnarray*}
\end{lem}
{\bf Proof}:
By Lemma~\ref{lem:DilatonCorr} (3), 
\begin{eqnarray*}
&& \left<\left< {\cal L}_{0} 
                \tau_{n}({\cal O}_{\beta}) 
		\tau_{m}({\cal O}_{\alpha}) \right>\right>_{0} 
	=	-\left<\left< {\cal X} 
                \tau_{n}({\cal O}_{\beta}) 
		\tau_{m}({\cal O}_{\alpha}) \right>\right>_{0} \\
&=&   - \delta_{m, 0} \delta_{n,0} {\cal C}_{\alpha \beta} -
        (m+n+b_{\alpha}+b_{\beta}) \left<\left<
                \tau_{n}({\cal O}_{\beta}) 
	\tau_{m}({\cal O}_{\alpha}) \right>\right>_{0} \\
 &  &   - \sum_{\gamma} {\cal C}_{\alpha}^{\gamma}
        \left<\left<
                \tau_{n}({\cal O}_{\beta}) 
		\tau_{m-1}({\cal O}_{\gamma}) \right>\right>_{0} 
        - \sum_{\gamma} {\cal C}_{\beta}^{\gamma}
        \left<\left<
                \tau_{n-1}({\cal O}_{\gamma}) 
		\tau_{m}({\cal O}_{\alpha}) \right>\right>_{0}.
\end{eqnarray*} 
Hence
\begin{eqnarray} \label{eqn:L0L0:2}
&&	\left<\left< {\cal L}_{0} ({\cal L}_{0} - {\cal D})
		\tau_{m}({\cal O}_{\alpha}) \right>\right>_{0} 
	\nonumber \\
&=&  \sum_{n, \beta} (n + b_{\beta}+1)\tilde{t}^{\beta}_{n}
        \left<\left< {\cal L}_{0} 
                \tau_{n}({\cal O}_{\beta}) 
		\tau_{m}({\cal O}_{\alpha}) \right>\right>_{0} 
        + \sum_{n, \beta, \sigma} 
        {\cal C}_{\beta}^{\sigma}\tilde{t}^{\beta}_{n}
        \left<\left< {\cal L}_{0} 
                \tau_{n-1}({\cal O}_{\sigma}) 
		\tau_{m}({\cal O}_{\alpha}) \right>\right>_{0} 
	\nonumber \\
&=& - \sum_{n, \beta} (n + b_{\beta}+1)
	 (m+n+b_{\alpha}+b_{\beta})
	\tilde{t}^{\beta}_{n}
	 \left<\left<
	  \tau_{n}({\cal O}_{\beta})
	\tau_{m}({\cal O}_{\alpha})  \right>\right>_{0} 
	\nonumber \\
&& - \sum_{n, \beta, \gamma} (n + b_{\beta}+1)
	 {\cal C}_{\alpha}^{\gamma}
	\tilde{t}^{\beta}_{n}
	\left<\left<
                \tau_{n}({\cal O}_{\beta}) 
		\tau_{m-1}({\cal O}_{\gamma}) \right>\right>_{0} 
	\nonumber  \\
&&  - \sum_{n, \beta, \gamma} (n + b_{\beta}+1)
	{\cal C}_{\beta}^{\gamma}
	\tilde{t}^{\beta}_{n}      
	\left<\left<
                \tau_{n-1}({\cal O}_{\gamma})
		\tau_{m}({\cal O}_{\alpha})  \right>\right>_{0}  
	\nonumber \\
&& - \delta_{m, 0} \delta_{n,0} 
	\sum_{n, \beta} (n + b_{\beta}+1){\cal C}_{\alpha \beta}
	\tilde{t}^{\beta}_{n}     
	\nonumber  \\     
 &  &   - \sum_{n, \beta, \sigma} 
	  (m+n+b_{\alpha}+b_{\sigma}-1) 
        {\cal C}_{\beta}^{\sigma}
	\tilde{t}^{\beta}_{n}
	\left<\left< \tau_{n-1}({\cal O}_{\sigma}) 
	\tau_{m}({\cal O}_{\alpha})         
		\right>\right>_{0} 
	\nonumber \\
&&	- \sum_{n, \beta, \sigma, \gamma}
	    {\cal C}_{\beta}^{\sigma}
	{\cal C}_{\alpha}^{\gamma}
	\tilde{t}^{\beta}_{n}
        \left<\left<
                \tau_{n-1}({\cal O}_{\sigma})
		\tau_{m-1}({\cal O}_{\gamma})  \right>\right>_{0}
	\nonumber \\
&&	- \sum_{n, \beta, \sigma, \gamma}
	    {\cal C}_{\beta}^{\sigma}
	{\cal C}_{\sigma}^{\gamma}
      \tilde{t}^{\beta}_{n}
        \left<\left<
                \tau_{n-2}({\cal O}_{\gamma}) 
		\tau_{m}({\cal O}_{\alpha}) \right>\right>_{0} 
	\nonumber \\
& &  - \delta_{m, 0} \delta_{n,1} 
	\sum_{n, \beta, \sigma}  {\cal C}_{\beta}^{\sigma}
	{\cal C}_{\alpha \sigma} 
      \tilde{t}^{\beta}_{n}     
	\nonumber \\
&=& - \sum_{n, \beta} (n+b_{\beta})(n + b_{\beta}+1)
	\tilde{t}^{\beta}_{n}
	 \left<\left<
	  \tau_{n}({\cal O}_{\beta})
	\tau_{m}({\cal O}_{\alpha})  \right>\right>_{0} 
	\nonumber \\
&&  - \sum_{n, \beta, \gamma} (2n + 2b_{\beta}+1)
	{\cal C}_{\beta}^{\gamma}
	\tilde{t}^{\beta}_{n}      
	\left<\left<
                \tau_{n-1}({\cal O}_{\gamma})
		\tau_{m}({\cal O}_{\alpha})  \right>\right>_{0}  
	\nonumber \\
&&	- \sum_{n, \beta, \gamma}
	    ({\cal C}^{2})_{\beta}^{\gamma}
      \tilde{t}^{\beta}_{n}
        \left<\left<
                \tau_{n-2}({\cal O}_{\gamma}) 
		\tau_{m}({\cal O}_{\alpha}) \right>\right>_{0} 
	\nonumber \\
&& - \delta_{m, 0}  
	\sum_{\beta} (b_{\beta}+1){\cal C}_{\alpha \beta}
	\tilde{t}^{\beta}_{0}      
	\nonumber \\     
& &  - \delta_{m, 0}  
	\sum_{\beta}  ({\cal C}^{2})_{\alpha \beta}
      \tilde{t}^{\beta}_{1}    
	\nonumber  \\
&& - \sum_{n, \beta} (n + b_{\beta}+1)
	 (m+b_{\alpha})
	\tilde{t}^{\beta}_{n}
	 \left<\left<
	  \tau_{n}({\cal O}_{\beta})
	\tau_{m}({\cal O}_{\alpha})  \right>\right>_{0} 
	\nonumber \\
 &  &   - \sum_{n, \beta, \sigma} 
	  (m+b_{\alpha}) 
        {\cal C}_{\beta}^{\sigma}
	\tilde{t}^{\beta}_{n}
	\left<\left< \tau_{n-1}({\cal O}_{\sigma}) 
	\tau_{m}({\cal O}_{\alpha})         
		\right>\right>_{0}  
	\nonumber \\
&& - \sum_{n, \beta, \gamma} (n + b_{\beta}+1)
	 {\cal C}_{\alpha}^{\gamma}
	\tilde{t}^{\beta}_{n}
	\left<\left<
                \tau_{n}({\cal O}_{\beta}) 
		\tau_{m-1}({\cal O}_{\gamma}) \right>\right>_{0}  
	\nonumber \\
&&	- \sum_{n, \beta, \sigma, \gamma}
	    {\cal C}_{\beta}^{\sigma}
	{\cal C}_{\alpha}^{\gamma}
	\tilde{t}^{\beta}_{n}
        \left<\left<
                \tau_{n-1}({\cal O}_{\sigma})
		\tau_{m-1}({\cal O}_{\gamma})  \right>\right>_{0}. 
\end{eqnarray}
On the other hand, by Lemma~\ref{lem:DilatonCorr} (2) and
Lemma~\ref{lem:EulerCorr} (2), we have
\begin{eqnarray} \label{eqn:2ptEuler}
&&  \sum_{n, \beta} (n + b_{\beta}+1)
	\tilde{t}^{\beta}_{n}
	 \left<\left<
	  \tau_{n}({\cal O}_{\beta})
	\tau_{m}({\cal O}_{\alpha})  \right>\right>_{0} 
  + \sum_{n, \beta, \sigma} 
	        {\cal C}_{\beta}^{\sigma}
	\tilde{t}^{\beta}_{n}
	\left<\left< \tau_{n-1}({\cal O}_{\sigma}) 
	\tau_{m}({\cal O}_{\alpha})         
		\right>\right>_{0}  
	\nonumber \\
& = & \left<\left< ({\cal L}_{0} - {\cal D}) 
		\tau_{m}({\cal O}_{\alpha})         
		\right>\right>_{0}   
	\nonumber \\
& = & - \left<\left< {\cal X}
		\tau_{m}({\cal O}_{\alpha})         
		\right>\right>_{0} 
	- \frac{5-d}{2}
		\left<\left< {\cal D}
		\tau_{m}({\cal O}_{\alpha})         
		\right>\right>_{0}   
	\nonumber \\
& = &  - \left(m + b_{\alpha} - 1 \right)
                \left<\left< \tau_{m}({\cal O}_{\alpha}) \right>\right>_{0}
       - \sum_{\sigma} {\cal C}_{\alpha}^{\sigma}
        \left<\left< \tau_{m-1}({\cal O}_{\sigma}) \right>\right>_{0}
        - \delta_{m, 0} \sum_{\sigma} {\cal C}_{\alpha \sigma} t^{\sigma}_{0}.
\end{eqnarray}
The lemma then follows by applying 
 (\ref{eqn:2ptEuler}) to the last 4 terms in (\ref{eqn:L0L0:2}).
$\Box$

Setting $m = 0$ in Lemma~\ref{lem:XXCorr} , multiplying both sides of the
equation by $\left<\left< {\cal O}^{\alpha} \tau_{k}({\cal O}_{\mu})
		\tau_{l}({\cal O}_{\nu}) \right>\right>_{0}$, and
summing over $\alpha$, then applying the genus-$0$ topological recursion
relation, we get
\begin{eqnarray}
	&& \sum_{\alpha}\left<\left< {\cal L}_{0} ({\cal L}_{0} - {\cal D})
		{\cal O}_{\alpha} \right>\right>_{0} 
	\left<\left< {\cal O}^{\alpha} \tau_{k}({\cal O}_{\mu})
		\tau_{l}({\cal O}_{\nu}) 
		\right>\right>_{0} \nonumber \\
	&=& - \sum_{n, \sigma} (n+b_{\sigma})(n+b_{\sigma}+1) 
	\tilde{t}^{\sigma}_{n} \left<\left< \tau_{n+1}({\cal O}_{\sigma})
		\tau_{k}({\cal O}_{\mu})
		\tau_{l}({\cal O}_{\nu}) 
		\right>\right>_{0}   \nonumber  \\
	&& - \sum_{n, \sigma, \rho} (2n+2b_{\sigma}+1) 
		{\cal C}_{\sigma}^{\rho}
	\tilde{t}^{\sigma}_{n} \left<\left< \tau_{n}({\cal O}_{\rho})
		\tau_{k}({\cal O}_{\mu})
		\tau_{l}({\cal O}_{\nu}) 
		\right>\right>_{0}   \nonumber \\
	&& - \sum_{n, \sigma, \rho} 
		({\cal C}^{2})_{\sigma}^{\rho}
	\tilde{t}^{\sigma}_{n} \left<\left< \tau_{n-1}({\cal O}_{\rho})
		\tau_{k}({\cal O}_{\mu})
		\tau_{l}({\cal O}_{\nu}) 
		\right>\right>_{0}   \nonumber  \\
	&& + \sum_{\alpha} b_{\alpha}(b_{\alpha}-1) 
	    \left<\left< {\cal O}_{\alpha} \right>\right>_{0}
		\left<\left< {\cal O}^{\alpha}
		\tau_{k}({\cal O}_{\mu})
		\tau_{l}({\cal O}_{\nu})  \right>\right>_{0}. 
	\label{eqn:1WDVVleft}
\end{eqnarray}
Notice that the ranges of summations may change when using the topological
recursion relation. Hence some scattered terms may be absorbed into a
big summation after using the topological recursion relation.

On the other hand, using Lemma~\ref{lem:DilatonCorr} (3) and
Lemma~\ref{lem:EulerCorr} (3), we have
\begin{eqnarray}
	&& \sum_{\alpha, \beta}\left<\left< {\cal L}_{0} 
			\tau_{k}({\cal O}_{\mu})
		{\cal O}_{\alpha} \right>\right>_{0} 
		\eta^{\alpha \beta} 
	\left<\left< {\cal O}_{\beta} ({\cal L}_{0} - {\cal D})
		\tau_{l}({\cal O}_{\nu}) 
		\right>\right>_{0} \nonumber \\
	&=& \sum_{\alpha, \beta}\left<\left< {\cal X} \tau_{k}({\cal O}_{\mu})
		{\cal O}_{\alpha} \right>\right>_{0} 
		\eta^{\alpha \beta} 
	\left<\left< {\cal O}_{\beta} {\cal X}
		\tau_{l}({\cal O}_{\nu}) 
		\right>\right>_{0} \nonumber \\
	&=& \sum_{\alpha, \beta}\left\{ \delta_{k, 0} {\cal C}_{\mu \alpha}
		+ (k+b_{\mu}+b_{\alpha})
		 \left<\left< \tau_{k}({\cal O}_{\mu})
		{\cal O}_{\alpha} \right>\right>_{0} 
		+ \sum_{\sigma}{\cal C}_{\mu}^{\sigma}
		\left<\left< \tau_{k-1}({\cal O}_{\sigma})
		{\cal O}_{\alpha} \right>\right>_{0} \right\} 
		\eta^{\alpha \beta}
		\nonumber \\
	&& \textrm{ \hspace{30pt}}
		\left\{ \delta_{l, 0} {\cal C}_{\nu \beta}
		+ (l+b_{\nu}+b_{\beta})
		 \left<\left< \tau_{l}({\cal O}_{\nu})
		{\cal O}_{\beta} \right>\right>_{0} 
		+ \sum_{\rho}{\cal C}_{\nu}^{\rho}
		\left<\left< \tau_{l-1}({\cal O}_{\rho})
		{\cal O}_{\beta} \right>\right>_{0} \right\}. 
		\label{eqn:1WDVVright}
\end{eqnarray}
The generalized WDVV equation implies that the left hand sides
of equations (\ref{eqn:1WDVVleft}) and (\ref{eqn:1WDVVright})
are equal. However, the right hand sides of these two equations
appear very different from each other. One obvious distinction
between them is that the right hand side of (\ref{eqn:1WDVVright})
has only finitely many terms, while the right hand side of
(\ref{eqn:1WDVVleft}) has infinitely many terms due to the existence
of infinitely many gravitational descendents. In the rest of this
section, we will show that the difference of these two expressions
is $\frac{\partial}{\partial t^{\nu}_{l}}
	\frac{\partial}{\partial t^{\mu}_{k}} \Psi_{0, 1}$. 

We first prove two lemmas which express certain
quadratic functions of 2-point
correlation functions in terms of  linear functions of correlation functions.
\begin{lem}   \label{lem:1QuadForm1}
\begin{eqnarray*}
	&& \sum_{\alpha} \left\{b_{\alpha}(k+b_{\mu}-l-b_{\nu})-
		(k+b_{\mu})(l+b_{\nu}+1)\right\}
	\left<\left< \tau_{k}({\cal O}_{\mu})
		{\cal O}_{\alpha} \right>\right>_{0} 
	\left<\left< {\cal O}^{\alpha}
		\tau_{l}({\cal O}_{\nu}) 
		\right>\right>_{0} \\
	&=& \left(k+b_{\mu}-l-b_{\nu}\right) \sum_{\alpha}\left\{ 
		{\cal C}_{\nu}^{\alpha} 
		\left<\left< \tau_{k}({\cal O}_{\mu})
		\tau_{l}({\cal O}_{\alpha}) 
		\right>\right>_{0} 
		- {\cal C}_{\mu}^{\alpha} 
		\left<\left< \tau_{k}({\cal O}_{\alpha})
		\tau_{l}({\cal O}_{\nu}) 
		\right>\right>_{0}  \right\} \\
	&& -(k+b_{\mu})(k+b_{\mu}+1) \left<\left< \tau_{k+1}({\cal O}_{\mu})
		\tau_{l}({\cal O}_{\nu}) 
		\right>\right>_{0}                   \\
	&&   - (l+b_{\nu})(l+b_{\nu}+1) \left<\left< \tau_{k}({\cal O}_{\mu})
		\tau_{l+1}({\cal O}_{\nu}) 
		\right>\right>_{0}.   
\end{eqnarray*}
\end{lem}
{\bf Proof}: Let
\begin{eqnarray*}
f & := & \sum_{\alpha} \left\{b_{\alpha}(k+b_{\mu}-l-b_{\nu})-
		(k+b_{\mu})(l+b_{\nu}+1)\right\}
	\left<\left< \tau_{k}({\cal O}_{\mu})
		{\cal O}_{\alpha} \right>\right>_{0} 
	\left<\left< {\cal O}^{\alpha}
		\tau_{l}({\cal O}_{\nu}) 
		\right>\right>_{0} \\
& = & (k+b_{\mu}-l-b_{\nu}) 
	\sum_{\alpha} (b_{\alpha}- b_{\nu} -l -1)
	\left<\left< \tau_{k}({\cal O}_{\mu})
		{\cal O}_{\alpha} \right>\right>_{0} 
	\left<\left< {\cal O}^{\alpha}
		\tau_{l}({\cal O}_{\nu}) 
		\right>\right>_{0} \\
&& - (b_{\nu} +l)(b_{\nu} +l +1)
	\sum_{\alpha} 
	\left<\left< \tau_{k}({\cal O}_{\mu})
		{\cal O}_{\alpha} \right>\right>_{0} 
	\left<\left< {\cal O}^{\alpha}
		\tau_{l}({\cal O}_{\nu}) 
		\right>\right>_{0}. 
\end{eqnarray*}
Applying  Lemma~\ref{lem:EulerCorr} (3) to the first term and 
the recursion formula (\ref{eqn:StringRec}) to the second term, we have
\begin{eqnarray*}
f & = & (k+b_{\mu}-l-b_{\nu}) 
	\sum_{\alpha} 
	\left<\left< \tau_{k}({\cal O}_{\mu})
		{\cal O}_{\alpha} \right>\right>_{0}  \\
&& \textrm{ \hspace{60pt}}	\cdot \left\{  
	- \left<\left< {\cal O}^{\alpha} {\cal X}
		\tau_{l}({\cal O}_{\nu}) 
		\right>\right>_{0} 
	+ \sum_{\sigma} {\cal C}_{\nu}^{\sigma} 
		\left<\left< {\cal O}^{\alpha}
		\tau_{l-1}({\cal O}_{\sigma}) 
		\right>\right>_{0}  
 	+ \delta_{l, 0}
	 {\cal C}_{\nu}^{\alpha} \right\}
		\\
&& - (b_{\nu} +l)(b_{\nu} +l +1) \left\{
		\left<\left< \tau_{k+1}({\cal O}_{\mu})
				\tau_{l}({\cal O}_{\nu}) 
		\right>\right>_{0} 
		+ \left<\left< \tau_{k}({\cal O}_{\mu})
				\tau_{l+1}({\cal O}_{\nu}) 
		\right>\right>_{0} \right\}. \\
\end{eqnarray*}
Using genus-$0$ topological recursion relation to the first term
and formula (\ref{eqn:StringRec}) to the second term, we have
\begin{eqnarray*}
f & = & - (k+b_{\mu}-l-b_{\nu}) 
	\left<\left< \tau_{k+1}({\cal O}_{\mu})
		  {\cal X}
		\tau_{l}({\cal O}_{\nu}) 
		\right>\right>_{0}     \\
&& +  (k+b_{\mu}-l-b_{\nu}) 
	 \sum_{\sigma} {\cal C}_{\nu}^{\sigma} 
	\left\{
		\left<\left< \tau_{k+1}({\cal O}_{\mu})
				\tau_{l-1}({\cal O}_{\sigma}) 
		\right>\right>_{0} 
		+ \left<\left< \tau_{k}({\cal O}_{\mu})
				\tau_{l}({\cal O}_{\sigma}) 
		\right>\right>_{0}   \right. \\
&& \textrm{ \hspace{150pt}} \left.
		- \delta_{l, 0}\left<\left< \tau_{k}({\cal O}_{\mu})
				{\cal O}_{\sigma} 
				\right>\right>_{0} 
		\right\} \\
&&	+ \delta_{l, 0}(k+b_{\mu}-l-b_{\nu}) 
	\sum_{\alpha} 
	\left<\left< \tau_{k}({\cal O}_{\mu})
		{\cal O}_{\alpha} \right>\right>_{0}  
	 {\cal C}_{\nu}^{\alpha} 
		\\
&& - (b_{\nu} +l)(b_{\nu} +l +1) \left\{
		\left<\left< \tau_{k+1}({\cal O}_{\mu})
				\tau_{l}({\cal O}_{\nu}) 
		\right>\right>_{0} 
		+ \left<\left< \tau_{k}({\cal O}_{\mu})
				\tau_{l+1}({\cal O}_{\nu}) 
		\right>\right>_{0} \right\}. \\
\end{eqnarray*}
Applying Lemma~\ref{lem:EulerCorr} (3) to the first term, we obtain
\begin{eqnarray*}
f & = & - (k+b_{\mu}-l-b_{\nu}) 
	\left\{ (k+b_{\mu}+l+b_{\nu}+ 1)
		\left<\left< \tau_{k+1}({\cal O}_{\mu})
		\tau_{l}({\cal O}_{\nu}) 
		\right>\right>_{0}    \right. \\
&& \textrm{ \hspace{100pt}} \left.
	+ \sum_{\sigma} {\cal C}_{\mu}^{\sigma}
	\left<\left< \tau_{k}({\cal O}_{\sigma})
		\tau_{l}({\cal O}_{\nu}) 
		\right>\right>_{0}    
	+ \sum_{\sigma} {\cal C}_{\nu}^{\sigma}
	\left<\left< \tau_{k+1}({\cal O}_{\mu})
		\tau_{l-1}({\cal O}_{\sigma}) 
		\right>\right>_{0}    \right\}  \\
&& +  (k+b_{\mu}-l-b_{\nu}) 
	 \sum_{\sigma} {\cal C}_{\nu}^{\sigma} 
	\left\{
		\left<\left< \tau_{k+1}({\cal O}_{\mu})
				\tau_{l-1}({\cal O}_{\sigma}) 
		\right>\right>_{0} 
		+ \left<\left< \tau_{k}({\cal O}_{\mu})
				\tau_{l}({\cal O}_{\sigma}) 
		\right>\right>_{0}  
		\right\} \\
&& - (b_{\nu} +l)(b_{\nu} +l +1) \left\{
		\left<\left< \tau_{k+1}({\cal O}_{\mu})
				\tau_{l}({\cal O}_{\nu}) 
		\right>\right>_{0} 
		+ \left<\left< \tau_{k}({\cal O}_{\mu})
				\tau_{l+1}({\cal O}_{\nu}) 
		\right>\right>_{0} \right\}. \\
\end{eqnarray*}
Simplifying this expression, we obtain the desired formula.  
$\Box$

\begin{lem}  \label{lem:1QuadForm2}
\begin{eqnarray*}
&& \sum_{\alpha, \beta}(k+b_{\alpha}+b_{\mu}) {\cal C}_{\nu}^{\beta}
	\left<\left< \tau_{k}({\cal O}_{\mu})
		{\cal O}_{\alpha} \right>\right>_{0} 
	\left<\left< {\cal O}^{\alpha}
		\tau_{l-1}({\cal O}_{\beta}) 
		\right>\right>_{0}   \\
	&& + \sum_{\alpha, \beta, \sigma}{\cal C}_{\mu}^{\alpha}
	{\cal C}_{\nu}^{\beta}
	\left<\left< \tau_{k-1}({\cal O}_{\alpha})
		{\cal O}_{\sigma} \right>\right>_{0} 
	\left<\left< {\cal O}^{\sigma}
		\tau_{l-1}({\cal O}_{\beta}) 
		\right>\right>_{0} \\
&=& \sum_{\alpha} (k+b_{\mu} + l + b_{\nu}+1) 
	{\cal C}_{\nu}^{\alpha}
	\left<\left< \tau_{k}({\cal O}_{\mu})
		\tau_{l}({\cal O}_{\alpha}) 
		\right>\right>_{0}   \\
	&&+\sum_{\alpha, \beta}{\cal C}_{\mu}^{\alpha}
	{\cal C}_{\nu}^{\beta}
	\left<\left< \tau_{k-1}({\cal O}_{\alpha})
		\tau_{l}({\cal O}_{\beta}) 
		\right>\right>_{0} 
	+\sum_{\alpha}({\cal C}^{2})_{\nu}^{\alpha}
	\left<\left< \tau_{k}({\cal O}_{\mu})
		\tau_{l-1}({\cal O}_{\alpha}) 
		\right>\right>_{0} \\
	&&- \delta_{k, 0} \sum_{\alpha, \beta}{\cal C}_{\mu}^{\alpha}
	{\cal C}_{\nu}^{\beta}
	\left<\left< {\cal O}_{\alpha}
		\tau_{l-1}({\cal O}_{\beta}) 
		\right>\right>_{0} 
	- \delta_{l, 0} \sum_{\alpha}{\cal C}_{\nu}^{\alpha}
	\left<\left< {\cal X} \tau_{k}({\cal O}_{\mu})
		{\cal O}_{\alpha}) 
		\right>\right>_{0} \\
	&& + \delta_{k, 0} \delta_{l, 0} ({\cal C}^{2})_{\mu \nu}.
\end{eqnarray*}
\end{lem}
{\bf Proof}: Let
\begin{eqnarray*}
f &:= & \sum_{\alpha, \beta}(k+b_{\alpha}+b_{\mu}) {\cal C}_{\nu}^{\beta}
	\left<\left< \tau_{k}({\cal O}_{\mu})
		{\cal O}_{\alpha} \right>\right>_{0} 
	\left<\left< {\cal O}^{\alpha}
		\tau_{l-1}({\cal O}_{\beta}) 
		\right>\right>_{0}   \\
	&& + \sum_{\alpha, \beta, \sigma}{\cal C}_{\mu}^{\alpha}
	{\cal C}_{\nu}^{\beta}
	\left<\left< \tau_{k-1}({\cal O}_{\alpha})
		{\cal O}_{\sigma} \right>\right>_{0} 
	\left<\left< {\cal O}^{\sigma}
		\tau_{l-1}({\cal O}_{\beta}) 
		\right>\right>_{0} \\
 & = & \sum_{\alpha, \beta} {\cal C}_{\nu}^{\beta}
	\left\{
		(k+b_{\alpha}+b_{\mu}) 
	\left<\left< \tau_{k}({\cal O}_{\mu})
		{\cal O}_{\alpha} \right>\right>_{0} 
	 + \sum_{\sigma}{\cal C}_{\mu}^{\sigma}
	\left<\left< \tau_{k-1}({\cal O}_{\sigma})
		{\cal O}_{\alpha} \right>\right>_{0} 
			\right\}
	\left<\left< {\cal O}^{\alpha}
		\tau_{l-1}({\cal O}_{\beta}) 
		\right>\right>_{0}.  
\end{eqnarray*}
Using Lemma~\ref{lem:EulerCorr} (3), we have
\begin{eqnarray*}
f &=& \sum_{\alpha, \beta} {\cal C}_{\nu}^{\beta}
	\left\{
	\left<\left< {\cal X}\tau_{k}({\cal O}_{\mu})
		{\cal O}_{\alpha} \right>\right>_{0} 
		- \delta_{k, 0} {\cal C}_{\mu \alpha}
			\right\}
	\left<\left< {\cal O}^{\alpha}
		\tau_{l-1}({\cal O}_{\beta}) 
		\right>\right>_{0}.   
\end{eqnarray*}
Using topological recursion relation to the first term, we have
\begin{eqnarray*}
f &=& \sum_{\beta} {\cal C}_{\nu}^{\beta}
	\left\{ \left<\left< {\cal X}\tau_{k}({\cal O}_{\mu})
		\tau_{l}({\cal O}_{\beta}) \right>\right>_{0} 
	  - \delta_{l, 0}\left<\left< {\cal X}\tau_{k}({\cal O}_{\mu})
		{\cal O}_{\beta} \right>\right>_{0} \right\} \\
&& 	- \delta_{k, 0} \sum_{\alpha, \beta} {\cal C}_{\nu}^{\beta}
		{\cal C}_{\mu \alpha}
	\left<\left< {\cal O}^{\alpha}
		\tau_{l-1}({\cal O}_{\beta}) 
		\right>\right>_{0}.   
\end{eqnarray*}
Using Lemma~\ref{lem:EulerCorr} (3) again to the first term, we have
\begin{eqnarray*}
f &=& \sum_{\beta} {\cal C}_{\nu}^{\beta}
	\left\{ \delta_{k, 0} \delta_{l, 0} {\cal C}_{\mu \beta} +
		(k+b_{\mu}+l+b_{\beta})
		\left<\left< \tau_{k}({\cal O}_{\mu})
		\tau_{l}({\cal O}_{\beta}) \right>\right>_{0} \right. \\
&& \textrm{ \hspace{50 pt}}
	 \left. + \sum_{\sigma} {\cal C}_{\mu}^{\sigma}
		\left<\left< \tau_{k-1}({\cal O}_{\sigma})
		\tau_{l}({\cal O}_{\beta}) \right>\right>_{0}  
		+ \sum_{\sigma} {\cal C}_{\beta}^{\sigma}
	\left<\left< \tau_{k}({\cal O}_{\mu})
		\tau_{l-1}({\cal O}_{\sigma}) \right>\right>_{0} \right\} \\
&&  - \delta_{l, 0}\sum_{\beta} {\cal C}_{\nu}^{\beta}
	\left<\left< {\cal X}\tau_{k}({\cal O}_{\mu})
		{\cal O}_{\beta} \right>\right>_{0} \\
&&	- \delta_{k, 0} \sum_{\alpha, \beta} {\cal C}_{\nu}^{\beta}
		{\cal C}_{\mu \alpha}
	\left<\left< {\cal O}^{\alpha}
		\tau_{l-1}({\cal O}_{\beta}) 
		\right>\right>_{0}.   
\end{eqnarray*}
The lemma then follows from the fact that 
${\cal C}_{\nu}^{\beta} \neq 0$ implies $b_{\beta} = b_{\nu} + 1$.
$\Box$

Now we can deduce from (\ref{eqn:1WDVVright}) the following
\begin{lem}  \label{lem:1WDVVright}
\begin{eqnarray*}
	&& \sum_{\alpha}\left<\left< {\cal L}_{0} \tau_{k}({\cal O}_{\mu})
		{\cal O}_{\alpha} \right>\right>_{0} 
	\left<\left< {\cal O}^{\alpha} ({\cal L}_{0} - {\cal D})
		\tau_{l}({\cal O}_{\nu}) 
		\right>\right>_{0}  \\
&=&  (k+b_{\mu})(k+b_{\mu}+1) \left<\left< \tau_{k+1}({\cal O}_{\mu})
		\tau_{l}({\cal O}_{\nu}) 
		\right>\right>_{0}                   \\
&& +(l+b_{\nu})(l+b_{\nu}+1) \left<\left< \tau_{k}({\cal O}_{\mu})
		\tau_{l+1}({\cal O}_{\nu}) 
		\right>\right>_{0}                   \\
&& + \sum_{\alpha}(2k+2b_{\mu}+1) {\cal C}_{\mu}^{\alpha}
	 \left<\left< \tau_{k}({\cal O}_{\alpha})
		\tau_{l}({\cal O}_{\nu}) 
		\right>\right>_{0}                   \\
&& + \sum_{\alpha}(2l+2b_{\nu}+1) {\cal C}_{\nu}^{\alpha}
	 \left<\left< \tau_{k}({\cal O}_{\mu})
		\tau_{l}({\cal O}_{\alpha}) 
		\right>\right>_{0}                   \\
&& + \sum_{\alpha}({\cal C}^{2})_{\mu}^{\alpha}
	 \left<\left< \tau_{k-1}({\cal O}_{\alpha})
		\tau_{l}({\cal O}_{\nu}) 
		\right>\right>_{0}                   \\
&& + \sum_{\alpha} ({\cal C}^{2})_{\nu}^{\alpha}
	 \left<\left< \tau_{k}({\cal O}_{\mu})
		\tau_{l-1}({\cal O}_{\alpha}) 
		\right>\right>_{0}                   \\
&& + \sum_{\alpha} b_{\alpha}(1- b_{\alpha}) 
	\left<\left< \tau_{k}({\cal O}_{\mu})
		{\cal O}_{\alpha}
		\right>\right>_{0}  
	\left<\left< {\cal O}^{\alpha}
		\tau_{l}({\cal O}_{\nu}) 
		\right>\right>_{0}   \\
&& + \delta_{k, 0} \delta_{l, 0} ({\cal C}^{2})_{\mu \nu}.
\end{eqnarray*}
\end{lem}
{\bf Proof}:
By (\ref{eqn:1WDVVright}),	
\begin{eqnarray*}
	&&  \sum_{\alpha, \beta}\left<\left< {\cal L}_{0} 
			\tau_{k}({\cal O}_{\mu})
		{\cal O}_{\alpha} \right>\right>_{0} 
		\eta^{\alpha \beta} 
	\left<\left< {\cal O}_{\beta} ({\cal L}_{0} - {\cal D})
		\tau_{l}({\cal O}_{\nu}) 
		\right>\right>_{0} \nonumber \\
	&=& -\sum_{\alpha} \left\{b_{\alpha}(k+b_{\mu}-l-b_{\nu})-
                (k+b_{\mu})(l+b_{\nu}+1)\right\}
        \left<\left< \tau_{k}({\cal O}_{\mu})
                {\cal O}_{\alpha} \right>\right>_{0} 
        \left<\left< {\cal O}^{\alpha}
                \tau_{l}({\cal O}_{\nu}) 
                \right>\right>_{0}   \nonumber \\
	&& + \sum_{\alpha} b_{\alpha}(1- b_{\alpha}) 
	\left<\left< \tau_{k}({\cal O}_{\mu})
		{\cal O}_{\alpha}
		\right>\right>_{0}  
	\left<\left< {\cal O}^{\alpha}
		\tau_{l}({\cal O}_{\nu}) 
		\right>\right>_{0}   \nonumber \\
	&& +  \sum_{\alpha, \rho}
		(k+b_{\mu}+b_{\alpha}){\cal C}_{\nu}^{\rho}
		 \left<\left< \tau_{k}({\cal O}_{\mu})
		{\cal O}_{\alpha} \right>\right>_{0} 
		\left<\left< 
			{\cal O}^{\alpha} 
			\tau_{l-1}({\cal O}_{\rho})
		\right>\right>_{0}  \nonumber \\
	&& + \sum_{\alpha, \sigma, \rho}{\cal C}_{\mu}^{\sigma}
		{\cal C}_{\nu}^{\rho}
		\left<\left< \tau_{k-1}({\cal O}_{\sigma})
		{\cal O}_{\alpha} \right>\right>_{0}
		\left<\left< {\cal O}^{\alpha} \tau_{l-1}({\cal O}_{\rho})
		 \right>\right>_{0}  \nonumber \\
	&& + \sum_{\beta, \sigma}(l+b_{\nu}+b_{\beta}) {\cal C}_{\mu}^{\sigma}
		\left<\left< \tau_{k-1}({\cal O}_{\sigma})
		{\cal O}^{\beta} \right>\right>_{0} 
		\left<\left< {\cal O}_{\beta} \tau_{l}({\cal O}_{\nu})
		\right>\right>_{0} \nonumber \\
	&& + \sum_{\beta, \sigma, \rho}{\cal C}_{\mu}^{\sigma}
		{\cal C}_{\nu}^{\rho}
		\left<\left< \tau_{k-1}({\cal O}_{\sigma})
		{\cal O}^{\beta} \right>\right>_{0} 
		\left<\left< {\cal O}_{\beta} \tau_{l-1}({\cal O}_{\rho})
		 \right>\right>_{0}  \nonumber \\
	&& -  \sum_{\beta, \sigma, \rho}{\cal C}_{\mu}^{\sigma}
		{\cal C}_{\nu}^{\rho}
		\left<\left< \tau_{k-1}({\cal O}_{\sigma})
		{\cal O}^{\beta} \right>\right>_{0} 
		\left<\left< {\cal O}_{\beta} \tau_{l-1}({\cal O}_{\rho})
		 \right>\right>_{0}  \nonumber \\
	&& + \delta_{k, 0} \sum_{\beta}{\cal C}_{\mu}^{\beta}
		\left\{ (l+b_{\nu}+b_{\beta})
		 \left<\left< \tau_{l}({\cal O}_{\nu})
		{\cal O}_{\beta} \right>\right>_{0} 
		+ \sum_{\rho}{\cal C}_{\nu}^{\rho}
		\left<\left< \tau_{l-1}({\cal O}_{\rho})
		{\cal O}_{\beta} \right>\right>_{0} \right\} 
		\nonumber \\
	&& + \delta_{l, 0} \sum_{\alpha}{\cal C}_{\nu}^{\alpha}
		\left\{ (k+b_{\mu}+b_{\alpha})
		 \left<\left< \tau_{k}({\cal O}_{\mu})
		{\cal O}_{\alpha} \right>\right>_{0} 
		+ \sum_{\sigma}{\cal C}_{\mu}^{\sigma}
		\left<\left< \tau_{k-1}({\cal O}_{\sigma})
		{\cal O}_{\alpha} \right>\right>_{0} \right\} 
		\nonumber \\
	&& + \delta_{k, 0}  \delta_{l, 0} 
		({\cal C}^{2})_{\mu \nu}.
\end{eqnarray*}
Applying Lemma~\ref{lem:1QuadForm1} to the first term,
Lemma~\ref{lem:1QuadForm2} to the third and fourth terms,
and an analogue of Lemma~\ref{lem:1QuadForm2} with
$(\mu, k)$ interchanged with $(\nu, l)$ to the fifth and sixth terms,
and formula (\ref{eqn:StringRec}) to the seventh term, we obtain
\begin{eqnarray*}
	&&  \sum_{\alpha, \beta}\left<\left< {\cal L}_{0} 
			\tau_{k}({\cal O}_{\mu})
		{\cal O}_{\alpha} \right>\right>_{0} 
		\eta^{\alpha \beta} 
	\left<\left< {\cal O}_{\beta} ({\cal L}_{0} - {\cal D})
		\tau_{l}({\cal O}_{\nu}) 
		\right>\right>_{0} \nonumber \\
	&=& -\left(k+b_{\mu}-l-b_{\nu}\right) \sum_{\alpha}\left\{ 
                {\cal C}_{\nu}^{\alpha} 
                \left<\left< \tau_{k}({\cal O}_{\mu})
                \tau_{l}({\cal O}_{\alpha}) 
                \right>\right>_{0} 
                - {\cal C}_{\mu}^{\alpha} 
                \left<\left< \tau_{k}({\cal O}_{\alpha})
                \tau_{l}({\cal O}_{\nu}) 
                \right>\right>_{0}  \right\} \\
        && +(k+b_{\mu})(k+b_{\mu}+1) \left<\left< \tau_{k+1}({\cal O}_{\mu})
                \tau_{l}({\cal O}_{\nu}) 
                \right>\right>_{0}                   \\
        && + (l+b_{\nu})(l+b_{\nu}+1) \left<\left< \tau_{k}({\cal O}_{\mu})
                \tau_{l+1}({\cal O}_{\nu}) 
                \right>\right>_{0}    \\
	&& + \sum_{\alpha} b_{\alpha}(1- b_{\alpha}) 
	\left<\left< \tau_{k}({\cal O}_{\mu})
		{\cal O}_{\alpha}
		\right>\right>_{0}  
	\left<\left< {\cal O}^{\alpha}
		\tau_{l}({\cal O}_{\nu}) 
		\right>\right>_{0}   \nonumber \\
	&& +\sum_{\alpha} (k+b_{\mu} + l + b_{\nu}+1) 
        {\cal C}_{\nu}^{\alpha}
        \left<\left< \tau_{k}({\cal O}_{\mu})
                \tau_{l}({\cal O}_{\alpha}) 
                \right>\right>_{0}   \\
        &&+\sum_{\alpha, \beta}{\cal C}_{\mu}^{\alpha}
        {\cal C}_{\nu}^{\beta}
        \left<\left< \tau_{k-1}({\cal O}_{\alpha})
                \tau_{l}({\cal O}_{\beta}) 
                \right>\right>_{0} 
        +\sum_{\alpha}({\cal C}^{2})_{\nu}^{\alpha}
        \left<\left< \tau_{k}({\cal O}_{\mu})
                \tau_{l-1}({\cal O}_{\alpha}) 
                \right>\right>_{0} \\
	&&- \delta_{k, 0} \sum_{\alpha, \beta}{\cal C}_{\mu}^{\alpha}
        {\cal C}_{\nu}^{\beta}
        \left<\left< {\cal O}_{\alpha}
                \tau_{l-1}({\cal O}_{\beta}) 
                \right>\right>_{0} 
        - \delta_{l, 0} \sum_{\alpha}{\cal C}_{\nu}^{\alpha}
        \left<\left< {\cal X} \tau_{k}({\cal O}_{\mu})
                {\cal O}_{\alpha}) 
                \right>\right>_{0} \\
        && + \delta_{k, 0} \delta_{l, 0} ({\cal C}^{2})_{\mu \nu} \\
	&& + \sum_{\alpha} (k+b_{\mu} + l + b_{\nu}+1) 
        {\cal C}_{\mu}^{\alpha}
        \left<\left< \tau_{l}({\cal O}_{\nu})
                \tau_{k}({\cal O}_{\alpha}) 
                \right>\right>_{0}   \\
        &&+\sum_{\alpha, \beta}{\cal C}_{\nu}^{\alpha}
        {\cal C}_{\mu}^{\beta}
        \left<\left< \tau_{l-1}({\cal O}_{\alpha})
                \tau_{k}({\cal O}_{\beta}) 
                \right>\right>_{0} 
        +\sum_{\alpha}({\cal C}^{2})_{\mu}^{\alpha}
        \left<\left< \tau_{l}({\cal O}_{\nu})
                \tau_{k-1}({\cal O}_{\alpha}) 
                \right>\right>_{0} \\
	&&- \delta_{l, 0} \sum_{\alpha, \beta}{\cal C}_{\nu}^{\alpha}
        {\cal C}_{\mu}^{\beta}
        \left<\left< {\cal O}_{\alpha}
                \tau_{k-1}({\cal O}_{\beta}) 
                \right>\right>_{0} 
        - \delta_{k, 0} \sum_{\alpha}{\cal C}_{\mu}^{\alpha}
        \left<\left< {\cal X} \tau_{l}({\cal O}_{\nu})
                {\cal O}_{\alpha}) 
                \right>\right>_{0} \\
        && + \delta_{k, 0} \delta_{l, 0} ({\cal C}^{2})_{\mu \nu} \\
	&& -  \sum_{\sigma, \rho}{\cal C}_{\mu}^{\sigma}
		{\cal C}_{\nu}^{\rho}
		\left\{ 
		\left<\left< \tau_{k}({\cal O}_{\sigma})
		 \tau_{l-1}({\cal O}_{\rho})
		 \right>\right>_{0} +
		\left<\left< \tau_{k-1}({\cal O}_{\sigma})
		 \tau_{l}({\cal O}_{\rho})
		 \right>\right>_{0}   \right.  \\
	&& \textrm{ \hspace{100pt}} \left.
		- \delta_{k, 0}\left<\left< {\cal O}_{\sigma}
		 \tau_{l-1}({\cal O}_{\rho})
		 \right>\right>_{0} 
		-\delta_{l, 0} \left<\left< \tau_{k-1}({\cal O}_{\sigma})
		 {\cal O}_{\rho}
		 \right>\right>_{0}  
			\right\}
		\nonumber \\
	&& + \delta_{k, 0} \sum_{\beta}{\cal C}_{\mu}^{\beta}
		\left\{ (l+b_{\nu}+b_{\beta})
		 \left<\left< \tau_{l}({\cal O}_{\nu})
		{\cal O}_{\beta} \right>\right>_{0} 
		+ \sum_{\rho}{\cal C}_{\nu}^{\rho}
		\left<\left< \tau_{l-1}({\cal O}_{\rho})
		{\cal O}_{\beta} \right>\right>_{0} \right\} 
		\nonumber \\
	&& + \delta_{l, 0} \sum_{\alpha}{\cal C}_{\nu}^{\alpha}
		\left\{ (k+b_{\mu}+b_{\alpha})
		 \left<\left< \tau_{k}({\cal O}_{\mu})
		{\cal O}_{\alpha} \right>\right>_{0} 
		+ \sum_{\sigma}{\cal C}_{\mu}^{\sigma}
		\left<\left< \tau_{k-1}({\cal O}_{\sigma})
		{\cal O}_{\alpha} \right>\right>_{0} \right\} 
		\nonumber \\
	&& + \delta_{k, 0}  \delta_{l, 0} 
		({\cal C}^{2})_{\mu \nu}.
\end{eqnarray*}
Simplifying this expression and applying Lemma~\ref{lem:EulerCorr} (3)
to the two terms containing the Euler vector field ${\cal X}$, we obtain
the desired formula.
$\Box$

Now it's straightforward to check that the difference
between the right hand side of Lemma~\ref{lem:1WDVVright} and
the right hand side of equation (\ref{eqn:1WDVVleft})
is $\frac{\partial}{\partial t^{\nu}_{l}}
	\frac{\partial}{\partial t^{\mu}_{k}} \Psi_{0, 1}$. 
Hence the generalized WDVV equation implies that 
$\frac{\partial}{\partial t^{\nu}_{l}}
	\frac{\partial}{\partial t^{\mu}_{k}} \Psi_{0, 1} = 0$. 
Therefore we proved that all the second derivatives of
$\Psi_{0, 1}$ vanish. As mentioned at the end of section~\ref{sec:VirOp},
this implies the following
\begin{pro} \label{pro:g0L1}
  The genus-$0$ free energy function $F_{0}$ 
	satisfies the $L_{1}$ constraint.
\end{pro}
The proof of other genus-$0$ Virasoro constraints has the similar
flavor, as we will see in section \ref{sec:g0L2} for the case of
$L_{2}$ constraint.

\section{$L_{2}$ constraint for genus zero free energy function}
\label{sec:g0L2}

The genus-$0$ $L_{2}$ constraint
is equivalent to the equation $\Psi_{0, 2} = 0$, where
\begin{eqnarray}
\Psi_{0, 2}& = & \sum_{m, \alpha} (m+b_{\alpha})(m+b_{\alpha}+1) 
		(m+b_{\alpha}+2) 
        \tilde{t}^{\alpha}_{m} \left<\left< \tau_{m+2}({\cal O}_{\alpha})
                \right>\right>_{0}  \nonumber \\
        && + \sum_{m, \alpha, \beta} 
	\left\{ 3(m+b_{\alpha})^{2} + 6(m+b_{\alpha}) +2 \right\} 
     {\cal C}_{\alpha}^{\beta}
        \tilde{t}^{\alpha}_{m} \left<\left< \tau_{m + 1}({\cal O}_{\beta})
                \right>\right>_{0}  \nonumber \\
        && + \sum_{m, \alpha, \beta} 3(m+b_{\alpha}+1) 
                ({\cal C}^{2})_{\alpha}^{\beta}
        \tilde{t}^{\alpha}_{m} \left<\left< \tau_{m}({\cal O}_{\beta})
                \right>\right>_{0}  \nonumber \\
	&& + \sum_{m, \alpha, \beta} 
                ({\cal C}^{3})_{\alpha}^{\beta}
        \tilde{t}^{\alpha}_{m} \left<\left< \tau_{m-1}({\cal O}_{\beta})
                \right>\right>_{0}  \nonumber \\
        && - \sum_{\alpha} (b_{\alpha}-1) b_{\alpha}(b_{\alpha}+1) 
            \left<\left< \tau_{1}({\cal O}_{\alpha})\right>\right>_{0} 
                \left<\left< {\cal O}^{\alpha}\right>\right>_{0} 
                 \nonumber \\
	&& - \frac{1}{2}\sum_{\alpha, \beta} (3b_{\alpha}^{2} -1)
		{\cal C}_{\alpha}^{\beta} 
            \left<\left< {\cal O}_{\beta}\right>\right>_{0} 
                \left<\left< {\cal O}^{\alpha}\right>\right>_{0} 
                 \nonumber \\
        && + \frac{1}{2}\sum_{\alpha, \beta} 
                ({\cal C}^{3})_{\alpha \beta} t^{\alpha}_{0} t^{\beta}_{0}.
\end{eqnarray}
As in the proof of the $L_{1}$ constraint, we only need to show that
all second derivatives of $\Psi_{0, 2}$ are equal to zero. This time
we need to compute 3-point correlation functions involving vector fields
\begin{eqnarray}
{\cal L}_{1} & := & \sum_{m, \alpha} (m+b_{\alpha})(m+b_{\alpha}+1) 
        \tilde{t}^{\alpha}_{m} 
		\frac{\partial}{\partial t^{\alpha}_{m+1}}  
		\nonumber \\
        && + \sum_{m, \alpha , \beta} (2m+2b_{\alpha}+1) 
		{\cal C}_{\alpha}^{\beta}
        \tilde{t}^{\alpha}_{m} 
		\frac{\partial}{\partial t^{\beta}_{m}} 
	  \nonumber \\
        && + \sum_{m, \alpha, \beta} 
                ({\cal C}^{2})_{\alpha}^{\beta}
        \tilde{t}^{\alpha}_{m} 
	\frac{\partial}{\partial t^{\beta}_{m-1}}
\end{eqnarray}
and
\[ {\cal L}_{0} - 2{\cal D}
  =  \sum_{m, \alpha} (m + b_{\alpha} +2)\tilde{t}^{\alpha}_{m}
        \frac{\partial}{\partial t^{\alpha}_{m}} 
        + \sum_{m, \alpha, \beta} 
        {\cal C}_{\alpha}^{\beta}\tilde{t}^{\alpha}_{m}
        \frac{\partial}{\partial t^{\beta}_{m-1}}.
\]

The genus-$0$ $L_{1}$ constraint, can be reformulated as
\begin{equation} \label{eqn:1ptL1}
	\left<\left< {\cal L}_{1} \right>\right>_{0}
	= - \frac{1}{2}\sum_{\alpha} b_{\alpha}(1-b_{\alpha}) 
            \left<\left< {\cal O}_{\alpha}\right>\right>_{0} 
                \left<\left< {\cal O}^{\alpha}\right>\right>_{0} 
                 \nonumber \\
          - \frac{1}{2}\sum_{\alpha, \beta} 
                ({\cal C}^{2})_{\alpha \beta} t^{\alpha}_{0} t^{\beta}_{0}.
\end{equation}
Using this equation and the fact that
\begin{eqnarray} 
	 [{\cal L}_{1}, \, \frac{\partial}{\partial t^{\alpha}_{m}} ] 
      &  = & - (m+b_{\alpha})(m+b_{\alpha}+1)  
		\frac{\partial}{\partial t^{\alpha}_{m+1}}  
	\nonumber \\
	&& - \sum_{\beta}
		(2m + 2b_{\alpha} +1) {\cal C}_{\alpha}^{\beta}
                \frac{\partial}{\partial t^{\beta}_{m}} 
		\nonumber \\
       &&  - \sum_{\beta} ({\cal C}^{2})_{\alpha}^{\beta}
                \frac{\partial}{\partial t^{\beta}_{m-1}},
	\label{eqn:L1Bracket}
\end{eqnarray}
we can prove the following
\begin{lem} \label{lem:L1Corr}
\begin{eqnarray*}
\left<\left< {\cal L}_{1} \tau_{m}({\cal O}_{\alpha}) 
                \tau_{n}({\cal O}_{\beta}) \right>\right>_{0}
& = & - (m+b_{\alpha})(m+b_{\alpha}+1) 
    \left<\left<\tau_{m+1}({\cal O}_{\alpha}) 
                \tau_{n}({\cal O}_{\beta}) \right>\right>_{0} 
                \\
	&& - \sum_{\sigma} (2m+2b_{\alpha}+1)
	{\cal C}_{\alpha}^{\sigma}
            \left<\left<\tau_{m}({\cal O}_{\sigma}) 
                \tau_{n}({\cal O}_{\beta}) \right>\right>_{0} 
           \\
        && - \sum_{\sigma} 
                ({\cal C}^{2})_{\alpha}^{\sigma}
            \left<\left<\tau_{m-1}({\cal O}_{\sigma}) 
                \tau_{n}({\cal O}_{\beta}) \right>\right>_{0} \\
	&& - (n+b_{\beta})(n+b_{\beta}+1)
	    \left<\left<\tau_{m}({\cal O}_{\alpha}) 
                \tau_{n+1}({\cal O}_{\beta}) \right>\right>_{0} \\
	&& - \sum_{\sigma} (2n+2b_{\beta}+1)
	{\cal C}_{\beta}^{\sigma}
            \left<\left<\tau_{m}({\cal O}_{\alpha}) 
                \tau_{n}({\cal O}_{\sigma}) \right>\right>_{0} 
           \\
        && - \sum_{\sigma} 
                ({\cal C}^{2})_{\beta}^{\sigma}
            \left<\left<\tau_{m}({\cal O}_{\alpha}) 
                \tau_{n-1}({\cal O}_{\sigma}) \right>\right>_{0} \\
	&& - \sum_{\sigma} b_{\sigma}(1-b_{\sigma})
		 \left<\left< \tau_{m}({\cal O}_{\alpha})
                \tau_{n}({\cal O}_{\beta}) 
                {\cal O}_{\sigma}\right>\right>_{0} 
                \left<\left< {\cal O}^{\sigma}\right>\right>_{0} 
                 \\
	&& - \sum_{\sigma} b_{\sigma}(1-b_{\sigma})
		 \left<\left< \tau_{m}({\cal O}_{\alpha})
                {\cal O}_{\sigma}\right>\right>_{0} 
                \left<\left< {\cal O}^{\sigma}
			\tau_{n}({\cal O}_{\beta}) \right>\right>_{0} 
                 \\
        &&  - \delta_{m, 0} \delta_{n, 0} 
                ({\cal C}^{2})_{\alpha \beta}.
\end{eqnarray*}
\end{lem}
{\bf Proof}:
\begin{eqnarray*}
\left<\left< {\cal L}_{1} \tau_{m}({\cal O}_{\alpha}) 
                \tau_{n}({\cal O}_{\beta}) \right>\right>_{0}
& = & {\cal L}_{1}\frac{\partial}{\partial t^{\alpha}_{m}}
	\frac{\partial}{\partial t^{\beta}_{n}} F_{0} \\
& = & \left\{ \frac{\partial}{\partial t^{\alpha}_{m}} {\cal L}_{1}
		+ [{\cal L}_{1}, \, \frac{\partial}{\partial t^{\alpha}_{m}} ]
	\right\} \frac{\partial}{\partial t^{\beta}_{n}} F_{0} \\
& = & \frac{\partial}{\partial t^{\alpha}_{m}}
	\left\{ \frac{\partial}{\partial t^{\beta}_{n}} {\cal L}_{1}
		+ [{\cal L}_{1}, \, \frac{\partial}{\partial t^{\beta}_{n}} ]
	\right\}  F_{0} 
	+ [{\cal L}_{1}, \, \frac{\partial}{\partial t^{\alpha}_{m}} ]
	\frac{\partial}{\partial t^{\beta}_{n}} F_{0} \\
& = & \frac{\partial}{\partial t^{\alpha}_{m}}
	\frac{\partial}{\partial t^{\beta}_{n}} 
	\left<\left< {\cal L}_{1} \right>\right>_{0}
	+ \frac{\partial}{\partial t^{\alpha}_{m}}
	  [{\cal L}_{1}, \, \frac{\partial}{\partial t^{\beta}_{n}} ]
	  F_{0} 
	+ [{\cal L}_{1}, \, \frac{\partial}{\partial t^{\alpha}_{m}} ]
	\frac{\partial}{\partial t^{\beta}_{n}} F_{0}. \\
\end{eqnarray*}
The lemma then follows from (\ref{eqn:1ptL1}) and 
(\ref{eqn:L1Bracket}).
$\Box$

We can then compute
\begin{lem}  \label{lem:L1L0Corr}
\begin{eqnarray*}
  && \left<\left< {\cal L}_{1} ({\cal L}_{0} - 2{\cal D})
                \tau_{n}({\cal O}_{\beta}) \right>\right>_{0} \\
 &=&	-\sum_{m, \alpha} (m+b_{\alpha})(m+b_{\alpha}+1) 
                (m+b_{\alpha}+2) 
        \tilde{t}^{\alpha}_{m} \left<\left< \tau_{m+1}({\cal O}_{\alpha})
		\tau_{n}({\cal O}_{\beta})
                \right>\right>_{0}   \\
        && - \sum_{m, \alpha, \sigma} 
        \left\{ 3(m+b_{\alpha})^{2} + 6(m+b_{\alpha}) +2 \right\} 
     		{\cal C}_{\alpha}^{\sigma}
        	\tilde{t}^{\alpha}_{m} 
		\left<\left< \tau_{m}({\cal O}_{\sigma})
		\tau_{n}({\cal O}_{\beta})
			\right>\right>_{0}   \\
        && - \sum_{m, \alpha, \sigma} 3(m+b_{\alpha}+1) 
                ({\cal C}^{2})_{\alpha}^{\sigma}
        	\tilde{t}^{\alpha}_{m} 
		\left<\left< \tau_{m-1}({\cal O}_{\sigma})
		\tau_{n}({\cal O}_{\beta})
                \right>\right>_{0}   \\
        && - \sum_{m, \alpha, \sigma} 
                ({\cal C}^{3})_{\alpha}^{\sigma}
        	\tilde{t}^{\alpha}_{m} 
		\left<\left< \tau_{m-2}({\cal O}_{\sigma})
		\tau_{n}({\cal O}_{\beta})
                \right>\right>_{0}   \\
	&&+ (n+b_{\beta})(n+b_{\beta}+1) 
                (n+b_{\beta}-1) 
         	\left<\left< \tau_{n+1}({\cal O}_{\beta})
                \right>\right>_{0}  \\
        && + \sum_{\sigma} 
        	\left\{ 3(n+b_{\beta})^{2} - 1 \right\} 
	    	 {\cal C}_{\beta}^{\sigma}
        	\left<\left< \tau_{n}({\cal O}_{\sigma})
		\right>\right>_{0}   \\
        && + \sum_{\sigma} 3(n+b_{\beta}) 
                ({\cal C}^{2})_{\beta}^{\sigma}
        	\left<\left< \tau_{n-1}({\cal O}_{\sigma})
                \right>\right>_{0}   \\
        && + \sum_{\sigma} 
                ({\cal C}^{3})_{\beta}^{\sigma}
         	\left<\left< \tau_{n-2}({\cal O}_{\sigma})
                \right>\right>_{0}  \\
	&&- \sum_{\sigma} (b_{\sigma}-1)b_{\sigma}(n+b_{\beta}-1)
         	\left<\left< {\cal O}^{\sigma}
                \right>\right>_{0} 
	 	\left<\left< {\cal O}_{\sigma} \tau_{n}({\cal O}_{\beta})
                \right>\right>_{0} \\
	&&- \sum_{\sigma, \rho} (b_{\sigma}-1)b_{\sigma}
		{\cal C}_{\beta}^{\rho}
         	 \left<\left< \tau_{n-1}({\cal O}_{\rho}) {\cal O}_{\sigma}
                \right>\right>_{0}
	 	\left<\left< {\cal O}^{\sigma}
                	\right>\right>_{0}  \\
	&& - \delta_{n, 0} \left\{
		\sum_{\sigma}b_{\beta}(b_{\beta}+1){\cal C}_{\beta\sigma}
		\left<\left< {\cal O}^{\sigma}
                \right>\right>_{0} 
		-3b_{\beta} \sum_{\sigma} ({\cal C}^{2})_{\beta \sigma}
		t^{\sigma}_{0} 
		+ \sum_{\sigma}({\cal C}^{3})_{\beta\sigma}
		\tilde{t}^{\sigma}_{1} \right\} \\
	&& + \delta_{n, 1} \sum_{\sigma}({\cal C}^{3})_{\beta\sigma} 
		t^{\sigma}_{0}.
\end{eqnarray*}  
\end{lem}
{\bf Proof}: Using Lemma~\ref{lem:L1Corr}, we have
\begin{eqnarray*}
  && \left<\left< {\cal L}_{1} ({\cal L}_{0} - 2{\cal D})
                \tau_{n}({\cal O}_{\beta}) \right>\right>_{0} \\
&=& \sum_{m, \alpha} (m + b_{\alpha} +2)\tilde{t}^{\alpha}_{m}
         \left<\left< {\cal L}_{1} \tau_{m}({\cal O}_{\alpha}) 
                \tau_{n}({\cal O}_{\beta}) \right>\right>_{0}
        + \sum_{m, \alpha, \sigma} 
        {\cal C}_{\alpha}^{\sigma}\tilde{t}^{\alpha}_{m}
        \left<\left< {\cal L}_{1} \tau_{m-1}({\cal O}_{\sigma}) 
                \tau_{n}({\cal O}_{\beta}) \right>\right>_{0}  \\
& = &  - \sum_{m, \alpha}
	(m+b_{\alpha})(m+b_{\alpha}+1)  
	(m + b_{\alpha} +2)\tilde{t}^{\alpha}_{m}
    \left<\left<\tau_{m+1}({\cal O}_{\alpha}) 
                \tau_{n}({\cal O}_{\beta}) \right>\right>_{0} 
                \\
 && - \sum_{m, \alpha, \sigma} (2m+2b_{\alpha}+1) (m + b_{\alpha} +2)
	   {\cal C}_{\alpha}^{\sigma} \tilde{t}^{\alpha}_{m} 
            \left<\left<\tau_{m}({\cal O}_{\sigma}) 
                \tau_{n}({\cal O}_{\beta}) \right>\right>_{0} 
           \\
 && - \sum_{m, \alpha, \sigma} (m + b_{\alpha} +2)
                ({\cal C}^{2})_{\alpha}^{\sigma} \tilde{t}^{\alpha}_{m}
            \left<\left<\tau_{m-1}({\cal O}_{\sigma}) 
                \tau_{n}({\cal O}_{\beta}) \right>\right>_{0} \\
 && - \sum_{m, \alpha} (m + b_{\alpha} +2)
	(n+b_{\beta})(n+b_{\beta}+1)\tilde{t}^{\alpha}_{m}
            \left<\left<\tau_{m}({\cal O}_{\alpha}) 
                \tau_{n+1}({\cal O}_{\beta}) \right>\right>_{0} \\
 && - \sum_{m, \alpha, \sigma} 
	(m + b_{\alpha} +2)
	(2n+2b_{\beta}+1)
        {\cal C}_{\beta}^{\sigma} \tilde{t}^{\alpha}_{m} 
            \left<\left<\tau_{m}({\cal O}_{\alpha}) 
                \tau_{n}({\cal O}_{\sigma}) \right>\right>_{0} 
           \\
&& - \sum_{m, \alpha, \sigma} (m + b_{\alpha} +2)
                ({\cal C}^{2})_{\beta}^{\sigma}\tilde{t}^{\alpha}_{m}
            \left<\left<\tau_{m}({\cal O}_{\alpha}) 
                \tau_{n-1}({\cal O}_{\sigma}) \right>\right>_{0} \\
 && - \sum_{m, \alpha, \sigma} (m + b_{\alpha} +2)
	 b_{\sigma}(1-b_{\sigma})\tilde{t}^{\alpha}_{m}
                 \left<\left< \tau_{m}({\cal O}_{\alpha})
                \tau_{n}({\cal O}_{\beta}) 
                {\cal O}_{\sigma}\right>\right>_{0} 
                \left<\left< {\cal O}^{\sigma}\right>\right>_{0} 
                 \\
 && - \sum_{m, \alpha, \sigma} (m + b_{\alpha} +2)
	b_{\sigma}(1-b_{\sigma}) \tilde{t}^{\alpha}_{m}
                 \left<\left< \tau_{m}({\cal O}_{\alpha})
                {\cal O}_{\sigma}\right>\right>_{0} 
                \left<\left< {\cal O}^{\sigma}
                        \tau_{n}({\cal O}_{\beta}) \right>\right>_{0} 
	        \\
&&  - \delta_{n, 0} \sum_{\alpha} (b_{\alpha} +2)
                ({\cal C}^{2})_{\alpha \beta} 
		\tilde{t}^{\alpha}_{0}   \\
&& - \sum_{m \geq 1} \sum_{\alpha, \sigma} (m+b_{\sigma}-1)(m+b_{\sigma})
		{\cal C}_{\alpha}^{\sigma} 
		\tilde{t}^{\alpha}_{m}
    		\left<\left<\tau_{m}({\cal O}_{\sigma}) 
                \tau_{n}({\cal O}_{\beta}) \right>\right>_{0} 
                \\
        && - \sum_{m, \alpha, \sigma, \rho} (2m+2b_{\sigma}-1)
        {\cal C}_{\alpha}^{\sigma} {\cal C}_{\sigma}^{\rho}
		\tilde{t}^{\alpha}_{m}
            \left<\left<\tau_{m-1}({\cal O}_{\rho}) 
                \tau_{n}({\cal O}_{\beta}) \right>\right>_{0} 
           \\
        && - \sum_{m, \alpha, \sigma, \rho} 
		{\cal C}_{\alpha}^{\sigma}
                ({\cal C}^{2})_{\sigma}^{\rho}
		\tilde{t}^{\alpha}_{m}
            \left<\left<\tau_{m-2}({\cal O}_{\rho}) 
                \tau_{n}({\cal O}_{\beta}) \right>\right>_{0} \\
        && - \sum_{m, \alpha, \sigma} (n+b_{\beta})(n+b_{\beta}+1)
		{\cal C}_{\alpha}^{\sigma}
		\tilde{t}^{\alpha}_{m}
            \left<\left<\tau_{m-1}({\cal O}_{\sigma}) 
                \tau_{n+1}({\cal O}_{\beta}) \right>\right>_{0} \\
        && - \sum_{m, \alpha, \sigma, \rho} (2n+2b_{\beta}+1)
        	{\cal C}_{\alpha}^{\sigma}
		{\cal C}_{\beta}^{\rho}
		\tilde{t}^{\alpha}_{m}
            \left<\left<\tau_{m-1}({\cal O}_{\sigma}) 
                \tau_{n}({\cal O}_{\rho}) \right>\right>_{0} 
           \\
        && -\sum_{m, \alpha, \sigma, \rho} 
		{\cal C}_{\alpha}^{\sigma}
                ({\cal C}^{2})_{\beta}^{\rho}
		\tilde{t}^{\alpha}_{m}
            \left<\left<\tau_{m-1}({\cal O}_{\sigma}) 
                \tau_{n-1}({\cal O}_{\rho}) \right>\right>_{0} \\
        && - \sum_{m, \alpha, \sigma, \rho} b_{\rho}(1-b_{\rho})
		{\cal C}_{\alpha}^{\sigma}
		\tilde{t}^{\alpha}_{m}
                 \left<\left< \tau_{m-1}({\cal O}_{\sigma})
                \tau_{n}({\cal O}_{\beta}) 
                {\cal O}_{\rho}\right>\right>_{0} 
                \left<\left< {\cal O}^{\rho}\right>\right>_{0} 
                 \\
        && - \sum_{m, \alpha, \sigma, \rho} b_{\rho}(1-b_{\rho})
		{\cal C}_{\alpha}^{\sigma}
			\tilde{t}^{\alpha}_{m}
                 \left<\left< \tau_{m-1}({\cal O}_{\sigma})
                {\cal O}_{\rho}\right>\right>_{0} 
                \left<\left< {\cal O}^{\rho}
                        \tau_{n}({\cal O}_{\beta}) \right>\right>_{0} 
                 \\
        &&  - \delta_{n, 0} \sum_{\alpha, \sigma} 
		{\cal C}_{\alpha}^{\sigma}
                ({\cal C}^{2})_{\sigma \beta}
		\tilde{t}^{\alpha}_{1}.
\end{eqnarray*}
The fourth term and the thirteenth term can be combined together to
produce a correlation  function involving ${\cal L}_{0} - 2{\cal D}$.
The same is true when combining together the fifth term and the
 fourteenth term,
the sixth term and the fifteenth term, the seventh term and the sixteenth term,
the eighth term and the seventeenth term.
Using the fact that $b_{\sigma} = b_{\alpha} + 1$ if 
${\cal C}_{\alpha}^{\sigma} \neq 0$, we can simplify the above expression as
\begin{eqnarray} 
  && \left<\left< {\cal L}_{1} ({\cal L}_{0} - 2{\cal D})
                \tau_{n}({\cal O}_{\beta}) \right>\right>_{0}
		\nonumber \\
& = &  - \sum_{m, \alpha}
	(m+b_{\alpha})(m+b_{\alpha}+1)  
	(m + b_{\alpha} +2)\tilde{t}^{\alpha}_{m}
    \left<\left<\tau_{m+1}({\cal O}_{\alpha}) 
                \tau_{n}({\cal O}_{\beta}) \right>\right>_{0} 
               \nonumber \\
 && - \sum_{m, \alpha, \sigma} 
	\left\{3(m+b_{\alpha})^{2} + 6(m + b_{\alpha}) +2 \right\}
	   {\cal C}_{\alpha}^{\sigma} \tilde{t}^{\alpha}_{m} 
            \left<\left<\tau_{m}({\cal O}_{\sigma}) 
                \tau_{n}({\cal O}_{\beta}) \right>\right>_{0} 
          \nonumber  \\
&& +  \sum_{\alpha, \sigma} 
	b_{\alpha}(b_{\alpha}+1)
	   {\cal C}_{\alpha}^{\sigma} \tilde{t}^{\alpha}_{0} 
            \left<\left<{\cal O}_{\sigma}
                \tau_{n}({\cal O}_{\beta}) \right>\right>_{0} 
           \nonumber \\
 && - \sum_{m, \alpha, \sigma} 3(m + b_{\alpha} + 1)
                ({\cal C}^{2})_{\alpha}^{\sigma} \tilde{t}^{\alpha}_{m}
            \left<\left<\tau_{m-1}({\cal O}_{\sigma}) 
                \tau_{n}({\cal O}_{\beta}) \right>\right>_{0}
	\nonumber \\
&& - \sum_{m, \alpha, \rho} 
		({\cal C}^{3})_{\alpha}^{\rho}
		\tilde{t}^{\alpha}_{m}
            \left<\left<\tau_{m-2}({\cal O}_{\rho}) 
                \tau_{n}({\cal O}_{\beta}) \right>\right>_{0}
		\nonumber \\
 && - (n+b_{\beta})(n+b_{\beta}+1)
            \left<\left< ({\cal L}_{0} - 2{\cal D})
                \tau_{n+1}({\cal O}_{\beta}) \right>\right>_{0}
		\nonumber  \\
 && - \sum_{\sigma} 
	(2n+2b_{\beta}+1)
        {\cal C}_{\beta}^{\sigma}
            \left<\left< ({\cal L}_{0} - 2{\cal D})
                \tau_{n}({\cal O}_{\sigma}) \right>\right>_{0} 
           \nonumber \\
&& - \sum_{\sigma}
                ({\cal C}^{2})_{\beta}^{\sigma}
            \left<\left<({\cal L}_{0} - 2{\cal D})
                \tau_{n-1}({\cal O}_{\sigma}) \right>\right>_{0} 
	\nonumber \\
 && - \sum_{\sigma} 
	 b_{\sigma}(1-b_{\sigma})
                 \left<\left< ({\cal L}_{0} - 2{\cal D})
                \tau_{n}({\cal O}_{\beta}) 
                {\cal O}_{\sigma}\right>\right>_{0} 
                \left<\left< {\cal O}^{\sigma}\right>\right>_{0} 
                \nonumber  \\
 && - \sum_{\sigma} 
	b_{\sigma}(1-b_{\sigma}) 
                 \left<\left< ({\cal L}_{0} - 2{\cal D})
                {\cal O}_{\sigma}\right>\right>_{0} 
                \left<\left< {\cal O}^{\sigma}
                        \tau_{n}({\cal O}_{\beta}) \right>\right>_{0} 
	        \nonumber \\
&&  - \delta_{n, 0} \left\{
	\sum_{\alpha} (b_{\alpha} +2)
                ({\cal C}^{2})_{\alpha \beta}\tilde{t}^{\alpha}_{0} +
        \sum_{\alpha} 
		({\cal C}^{3})_{\alpha \beta}
		\tilde{t}^{\alpha}_{1} \right\}. \label{eqn:L1L0Corr:2} 
\end{eqnarray}
By Lemma~\ref{lem:DilatonCorr} (2) and Lemma~\ref{lem:EulerCorr} (2),
\begin{eqnarray*}
&& \left<\left< ({\cal L}_{0} - 2{\cal D})
		 \tau_{n}({\cal O}_{\beta}) \right>\right>_{0} \\
& = &    - \left(n + b_{\beta} -2 \right)
                \left<\left< \tau_{n}({\cal O}_{\beta}) \right>\right>_{0}
                - \sum_{\sigma} {\cal C}_{\beta}^{\sigma}
        \left<\left< \tau_{n-1}({\cal O}_{\sigma}) \right>\right>_{0}
        - \delta_{n, 0} \sum_{\sigma} {\cal C}_{\beta \sigma} t^{\sigma}_{0},
\end{eqnarray*}
and by Lemma~\ref{lem:DilatonCorr} (3) and Lemma~\ref{lem:EulerCorr} (3),
\begin{eqnarray*}
&&  \left<\left< ({\cal L}_{0} - 2 {\cal D})
	 \tau_{n}({\cal O}_{\beta}) 
                {\cal O}_{\rho}) \right>\right>_{0}  \\
&  = & - (n+b_{\beta}+b_{\rho}) \left<\left<\tau_{n}({\cal O}_{\beta}) 
                {\cal O}_{\rho} \right>\right>_{0} 
 	- \sum_{\gamma} {\cal C}_{\beta}^{\gamma}
        \left<\left<\tau_{n-1}({\cal O}_{\gamma}) 
	{\cal O}_{\rho} \right>\right>_{0} 
	-   \delta_{n,0} {\cal C}_{\beta \rho}.
\end{eqnarray*} 
Applying these two formulas to the right hand side of (\ref{eqn:L1L0Corr:2})
and simplying, we obtain the desired formula.
$\Box$

Setting $n = 0$ in Lemma~\ref{lem:L1L0Corr},
multiplying both sides of the
equation by $\left<\left< {\cal O}^{\beta} \tau_{k}({\cal O}_{\mu})
                \tau_{l}({\cal O}_{\nu}) \right>\right>_{0}$, and
summing over $\beta$, then applying the genus-$0$ topological recursion
relation, we get
\begin{eqnarray} \label{eqn:2WDVVleft}
  &&  \sum_{\beta}
	\left<\left< {\cal L}_{1} ({\cal L}_{0} - 2{\cal D})
                {\cal O}_{\beta} \right>\right>_{0}
	\left<\left< {\cal O}^{\beta} 
		\tau_{k}({\cal O}_{\mu})
                \tau_{l}({\cal O}_{\nu}) 
                \right>\right>_{0} \nonumber \\
 &=&	- \sum_{m, \alpha} (m+b_{\alpha})(m+b_{\alpha}+1) 
                (m+b_{\alpha}+2) 
        \tilde{t}^{\alpha}_{m} \left<\left< \tau_{m+2}({\cal O}_{\alpha})
		\tau_{k}({\cal O}_{\mu})
                \tau_{l}({\cal O}_{\nu}) 
                \right>\right>_{0}   \nonumber \\
        && - \sum_{m, \alpha, \sigma} 
        \left\{ 3(m+b_{\alpha})^{2} + 6(m+b_{\alpha}) +2 \right\} 
     {\cal C}_{\alpha}^{\sigma}
        \tilde{t}^{\alpha}_{m} \left<\left< \tau_{m+1}({\cal O}_{\sigma})
		\tau_{k}({\cal O}_{\mu})
                \tau_{l}({\cal O}_{\nu}) 
	\right>\right>_{0}  \nonumber  \\
        && - \sum_{m, \alpha, \sigma} 3(m+b_{\alpha}+1) 
                ({\cal C}^{2})_{\alpha}^{\sigma}
        \tilde{t}^{\alpha}_{m} \left<\left< \tau_{m}({\cal O}_{\sigma})
		\tau_{k}({\cal O}_{\mu})
                \tau_{l}({\cal O}_{\nu}) 
                \right>\right>_{0}  \nonumber \\
        && - \sum_{m, \alpha, \sigma} 
                ({\cal C}^{3})_{\alpha}^{\sigma}
        \tilde{t}^{\alpha}_{m} \left<\left< \tau_{m-1}({\cal O}_{\sigma})
		\tau_{k}({\cal O}_{\mu})
                \tau_{l}({\cal O}_{\nu}) 
                \right>\right>_{0}  \nonumber  \\
	&&+ \sum_{\beta} (b_{\beta}-1) b_{\beta}
                (b_{\beta}+1) 
         \left<\left< \tau_{1}({\cal O}_{\beta})
                \right>\right>_{0}  
		\left<\left< {\cal O}^{\beta} 
		\tau_{k}({\cal O}_{\mu})
                \tau_{l}({\cal O}_{\nu}) 
                \right>\right>_{0} \nonumber \\
        && + \sum_{\beta, \sigma} 
          (3 b_{\beta}^{2} - 1)
	     {\cal C}_{\beta}^{\sigma}
        \left<\left< {\cal O}_{\sigma}
	\right>\right>_{0} 
		\left<\left< {\cal O}^{\beta} 
		\tau_{k}({\cal O}_{\mu})
                \tau_{l}({\cal O}_{\nu}) 
                \right>\right>_{0} \nonumber \\
     	&&- \sum_{\beta, \sigma} (b_{\sigma}-1)b_{\sigma}(b_{\beta}-1)
         \left<\left< {\cal O}^{\sigma}
                \right>\right>_{0} 
	 \left<\left< {\cal O}_{\sigma} {\cal O}_{\beta}
                \right>\right>_{0} 
		\left<\left< {\cal O}^{\beta} 
		\tau_{k}({\cal O}_{\mu})
                \tau_{l}({\cal O}_{\nu}) 
                \right>\right>_{0} \nonumber \\
	&& - 	\sum_{\beta, \sigma}
		b_{\beta}(b_{\beta}+1){\cal C}_{\beta\sigma}
		\left<\left< {\cal O}^{\sigma}
                \right>\right>_{0} 
		\left<\left< {\cal O}^{\beta} 
		\tau_{k}({\cal O}_{\mu})
                \tau_{l}({\cal O}_{\nu}) 
                \right>\right>_{0}. 
\end{eqnarray}  
Notice that the ranges of summations may change when using the topological
recursion relation. Hence some scattered terms may be absorbed into a
big summation after using the topological recursion relation.

By Lemma~\ref{lem:DilatonCorr} (3) and the definition of ${\cal L}_{0}$,
\[	\left<\left< ({\cal L}_{0} - 2{\cal D})
               	\tau_{m}({\cal O}_{\alpha})
                \tau_{n}({\cal O}_{\beta}) 
                \right>\right>_{0}  
	= - \left<\left< {\cal X}
               	\tau_{m}({\cal O}_{\alpha})
                \tau_{n}({\cal O}_{\beta}) 
                \right>\right>_{0},    \]
for any $m, n, \alpha, \beta$. Hence by Lemma~\ref{lem:L1Corr}, we 
have
\begin{eqnarray} \label{eqn:2WDVVright}
  &&  \sum_{\beta}
	\left<\left< {\cal L}_{1} \tau_{k}({\cal O}_{\mu})
                {\cal O}_{\beta} \right>\right>_{0}
	\left<\left< {\cal O}^{\beta} 
		({\cal L}_{0} - 2{\cal D})
                \tau_{l}({\cal O}_{\nu}) 
                \right>\right>_{0} \nonumber \\
&=&  \sum_{\beta} (k+b_{\mu})(k+b_{\mu}+1) 
    \left<\left<\tau_{k+1}({\cal O}_{\mu}) 
                {\cal O}_{\beta} \right>\right>_{0} 
	 \left<\left< {\cal O}^{\beta} {\cal X}
                \tau_{l}({\cal O}_{\nu}) 
                \right>\right>_{0} 
               \nonumber   \\
	&& + \sum_{\beta, \sigma} (2k+2b_{\mu}+1)
	{\cal C}_{\mu}^{\sigma}
            \left<\left<\tau_{k}({\cal O}_{\sigma}) 
                {\cal O}_{\beta} \right>\right>_{0} 
	 \left<\left< {\cal O}^{\beta} {\cal X}
                \tau_{l}({\cal O}_{\nu}) 
                \right>\right>_{0} 
           \nonumber  \\
        && + \sum_{\beta, \sigma} 
                ({\cal C}^{2})_{\mu}^{\sigma}
            \left<\left<\tau_{k-1}({\cal O}_{\sigma}) 
                {\cal O}_{\beta} \right>\right>_{0} 
	 \left<\left< {\cal O}^{\beta} {\cal X}
                \tau_{l}({\cal O}_{\nu}) 
                \right>\right>_{0} 
		\nonumber \\
	&& + \sum_{\beta} b_{\beta}(b_{\beta}+1)
	    \left<\left<\tau_{k}({\cal O}_{\mu}) 
                \tau_{1}({\cal O}_{\beta}) \right>\right>_{0} 
	 \left<\left< {\cal O}^{\beta} {\cal X}
                \tau_{l}({\cal O}_{\nu}) 
                \right>\right>_{0} 
		\nonumber  \\
	&& + \sum_{\beta, \sigma} (2b_{\beta}+1)
	{\cal C}_{\beta}^{\sigma}
            \left<\left<\tau_{k}({\cal O}_{\mu}) 
                {\cal O}_{\sigma} \right>\right>_{0} 
 	\left<\left< {\cal O}^{\beta} {\cal X}
                \tau_{l}({\cal O}_{\nu}) 
                \right>\right>_{0} 
          \nonumber   \\
     	&& + \sum_{\beta, \sigma} b_{\sigma}(1-b_{\sigma})
		 \left<\left< \tau_{k}({\cal O}_{\mu})
                {\cal O}_{\beta}
                {\cal O}_{\sigma}\right>\right>_{0} 
                \left<\left< {\cal O}^{\sigma}\right>\right>_{0} 
 	\left<\left< {\cal O}^{\beta} {\cal X}
                \tau_{l}({\cal O}_{\nu}) 
                \right>\right>_{0} 
                \nonumber  \\
	&& + \sum_{\beta, \sigma} b_{\sigma}(1-b_{\sigma})
		 \left<\left< \tau_{k}({\cal O}_{\mu})
                {\cal O}_{\sigma}\right>\right>_{0} 
                \left<\left< {\cal O}^{\sigma}
			{\cal O}_{\beta} \right>\right>_{0} 
		 \left<\left< {\cal O}^{\beta} {\cal X}
                \tau_{l}({\cal O}_{\nu}) 
                \right>\right>_{0} 
                 \nonumber  \\
        &&   + \delta_{k, 0}\sum_{\beta}
                ({\cal C}^{2})_{\mu \beta}
	 \left<\left< {\cal O}^{\beta} {\cal X}
                \tau_{l}({\cal O}_{\nu}) 
                \right>\right>_{0}.
	\end{eqnarray}
By Lemma~\ref{lem:EulerCorr} (3), 
we know that the right hand side
of (\ref{eqn:2WDVVright}) has only finitely many terms. As in the proof
of the $L_{1}$ constraint, we will show that the difference
of the right hand sides of (\ref{eqn:2WDVVright}) 
and (\ref{eqn:2WDVVleft})
is equal to $\frac{\partial}{\partial t^{\nu}_{l}}
        \frac{\partial}{\partial t^{\mu}_{k}} \Psi_{0, 2}$.
For this purpose, we need to express the products of correlation
functions in  (\ref{eqn:2WDVVright}) as summations of correlation
functions. 

Using the generalized WDVV equation to the sixth term and 
the topological recursion relation to
the first three terms and the seventh term on the right hand side
of (\ref{eqn:2WDVVright}), we obtain
\begin{eqnarray*} 
  &&  \sum_{\beta}
	\left<\left< {\cal L}_{1} \tau_{k}({\cal O}_{\mu})
                {\cal O}_{\beta} \right>\right>_{0}
	\left<\left< {\cal O}^{\beta} 
		({\cal L}_{0} - 2{\cal D})
                \tau_{l}({\cal O}_{\nu}) 
                \right>\right>_{0} \nonumber \\
&=&    (k+b_{\mu})(k+b_{\mu}+1) 
    \left<\left<\tau_{k+2}({\cal O}_{\mu}) 
                 {\cal X}
                \tau_{l}({\cal O}_{\nu}) 
                \right>\right>_{0} 
               \nonumber   \\
	&& + \sum_{\sigma} (2k+2b_{\mu}+1)
	{\cal C}_{\mu}^{\sigma}
            \left<\left<\tau_{k+1}({\cal O}_{\sigma}) 
                {\cal X}
                \tau_{l}({\cal O}_{\nu}) 
                \right>\right>_{0} 
           \nonumber  \\
        && + \sum_{\sigma} 
                ({\cal C}^{2})_{\mu}^{\sigma}
            \left<\left<\tau_{k}({\cal O}_{\sigma}) 
                 {\cal X}
                \tau_{l}({\cal O}_{\nu}) 
                \right>\right>_{0} 
	  - \delta_{k,0}\sum_{\sigma} 
                ({\cal C}^{2})_{\mu}^{\sigma}
            \left<\left<{\cal O}_{\sigma}
                 {\cal X}
                \tau_{l}({\cal O}_{\nu}) 
                \right>\right>_{0} 
		\nonumber \\
	&& + \sum_{\beta} b_{\beta}(b_{\beta}+1)
	    \left<\left<\tau_{k}({\cal O}_{\mu}) 
                \tau_{1}({\cal O}_{\beta}) \right>\right>_{0} 
	 \left<\left< {\cal O}^{\beta} {\cal X}
                \tau_{l}({\cal O}_{\nu}) 
                \right>\right>_{0} 
		\nonumber  \\
	&& + \sum_{\beta, \sigma} (2b_{\beta}+1)
	{\cal C}_{\beta}^{\sigma}
            \left<\left<\tau_{k}({\cal O}_{\mu}) 
                {\cal O}_{\sigma} \right>\right>_{0} 
 	\left<\left< {\cal O}^{\beta} {\cal X}
                \tau_{l}({\cal O}_{\nu}) 
                \right>\right>_{0} 
          \nonumber   \\
     	&& + \sum_{\beta, \sigma} b_{\sigma}(1-b_{\sigma})
		 \left<\left< {\cal X}
                {\cal O}_{\beta}
                {\cal O}_{\sigma}\right>\right>_{0} 
                \left<\left< {\cal O}^{\sigma}\right>\right>_{0} 
 	\left<\left< {\cal O}^{\beta} \tau_{k}({\cal O}_{\mu})
                \tau_{l}({\cal O}_{\nu}) 
                \right>\right>_{0} 
                \nonumber  \\
	&& + \sum_{\sigma} b_{\sigma}(1-b_{\sigma})
		 \left<\left< \tau_{k}({\cal O}_{\mu})
                {\cal O}^{\sigma}\right>\right>_{0} 
                \left<\left< \tau_{1}({\cal O}_{\sigma})
			{\cal X}
                \tau_{l}({\cal O}_{\nu}) 
                \right>\right>_{0} 
                 \nonumber  \\
        &&   + \delta_{k, 0}\sum_{\beta}
                ({\cal C}^{2})_{\mu \beta}
	 \left<\left< {\cal O}^{\beta} {\cal X}
                \tau_{l}({\cal O}_{\nu}) 
                \right>\right>_{0}.
	\end{eqnarray*}
We then use Lemma~\ref{lem:EulerCorr} (3)
to get rid of the Euler vector fields in the above expression.
After simplification, we will find many terms which also appear
in $\frac{\partial}{\partial t^{\nu}_{l}}
        \frac{\partial}{\partial t^{\mu}_{k}} \Psi_{0, 2}$
or in the right hand side of (\ref{eqn:2WDVVleft}).
We call such terms {\it good terms}.
There are also many terms which do not appear in 
$\frac{\partial}{\partial t^{\nu}_{l}}
        \frac{\partial}{\partial t^{\mu}_{k}} \Psi_{0, 2}$,
nor do they appear in the right hand side of (\ref{eqn:2WDVVleft}).
We call such terms {\it bad terms}. Grouping good terms together
and bad terms together, we obtain the following
\begin{eqnarray} \label{eqn:2WDVVright:2}
  &&  \sum_{\beta}
	\left<\left< {\cal L}_{1} \tau_{k}({\cal O}_{\mu})
                {\cal O}_{\beta} \right>\right>_{0}
	\left<\left< {\cal O}^{\beta} 
		({\cal L}_{0} - 2{\cal D})
                \tau_{l}({\cal O}_{\nu}) 
                \right>\right>_{0} \nonumber \\
&=&    (k+b_{\mu})(k+b_{\mu}+1) (k+b_{\mu}+2) 
   		 \left<\left<\tau_{k+2}({\cal O}_{\mu}) 
                \tau_{l}({\cal O}_{\nu}) 
                \right>\right>_{0} 
               \nonumber   \\
	&& + \sum_{\sigma} \{3(k+b_{\mu})^{2} + 6(k+b_{\mu})+ 2 \}
	{\cal C}_{\mu}^{\sigma}
            \left<\left<\tau_{k+1}({\cal O}_{\sigma}) 
                \tau_{l}({\cal O}_{\nu}) 
                \right>\right>_{0} 
           \nonumber  \\
	&& + \sum_{\sigma} 3(k+b_{\mu}+1)
                ({\cal C}^{2})_{\mu}^{\sigma}
            \left<\left<\tau_{k}({\cal O}_{\sigma}) 
                \tau_{l}({\cal O}_{\nu}) 
                \right>\right>_{0}  \nonumber  \\
	&& + \sum_{\sigma}({\cal C}^{3})_{\mu}^{\sigma} 
	 	\left<\left<\tau_{k-1}({\cal O}_{\sigma}) 
                \tau_{l}({\cal O}_{\nu}) 
                \right>\right>_{0}  \nonumber  \\
	&& + \sum_{\beta} (1-b_{\beta})b_{\beta}(1+b_{\beta})
 		\left<\left< \tau_{k}({\cal O}_{\mu})
                \tau_{1}({\cal O}_{\beta})\right>\right>_{0} 
                \left<\left< {\cal O}^{\beta}) \tau_{l}({\cal O}_{\nu}) 
                \right>\right>_{0} 
                 \nonumber  \\
        && + \sum_{\alpha, \beta} b_{\alpha}(1-b_{\alpha})
		(b_{\alpha}+b_{\beta})
        \left<\left< {\cal O}^{\alpha}\right>\right>_{0}         
	 \left<\left<  {\cal O}_{\alpha}
                {\cal O}_{\beta}\right>\right>_{0} 
        \left<\left< {\cal O}^{\beta} \tau_{k}({\cal O}_{\mu})
                \tau_{l}({\cal O}_{\nu}) 
                \right>\right>_{0} 
                \nonumber  \\
        && +    \sum_{\alpha, \beta}
                b_{\alpha}(1-b_{\alpha}){\cal C}_{\alpha \beta}
                \left<\left< {\cal O}^{\alpha}
                \right>\right>_{0} 
                \left<\left< {\cal O}^{\beta} 
                \tau_{k}({\cal O}_{\mu})
                \tau_{l}({\cal O}_{\nu}) 
                \right>\right>_{0} \nonumber  \\
	&& + \sum_{\alpha}
                b_{\alpha}(1-b_{\alpha})(b_{\alpha}+1)
		\left<\left< \tau_{k}({\cal O}_{\mu}) {\cal O}^{\alpha}
                \right>\right>_{0} 
                \left<\left< \tau_{1}({\cal O}_{\alpha} )
                \tau_{l}({\cal O}_{\nu}) 
                \right>\right>_{0} \nonumber  \\
	&& - \sum_{\alpha, \beta}
                (3b_{\alpha}^{2}-1) {\cal C}_{\alpha}^{\beta}
		\left<\left< \tau_{k}({\cal O}_{\mu}) {\cal O}_{\beta}
                \right>\right>_{0} 
                \left<\left< {\cal O}^{\alpha}
                \tau_{l}({\cal O}_{\nu}) 
                \right>\right>_{0} \nonumber  \\
	&& + \delta_{k, 0} \delta_{l, 0} ({\cal C}^{3})_{\mu \nu} 
	\nonumber  \\
	&& +  (k+b_{\mu})(k+b_{\mu}+1) (l+b_{\nu}) 
   		 \left<\left<\tau_{k+2}({\cal O}_{\mu})
		\tau_{l}({\cal O}_{\nu}) 
                \right>\right>_{0} 
               \nonumber   \\
	&& + \sum_{\sigma} (2k+2b_{\mu}+1)(l+b_{\nu})
        {\cal C}_{\mu}^{\sigma}
            \left<\left<\tau_{k+1}({\cal O}_{\sigma}) 
                \tau_{l}({\cal O}_{\nu}) 
                \right>\right>_{0} 
           \nonumber  \\
	&& + \sum_{\sigma} (k+b_{\mu})(k+b_{\mu}+1)
		{\cal C}_{\nu}^{\sigma}
            \left<\left<\tau_{k+2}({\cal O}_{\mu}) 
                \tau_{l-1}({\cal O}_{\sigma}) 
                \right>\right>_{0} 
           \nonumber  \\
        && + \sum_{\sigma} (l+b_{\nu})
                ({\cal C}^{2})_{\mu}^{\sigma}
            \left<\left<\tau_{k}({\cal O}_{\sigma}) 
                \tau_{l}({\cal O}_{\nu}) 
                \right>\right>_{0}  \nonumber  \\
	&& +\sum_{\alpha, \beta}(2k+2b_{\mu}+1) 
		{\cal C}_{\mu}^{\alpha}{\cal C}_{\nu}^{\beta}
		\left<\left<\tau_{k+1}({\cal O}_{\alpha}) 
                \tau_{l-1}({\cal O}_{\beta}) 
                \right>\right>_{0} 
           \nonumber  \\
        && +\sum_{\alpha, \beta}
		({\cal C}^{2})_{\mu}^{\alpha}{\cal C}_{\nu}^{\beta}
		\left<\left<\tau_{k}({\cal O}_{\alpha}) 
                \tau_{l-1}({\cal O}_{\beta}) 
                \right>\right>_{0} 
           \nonumber  \\
	&& + \sum_{\beta} b_{\beta}(b_{\beta}+1)(l+b_{\nu})
 		\left<\left< \tau_{k}({\cal O}_{\mu})
                \tau_{1}({\cal O}_{\beta})\right>\right>_{0} 
                \left<\left< {\cal O}^{\beta} \tau_{l}({\cal O}_{\nu}) 
                \right>\right>_{0} 
                 \nonumber  \\
	&& + \sum_{\beta} b_{\beta}(1-b_{\beta})(l+b_{\nu})
 		\left<\left< \tau_{k}({\cal O}_{\mu})
                {\cal O}^{\beta})\right>\right>_{0} 
                \left<\left< \tau_{1}({\cal O}_{\beta}) 
		\tau_{l}({\cal O}_{\nu}) 
                \right>\right>_{0} 
                 \nonumber  \\
	&& + \sum_{\alpha, \beta} (2b_{\alpha}+1)(l+b_{\nu})
		{\cal C}_{\alpha}^{\beta}
 		\left<\left< \tau_{k}({\cal O}_{\mu})
                {\cal O}_{\beta}\right>\right>_{0} 
                \left<\left< {\cal O}^{\alpha} \tau_{l}({\cal O}_{\nu}) 
                \right>\right>_{0} 
                 \nonumber  \\
	&& + \delta_{l, 0} \left\{ \sum_{\beta} b_{\beta}(b_{\beta}+1)
			{\cal C}_{\nu}^{\beta} 
		\left<\left< \tau_{k}({\cal O}_{\mu}) 
		\tau_{1}({\cal O}_{\beta})
                \right>\right>_{0} 
		+ \sum_{\beta}(2b_{\nu}+3)({\cal C}^{2})_{\nu}^{\beta}
		\left<\left< \tau_{k}({\cal O}_{\mu}) {\cal O}_{\beta}
                \right>\right>_{0}  \right\}
                 \nonumber  \\
	&& + \sum_{\beta, \sigma} b_{\beta}(b_{\beta}+1)
		{\cal C}_{\nu}^{\sigma}
 		\left<\left< \tau_{k}({\cal O}_{\mu})
                \tau_{1}({\cal O}_{\beta})\right>\right>_{0} 
                \left<\left< {\cal O}^{\beta} \tau_{l-1}({\cal O}_{\sigma}) 
                \right>\right>_{0} 
                 \nonumber  \\
	&& + \sum_{\beta, \sigma} b_{\beta}(1- b_{\beta})
		{\cal C}_{\nu}^{\sigma}
 		\left<\left< \tau_{k}({\cal O}_{\mu})
                {\cal O}^{\beta}\right>\right>_{0} 
                \left<\left< \tau_{1}({\cal O}_{\beta}) 
		\tau_{l-1}({\cal O}_{\sigma}) 
                \right>\right>_{0} \nonumber  \\
	&& + \sum_{\alpha, \beta, \sigma} (2b_{\beta}+1)
		{\cal C}_{\beta}^{\alpha} {\cal C}_{\nu}^{\sigma} 
 		\left<\left< \tau_{k}({\cal O}_{\mu})
                {\cal O}_{\alpha}\right>\right>_{0} 
                \left<\left< {\cal O}^{\beta} \tau_{l-1}({\cal O}_{\sigma}) 
                \right>\right>_{0}, 
\end{eqnarray}
where the first 10 terms are good terms and  rest terms are bad terms.
In order to transform bad terms into good terms, we need more
properties for correlation functions.
\begin{lem} \label{lem:QuadRel}
\begin{eqnarray*}
&\textrm{(i)}& \textrm{ \hspace{10pt}}
	\sum_{\beta}\left<\left< {\cal X} \tau_{k}({\cal O}_{\mu})
                \tau_{1}({\cal O}_{\beta})\right>\right>_{0} 
                \left<\left< {\cal O}^{\beta}
		{\cal X} 
		\tau_{l}({\cal O}_{\nu}) 
                \right>\right>_{0}    \\
   &&	= \sum_{\beta}\left<\left< {\cal X} \tau_{k}({\cal O}_{\mu})
                {\cal O}^{\beta}\right>\right>_{0} 
                \left<\left< \tau_{1}({\cal O}_{\beta})
		{\cal X} 
		\tau_{l}({\cal O}_{\nu}) 
                \right>\right>_{0},   \\
&\textrm{(ii)}& \textrm{ \hspace{10pt}}
	\sum_{\beta}\left<\left< \tau_{k-1}({\cal O}_{\mu})
                {\cal O}^{\beta} \right>\right>_{0} 
                \left<\left< \tau_{1}({\cal O}_{\beta})
		{\cal X} 
		\tau_{l}({\cal O}_{\nu}) 
                \right>\right>_{0}   \\
   &&	= \sum_{\beta}\left<\left<  \tau_{k-1}({\cal O}_{\mu})
                \tau_{1}({\cal O}_{\beta})\right>\right>_{0} 
                \left<\left< {\cal O}^{\beta}
		{\cal X} 
		\tau_{l}({\cal O}_{\nu}) 
                \right>\right>_{0}  
	+ \left<\left<  \tau_{k+1}({\cal O}_{\mu})
               		{\cal X} 
		\tau_{l}({\cal O}_{\nu}) 
                \right>\right>_{0}    \\
  && \textrm{ \hspace{15pt}} 
	- \delta_{k, 0}  \left<\left<  \tau_{1}({\cal O}_{\mu})
               		{\cal X} 
		\tau_{l}({\cal O}_{\nu}) 
                \right>\right>_{0} 
	- \delta_{k, -1}  \left<\left<  {\cal O}_{\mu}
               		{\cal X} 
		\tau_{l}({\cal O}_{\nu}) 
                \right>\right>_{0},  \\
&\textrm{(iii)} &  \textrm{ \hspace{10pt}} 
	\sum_{\beta}\left<\left< \tau_{k}({\cal O}_{\mu})
                {\cal O}^{\beta} \right>\right>_{0} 
                \left<\left< \tau_{1}({\cal O}_{\beta})
		\tau_{l}({\cal O}_{\nu}) 
                \right>\right>_{0}    \\
&& =  \sum_{\beta}\left<\left< \tau_{k}({\cal O}_{\mu})
                \tau_{1}({\cal O}_{\beta}) \right>\right>_{0} 
                \left<\left< {\cal O}^{\beta}
		\tau_{l}({\cal O}_{\nu}) 
                \right>\right>_{0}    \\
&&  \textrm{ \hspace{15pt}} 
	+ \left<\left<  \tau_{k+2}({\cal O}_{\mu})
		\tau_{l}({\cal O}_{\nu}) 
                \right>\right>_{0}   
	- \left<\left<  \tau_{k}({\cal O}_{\mu})
		\tau_{l+2}({\cal O}_{\nu}) 
                \right>\right>_{0}    \\
&&  \textrm{ \hspace{15pt}} 
	- \delta_{k, -1}  \left<\left<  \tau_{1}({\cal O}_{\mu})
		\tau_{l}({\cal O}_{\nu}) 
                \right>\right>_{0}   
	- \delta_{k, -2}  \left<\left<  {\cal O}_{\mu}
		\tau_{l}({\cal O}_{\nu}) 
                \right>\right>_{0}    \\
&&  \textrm{ \hspace{15pt}} 
	+ \delta_{l, -1} \left<\left<  \tau_{k}({\cal O}_{\mu})
		\tau_{1}({\cal O}_{\nu}) 
                \right>\right>_{0}   
	+ \delta_{l, -2} \left<\left<  \tau_{k}({\cal O}_{\mu})
		{\cal O}_{\nu}
                \right>\right>_{0}.   
\end{eqnarray*}
\end{lem}
{\bf Proof}:
(i) follows by applying topological recursion relation to both sides
of the equation for the terms which contain $\tau_{1}$.

To prove (ii), we first use topological recursion relation,
then use formula (\ref{eqn:StringRec}). We have
\begin{eqnarray*}
&&	\sum_{\beta}\left<\left< \tau_{k-1}({\cal O}_{\mu})
                {\cal O}^{\beta} \right>\right>_{0} 
                \left<\left< \tau_{1}({\cal O}_{\beta})
		{\cal X} 
		\tau_{l}({\cal O}_{\nu}) 
                \right>\right>_{0}   \\
&=&	\sum_{\alpha, \beta}\left<\left< \tau_{k-1}({\cal O}_{\mu})
                {\cal O}^{\beta} \right>\right>_{0} 
                \left<\left< {\cal O}_{\beta}
		{\cal O}_{\alpha}
                \right>\right>_{0} 
                \left<\left< {\cal O}^{\alpha}
		{\cal X} 
		\tau_{l}({\cal O}_{\nu}) 
                \right>\right>_{0}   \\
& = & \sum_{\alpha} \left\{
		\left<\left< \tau_{k}({\cal O}_{\mu})
                {\cal O}_{\alpha}
                \right>\right>_{0} 
		+ \left<\left< \tau_{k-1}({\cal O}_{\mu})
                \tau_{1}({\cal O}_{\alpha}) 
                \right>\right>_{0} 
		- \delta_{k, 0} \left<\left< {\cal O}_{\mu}
                {\cal O}_{\alpha}
                \right>\right>_{0} 
		\right\}
                \left<\left< {\cal O}^{\alpha}
		{\cal X} 
		\tau_{l}({\cal O}_{\nu}) 
                \right>\right>_{0}.   
\end{eqnarray*}
Applying topological recursion relation again to the first and the third
terms, we obtain (ii).

We now prove (iii). Let
\begin{eqnarray*}
f &:= & \sum_{\beta}\left<\left< \tau_{k}({\cal O}_{\mu})
                {\cal O}^{\beta} \right>\right>_{0} 
                \left<\left< \tau_{1}({\cal O}_{\beta})
		\tau_{l}({\cal O}_{\nu}) 
                \right>\right>_{0}.    
\end{eqnarray*}
By (\ref{eqn:StringRec}),
\begin{eqnarray*}
f & = & \sum_{\beta}\left<\left< \tau_{k}({\cal O}_{\mu})
                {\cal O}^{\beta} \right>\right>_{0} 
	\left\{ - \left<\left< {\cal O}_{\beta}
		\tau_{l+1}({\cal O}_{\nu}) 
                \right>\right>_{0}
		+ \delta_{l, -1}
		 \left<\left< {\cal O}_{\beta}{\cal O}_{\nu} 
                \right>\right>_{0}    \right. \\
&& \textrm{ \hspace{120pt}} \left.
		+  \sum_{\alpha}\left<\left< 
		\tau_{l}({\cal O}_{\nu}) 
		{\cal O}^{\alpha}
                \right>\right>_{0}
		\left<\left< {\cal O}_{\alpha} 
		{\cal O}_{\beta}
                \right>\right>_{0} \right\}.       
\end{eqnarray*}
Using (\ref{eqn:StringRec}) reversely to each term, we have
\begin{eqnarray*}
f & = & - \left<\left< \tau_{k+1}({\cal O}_{\mu})
               	\tau_{l+1}({\cal O}_{\nu}) 
                \right>\right>_{0}
	- \left<\left< \tau_{k}({\cal O}_{\mu})
               	\tau_{l+2}({\cal O}_{\nu}) 
                \right>\right>_{0}      \\    
&&	+ \delta_{k, -1} \left<\left< {\cal O}_{\mu}
               	\tau_{l+1}({\cal O}_{\nu}) 
                \right>\right>_{0}
	+ \delta_{l, -2} \left<\left< \tau_{k}({\cal O}_{\mu})
               	{\cal O}_{\nu}
                \right>\right>_{0}  \\
&& + \delta_{l, -1}
	\left\{ \left<\left< \tau_{k+1}({\cal O}_{\mu})
                {\cal O}_{\nu} 
                \right>\right>_{0} 
		+ \left<\left< \tau_{k}({\cal O}_{\mu})
                \tau_{1}({\cal O}_{\nu})
                \right>\right>_{0} 
		- \delta_{k, -1} \left<\left< {\cal O}_{\mu}
                {\cal O}_{\nu} 
                \right>\right>_{0}  \right\} \\
&& + \sum_{\alpha}\left<\left< 
		\tau_{l}({\cal O}_{\nu}) 
		{\cal O}^{\alpha}
                \right>\right>_{0}   
		\left\{ \left<\left< 
                 \tau_{1}({\cal O}_{\alpha}) \tau_{k}({\cal O}_{\mu})
		       \right>\right>_{0}
		+ \left<\left< 
		{\cal O}_{\alpha} \tau_{k+1}({\cal O}_{\mu})
		       \right>\right>_{0} 
		- \delta_{k, -1} \left<\left< 
                 {\cal O}_{\alpha} {\cal O}_{\mu}
		       \right>\right>_{0} \right\}.       
\end{eqnarray*}
Applying (\ref{eqn:StringRec}) again to the last two terms and simplifying,
we obtain (iii).
$\Box$

Using Lemma~\ref{lem:QuadRel}, we can show
\begin{lem} \label{lem:QuadForm}
\begin{eqnarray*}
 && \sum_{\beta} b_{\beta}(b_{\beta}+1)
                \left<\left< \tau_{k}({\cal O}_{\mu})
                \tau_{1}({\cal O}_{\beta})\right>\right>_{0} 
                \left<\left< {\cal O}^{\beta} \tau_{l}({\cal O}_{\nu}) 
                \right>\right>_{0} 
                 \nonumber  \\
&&  + \sum_{\beta} b_{\beta}(1-b_{\beta})
                \left<\left< \tau_{k}({\cal O}_{\mu})
                {\cal O}^{\beta}\right>\right>_{0} 
                \left<\left< \tau_{1}({\cal O}_{\beta}) 
                \tau_{l}({\cal O}_{\nu}) 
                \right>\right>_{0} 
                 \nonumber  \\
&& + \sum_{\alpha, \beta}(2b_{\beta} +1) {\cal C}_{\alpha \beta}
                \left<\left< \tau_{k}({\cal O}_{\mu})
                {\cal O}^{\alpha} \right>\right>_{0} 
                \left<\left< {\cal O}^{\beta}
                \tau_{l}({\cal O}_{\nu}) 
                \right>\right>_{0} \\
&=&  - (k+b_{\mu})(k + b_{\mu}+1) 
	 \left<\left< \tau_{k+2}({\cal O}_{\mu})
                 \tau_{l}({\cal O}_{\nu}) 
                \right>\right>_{0}    \\
&& + (l+b_{\nu}+1)(l + b_{\nu}+2) 
	 \left<\left< \tau_{k}({\cal O}_{\mu})
                 \tau_{l+2}({\cal O}_{\nu}) 
                \right>\right>_{0}  \\
&&-\sum_{\alpha} (2k + 2b_{\mu}+1) {\cal C}_{\mu}^{\alpha} 
	 \left<\left< \tau_{k+1}({\cal O}_{\alpha})
                 \tau_{l}({\cal O}_{\nu}) 
                \right>\right>_{0}    \\
&& + \sum_{\alpha} (2l + 2b_{\nu}+3) {\cal C}_{\nu}^{\alpha} 
	 \left<\left< \tau_{k}({\cal O}_{\mu})
                 \tau_{l+1}({\cal O}_{\alpha}) 
                \right>\right>_{0}    \\
&& - \sum_{\alpha}({\cal C}^{2})_{\mu}^{\alpha} 
	 \left<\left< \tau_{k}({\cal O}_{\alpha})
                 \tau_{l}({\cal O}_{\nu}) 
                \right>\right>_{0}    \\
&& +  \sum_{\alpha}({\cal C}^{2})_{\nu}^{\alpha} 
	 \left<\left< \tau_{k}({\cal O}_{\mu})
                 \tau_{l}({\cal O}_{\alpha}) 
                \right>\right>_{0}    \\
&& - \delta_{l, -1}\left\{ b_{\nu}(b_{\nu}+1)
		\left<\left< \tau_{k}({\cal O}_{\mu}) 
		\tau_{1}({\cal O}_{\nu})
                \right>\right>_{0} 
		+ \sum_{\alpha}(2b_{\nu}+1){\cal C}_{\nu}^{\alpha}
		\left<\left< \tau_{k}({\cal O}_{\mu}) {\cal O}_{\alpha}
                \right>\right>_{0}  \right\}.
\end{eqnarray*}
\end{lem}
{\bf Proof}: By Lemma~\ref{lem:EulerCorr} (3), 
\begin{eqnarray}
 && \sum_{\beta} b_{\beta}(b_{\beta}+1)
                \left<\left< \tau_{k}({\cal O}_{\mu})
                \tau_{1}({\cal O}_{\beta})\right>\right>_{0} 
                \left<\left< {\cal O}^{\beta} \tau_{l}({\cal O}_{\nu}) 
                \right>\right>_{0} 
                 \nonumber  \\
& = & - \sum_{\beta} \left\{ (k+b_{\mu}+b_{\beta}+1)
                \left<\left< \tau_{k}({\cal O}_{\mu})
                \tau_{1}({\cal O}_{\beta})\right>\right>_{0}
		- (k+b_{\mu})\left<\left< \tau_{k}({\cal O}_{\mu})
                \tau_{1}({\cal O}_{\beta})\right>\right>_{0}
			\right\}  \nonumber \\
&& \textrm{ \hspace{30pt}}
		\left\{ (l+b_{\nu}+1-b_{\beta}) 
                \left<\left< {\cal O}^{\beta} \tau_{l}({\cal O}_{\nu}) 
                \right>\right>_{0}  
		- (l+ b_{\nu} +1)
		\left<\left< {\cal O}^{\beta} \tau_{l}({\cal O}_{\nu}) 
                \right>\right>_{0}  
		\right\}
                 \nonumber  \\
& = &  - \sum_{\beta} \left\{ 
                \left<\left< {\cal X}  \tau_{k}({\cal O}_{\mu})
                \tau_{1}({\cal O}_{\beta})\right>\right>_{0}
		-\sum_{\sigma}{\cal C}_{\mu}^{\sigma}
		\left<\left< \tau_{k-1}({\cal O}_{\sigma})
                \tau_{1}({\cal O}_{\beta})\right>\right>_{0}
				 \right.
                 \nonumber  \\
&& \textrm{ \hspace{80pt}} \left.
		-\sum_{\sigma}{\cal C}_{\beta}^{\sigma}
		\left<\left< \tau_{k}({\cal O}_{\mu})
                {\cal O}_{\sigma} \right>\right>_{0}
		- (k+b_{\mu})\left<\left< \tau_{k}({\cal O}_{\mu})
                \tau_{1}({\cal O}_{\beta})\right>\right>_{0}
			\right\}  \nonumber \\
&& \textrm{ \hspace{30pt}}
		\left\{ 
                \left<\left< {\cal X} 
			{\cal O}^{\beta} \tau_{l}({\cal O}_{\nu}) 
                \right>\right>_{0}
		-\sum_{\rho}{\cal C}_{\nu}^{\rho}
		\left<\left<{\cal O}^{\beta}
		 \tau_{l-1}({\cal O}_{\rho})
		\right>\right>_{0}  
		-\delta_{l, 0} {\cal C}^{\beta}_{\nu} \right.
                 \nonumber  \\
&& \textrm{ \hspace{180pt}} \left.
		- (l+ b_{\nu} +1)
		\left<\left< {\cal O}^{\beta} \tau_{l}({\cal O}_{\nu}) 
                \right>\right>_{0}  
		\right\}.   \label{eqn:QuadForm:1}
\end{eqnarray}
Similarly, we have
\begin{eqnarray}
&&   \sum_{\beta} b_{\beta}(1-b_{\beta})
                \left<\left< \tau_{k}({\cal O}_{\mu})
                {\cal O}^{\beta}\right>\right>_{0} 
                \left<\left< \tau_{1}({\cal O}_{\beta}) 
                \tau_{l}({\cal O}_{\nu}) 
                \right>\right>_{0} 
                 \nonumber  \\
& = &  \sum_{\beta} \left\{ (k+b_{\mu}+1-b_{\beta})
                \left<\left< \tau_{k}({\cal O}_{\mu})
                {\cal O}^{\beta}\right>\right>_{0}
		- (k+b_{\mu})\left<\left< \tau_{k}({\cal O}_{\mu})
                {\cal O}^{\beta}\right>\right>_{0}
			\right\}  \nonumber \\
&& \textrm{ \hspace{30pt}}
		\left\{ (l+b_{\nu}+1+b_{\beta}) 
                \left<\left< \tau_{1}({\cal O}_{\beta}) 
			\tau_{l}({\cal O}_{\nu}) 
                \right>\right>_{0}  
		- (l+ b_{\nu} +1)
		\left<\left< \tau_{1}({\cal O}_{\beta}) 
		\tau_{l}({\cal O}_{\nu}) 
                \right>\right>_{0}  
		\right\}
                 \nonumber  \\
& = &   \sum_{\beta} \left\{ 
                \left<\left< {\cal X}  \tau_{k}({\cal O}_{\mu})
                {\cal O}^{\beta}\right>\right>_{0}
		-\sum_{\sigma}{\cal C}_{\mu}^{\sigma}
		\left<\left< \tau_{k-1}({\cal O}_{\sigma})
                {\cal O}^{\beta}\right>\right>_{0}
                 -\delta_{k, 0} {\cal C}_{\mu}^{\beta}  \right.
		\nonumber  \\
&& \textrm{ \hspace{180pt}} \left.
	- (k+b_{\mu})\left<\left< \tau_{k}({\cal O}_{\mu})
                {\cal O}^{\beta}\right>\right>_{0}
			\right\}  \nonumber \\
&& \textrm{ \hspace{30pt}}
		\left\{ 
                \left<\left< {\cal X} 
			\tau_{1}({\cal O}_{\beta}) \tau_{l}({\cal O}_{\nu}) 
                \right>\right>_{0}
		-\sum_{\rho}{\cal C}_{\nu}^{\rho}
		\left<\left<\tau_{1}({\cal O}_{\beta})
		 \tau_{l-1}({\cal O}_{\rho})
		\right>\right>_{0}  \right.
                 \nonumber  \\
&& \textrm{ \hspace{80pt}} \left.
		-\sum_{\rho} {\cal C}_{\beta}^{\rho}
		\left<\left< 
			{\cal O}_{\rho} \tau_{l}({\cal O}_{\nu}) 
                \right>\right>_{0}
		- (l+ b_{\nu} +1)
		\left<\left< \tau_{1}({\cal O}_{\beta}) 
			\tau_{l}({\cal O}_{\nu}) 
                \right>\right>_{0}  
		\right\}.   \label{eqn:QuadForm:2}
\end{eqnarray}
Expanding both (\ref{eqn:QuadForm:1}) and (\ref{eqn:QuadForm:2}),
summing them together, then using Lemma~\ref{lem:QuadRel} to simplify,
we obtain that, for $k \geq 0$ and $\l \geq -1$,
\begin{eqnarray*}
 && \sum_{\beta} b_{\beta}(b_{\beta}+1)
                \left<\left< \tau_{k}({\cal O}_{\mu})
                \tau_{1}({\cal O}_{\beta})\right>\right>_{0} 
                \left<\left< {\cal O}^{\beta} \tau_{l}({\cal O}_{\nu}) 
                \right>\right>_{0} 
                 \nonumber  \\
&&  + \sum_{\beta} b_{\beta}(1-b_{\beta})
                \left<\left< \tau_{k}({\cal O}_{\mu})
                {\cal O}^{\beta}\right>\right>_{0} 
                \left<\left< \tau_{1}({\cal O}_{\beta}) 
                \tau_{l}({\cal O}_{\nu}) 
                \right>\right>_{0} 
                 \nonumber  \\
&=& - \sum_{\sigma}{\cal C}_{\mu}^{\sigma}
	\left\{ \left<\left<  \tau_{k+1}({\cal O}_{\sigma})
                        {\cal X} 
                \tau_{l}({\cal O}_{\nu}) 
                \right>\right>_{0}
	- \delta_{k, 0}  \left<\left<  \tau_{1}({\cal O}_{\sigma})
                        {\cal X} 
                \tau_{l}({\cal O}_{\nu}) 
                \right>\right>_{0}  \right\}  \\
&& - (k + b_{\mu}) \left<\left<  \tau_{k+2}({\cal O}_{\mu})
                        {\cal X} 
                \tau_{l}({\cal O}_{\nu}) 
                \right>\right>_{0} \\
&& + \sum_{\sigma}{\cal C}_{\nu}^{\sigma}
	\left\{ \left<\left<  \tau_{k}({\cal O}_{\mu})
                        {\cal X} 
                \tau_{l+1}({\cal O}_{\sigma}) 
                \right>\right>_{0}
	- \delta_{l, 0}  \left<\left<  \tau_{k}({\cal O}_{\mu})
                        {\cal X} 
                \tau_{1}({\cal O}_{\sigma}) 
                \right>\right>_{0}  
	- \delta_{l, -1}  \left<\left<  \tau_{k}({\cal O}_{\mu})
                        {\cal X} 
                {\cal O}_{\sigma}
                \right>\right>_{0}  
	\right\}  \\
&& - \sum_{\sigma, \rho}{\cal C}_{\mu}^{\sigma}{\cal C}_{\nu}^{\rho}
	\left\{ \left<\left<  \tau_{k-1}({\cal O}_{\sigma})
                \tau_{l+1}({\cal O}_{\rho}) 
                \right>\right>_{0}
		- \left<\left<  \tau_{k+1}({\cal O}_{\sigma})
                \tau_{l-1}({\cal O}_{\rho}) 
                \right>\right>_{0}     \right.        \\
&& \textrm{ \hspace{20pt}}   \left.
		+ \delta_{k, 0}
		 \left<\left<  \tau_{1}({\cal O}_{\sigma})
                \tau_{l-1}({\cal O}_{\rho}) 
                \right>\right>_{0}
		- \delta_{l, 0}  \left<\left<  \tau_{k-1}({\cal O}_{\sigma})
                \tau_{1}({\cal O}_{\rho}) 
                \right>\right>_{0}  
		- \delta_{l, -1}  \left<\left<  \tau_{k-1}({\cal O}_{\sigma})
                {\cal O}_{\rho}
                \right>\right>_{0}  
	\right\}  \\
&& - \sum_{\sigma} (k+b_{\mu}){\cal C}_{\nu}^{\sigma}
	\left\{ \left<\left<  \tau_{k}({\cal O}_{\mu})
                \tau_{l+1}({\cal O}_{\sigma}) 
                \right>\right>_{0}
		- \left<\left<  \tau_{k+2}({\cal O}_{\mu})
                \tau_{l-1}({\cal O}_{\sigma}) 
                \right>\right>_{0}  \right.  \\
&& \textrm{ \hspace{100pt}} \left.
		- \delta_{l, 0}  \left<\left<  \tau_{k}({\cal O}_{\mu})
                \tau_{1}({\cal O}_{\sigma}) 
                \right>\right>_{0}  
		- \delta_{l, -1}  \left<\left<  \tau_{k}({\cal O}_{\mu})
                {\cal O}_{\sigma}
                \right>\right>_{0}  
	\right\}  \\
&& - \sum_{\sigma, \rho}
	{\cal C}_{\sigma}^{\rho}
	\left<\left<   {\cal X}\tau_{k}({\cal O}_{\mu}) 
                {\cal O}^{\sigma}
                \right>\right>_{0} 
	 \left<\left<   
                {\cal O}_{\rho} \tau_{l}({\cal O}_{\nu}) 
                \right>\right>_{0}    \\
&& + \sum_{\alpha, \sigma, \rho}
	{\cal C}_{\sigma}^{\rho}{\cal C}_{\mu}^{\alpha}
	\left<\left<   \tau_{k-1}({\cal O}_{\alpha}) 
                {\cal O}^{\sigma}
                \right>\right>_{0} 
	 \left<\left<   
                {\cal O}_{\rho} \tau_{l}({\cal O}_{\nu}) 
                \right>\right>_{0}    \\
&& + \sum_{\sigma, \rho}
	{\cal C}_{\sigma}^{\rho}
	\left<\left<  
		\tau_{k}({\cal O}_{\mu}) 
                {\cal O}_{\rho} 
                \right>\right>_{0} 
	 \left<\left<   
                {\cal O}^{\sigma}  {\cal X} \tau_{l}({\cal O}_{\nu}) 
                \right>\right>_{0}    \\
&& - \sum_{\alpha, \sigma, \rho}
	{\cal C}_{\sigma}^{\rho}{\cal C}_{\nu}^{\alpha}
	\left<\left<  
		\tau_{k}({\cal O}_{\mu}) 
                {\cal O}_{\rho} 
                \right>\right>_{0} 
	 \left<\left<   
                {\cal O}^{\sigma}  \tau_{l-1}({\cal O}_{\alpha}) 
                \right>\right>_{0}    \\
&& + \sum_{\sigma, \rho} (k+b_{\mu})
	{\cal C}_{\sigma}^{\rho}
	\left<\left<  
		\tau_{k}({\cal O}_{\mu}) 
                {\cal O}^{\sigma} 
                \right>\right>_{0} 
	 \left<\left<   
                {\cal O}_{\rho} \tau_{l}({\cal O}_{\nu}) 
                \right>\right>_{0}    \\
&&  + (l + b_{\nu} +1)
	\left\{ \left<\left<  \tau_{k}({\cal O}_{\mu})
                        {\cal X} 
                \tau_{l+2}({\cal O}_{\nu}) 
                \right>\right>_{0}
	- \delta_{l, -1}  \left<\left<  \tau_{k}({\cal O}_{\mu})
                        {\cal X} 
                \tau_{1}({\cal O}_{\nu}) 
                \right>\right>_{0}  \right\}  \\
&&  - \sum_{\sigma} (l+b_{\nu}+1){\cal C}_{\mu}^{\sigma}
	\left\{ \left<\left<  \tau_{k-1}({\cal O}_{\sigma})
                \tau_{l+2}({\cal O}_{\nu}) 
                \right>\right>_{0}
		- \left<\left<  \tau_{k+1}({\cal O}_{\sigma})
                \tau_{l}({\cal O}_{\nu}) 
                \right>\right>_{0} \right. \\
&& \textrm{ \hspace{100pt}} \left.
		+ \delta_{k, 0}  \left<\left<  \tau_{1}({\cal O}_{\sigma})
                \tau_{l}({\cal O}_{\nu}) 
                \right>\right>_{0}  
		- \delta_{l, -1}  \left<\left<  \tau_{k-1}({\cal O}_{\sigma})
                \tau_{1}({\cal O}_{\nu}) 
                \right>\right>_{0}  
		\right\}  \\
&&  - (l+b_{\nu}+1)(k+b_{\mu})
	\left\{ \left<\left<  \tau_{k}({\cal O}_{\mu})
                \tau_{l+2}({\cal O}_{\nu}) 
                \right>\right>_{0}
		- \left<\left<  \tau_{k+2}({\cal O}_{\mu})
                \tau_{l}({\cal O}_{\nu}) 
                \right>\right>_{0}  \right.\\
&& \textrm{ \hspace{120pt}} \left.
		- \delta_{l, -1}  \left<\left<  \tau_{k}({\cal O}_{\mu})
                \tau_{1}({\cal O}_{\nu}) 
                \right>\right>_{0}  
		\right\}  \\
&&  - (l+b_{\nu}+1)\sum_{\sigma, \rho} {\cal C}_{\sigma}^{\rho}
	\left<\left<  \tau_{k}({\cal O}_{\mu})
                {\cal O}_{\rho} 
                \right>\right>_{0}
		\left<\left< {\cal O}^{\sigma}
                \tau_{l}({\cal O}_{\nu}) 
                \right>\right>_{0}  \\
&& - \delta_{k, 0}  \sum_{\sigma} 
	b_{\sigma} {\cal C}_{\mu}^{\sigma}
		\left<\left<  \tau_{1}({\cal O}_{\sigma})
                \tau_{l}({\cal O}_{\nu}) 
                \right>\right>_{0}  
		  \\
&& + \delta_{l, 0}  \sum_{\sigma} (b_{\sigma}+1)
		{\cal C}_{\nu}^{\sigma}
		\left<\left<  \tau_{k}({\cal O}_{\mu})
                \tau_{1}({\cal O}_{\sigma}) 
                \right>\right>_{0}.   
\end{eqnarray*}
Using Lemma~\ref{lem:EulerCorr} to remove ${\cal X}$ in
this expression and simplifying, we obtain the desired 
formula.
$\Box$

Applying Lemma~\ref{lem:QuadForm} to the last 7 terms of 
(\ref{eqn:2WDVVright:2}) and simplifying, we obtain
\begin{eqnarray} \label{eqn:2WDVVright:3}
  &&  \sum_{\beta}
	\left<\left< {\cal L}_{1} \tau_{k}({\cal O}_{\mu})
                {\cal O}_{\beta} \right>\right>_{0}
	\left<\left< {\cal O}^{\beta} 
		({\cal L}_{0} - 2{\cal D})
                \tau_{l}({\cal O}_{\nu}) 
                \right>\right>_{0} \nonumber \\
&=&    (k+b_{\mu})(k+b_{\mu}+1) (k+b_{\mu}+2) 
   		 \left<\left<\tau_{k+2}({\cal O}_{\mu}) 
                \tau_{l}({\cal O}_{\nu}) 
                \right>\right>_{0} 
               \nonumber   \\
	&& + \sum_{\sigma} \{3(k+b_{\mu})^{2} + 6(k+b_{\mu})+ 2 \}
	{\cal C}_{\mu}^{\sigma}
            \left<\left<\tau_{k+1}({\cal O}_{\sigma}) 
                \tau_{l}({\cal O}_{\nu}) 
                \right>\right>_{0} 
           \nonumber  \\
	&& + \sum_{\sigma} 3(k+b_{\mu}+1)
                ({\cal C}^{2})_{\mu}^{\sigma}
            \left<\left<\tau_{k}({\cal O}_{\sigma}) 
                \tau_{l}({\cal O}_{\nu}) 
                \right>\right>_{0}  \nonumber  \\
	&& + \sum_{\sigma}({\cal C}^{3})_{\mu}^{\sigma} 
	 	\left<\left<\tau_{k-1}({\cal O}_{\sigma}) 
                \tau_{l}({\cal O}_{\nu}) 
                \right>\right>_{0}  \nonumber  \\
	&& + \sum_{\beta} (1-b_{\beta})b_{\beta}(1+b_{\beta})
 		\left<\left< \tau_{k}({\cal O}_{\mu})
                \tau_{1}({\cal O}_{\beta})\right>\right>_{0} 
                \left<\left< {\cal O}^{\beta}) \tau_{l}({\cal O}_{\nu}) 
                \right>\right>_{0} 
                 \nonumber  \\
        && + \sum_{\alpha, \beta} b_{\alpha}(1-b_{\alpha})
		(b_{\alpha}+b_{\beta})
        \left<\left< {\cal O}^{\alpha}\right>\right>_{0}         
	 \left<\left<  {\cal O}_{\alpha}
                {\cal O}_{\beta}\right>\right>_{0} 
        \left<\left< {\cal O}^{\beta} \tau_{k}({\cal O}_{\mu})
                \tau_{l}({\cal O}_{\nu}) 
                \right>\right>_{0} 
                \nonumber  \\
        && +    \sum_{\alpha, \beta}
                b_{\alpha}(1-b_{\alpha}){\cal C}_{\alpha \beta}
                \left<\left< {\cal O}^{\alpha}
                \right>\right>_{0} 
                \left<\left< {\cal O}^{\beta} 
                \tau_{k}({\cal O}_{\mu})
                \tau_{l}({\cal O}_{\nu}) 
                \right>\right>_{0} \nonumber  \\
	&& + \sum_{\alpha}
                b_{\alpha}(1-b_{\alpha})(b_{\alpha}+1)
		\left<\left< \tau_{k}({\cal O}_{\mu}) {\cal O}^{\alpha}
                \right>\right>_{0} 
                \left<\left< \tau_{1}({\cal O}_{\alpha} )
                \tau_{l}({\cal O}_{\nu}) 
                \right>\right>_{0} \nonumber  \\
	&& - \sum_{\alpha, \beta}
                (3b_{\alpha}^{2}-1) {\cal C}_{\alpha}^{\beta}
		\left<\left< \tau_{k}({\cal O}_{\mu}) {\cal O}_{\beta}
                \right>\right>_{0} 
                \left<\left< {\cal O}^{\alpha}
                \tau_{l}({\cal O}_{\nu}) 
                \right>\right>_{0} \nonumber  \\
	&& + \delta_{k, 0} \delta_{l, 0} ({\cal C}^{3})_{\mu \nu} 
	\nonumber  \\
	&& +  (l+b_{\nu})(l+b_{\nu}+1) (l+b_{\nu}+2) 
   		 \left<\left<\tau_{k}({\cal O}_{\mu}) 
                \tau_{l+2}({\cal O}_{\nu}) 
                \right>\right>_{0} 
               \nonumber   \\
	&& + \sum_{\sigma} \{3(l+b_{\nu})^{2} + 6(l+b_{\nu})+ 2 \}
		{\cal C}_{\nu}^{\sigma}
            \left<\left<\tau_{k}({\cal O}_{\mu}) 
                \tau_{l+1}({\cal O}_{\sigma}) 
                \right>\right>_{0} 
           \nonumber  \\
	&& + \sum_{\sigma} 3(l+b_{\nu}+1)
                ({\cal C}^{2})_{\nu}^{\sigma}
            \left<\left<\tau_{k}({\cal O}_{\mu}) 
                \tau_{l}({\cal O}_{\sigma}) 
                \right>\right>_{0}  \nonumber  \\
	&& + \sum_{\sigma}({\cal C}^{3})_{\nu}^{\sigma} 
	 	\left<\left<\tau_{k}({\cal O}_{\mu}) 
                \tau_{l-1}({\cal O}_{\sigma}) 
                \right>\right>_{0}. 
\end{eqnarray}
Now it is straightforward to check that the difference of the right
hand sides of (\ref{eqn:2WDVVright:3}) and (\ref{eqn:2WDVVleft})
is
$\frac{\partial}{\partial t^{\nu}_{l}}
        \frac{\partial}{\partial t^{\mu}_{k}} \Psi_{0, 2}$,
by using topological recursion relation and the relationships between
$b_{\alpha}$'s stated in section~\ref{sec:notation}.
The generalized WDVV equation implies that the left hand sides of
(\ref{eqn:2WDVVright:3}) and (\ref{eqn:2WDVVleft}) are equal.
Therefore 
$\frac{\partial}{\partial t^{\nu}_{l}}
        \frac{\partial}{\partial t^{\mu}_{k}} \Psi_{0, 2}=0$,
for arbitrary $(k, \mu)$ and $(l, \nu)$. As noted at the end of 
section~\ref{sec:VirOp},
this implies the following
\begin{pro}  \label{pro:g0L2}
	The genus-$0$ free energy function $F_{0}$ satisfies the $L_{2}$
	constraint.
\end{pro}

\section{$\widetilde{L}_{n}$ constraints for genus-$0$ free energy function}
\label{sec:tilde}

Besides $L_{n}$ constraints, Eguchi, Hori, and Xiong 
(\cite{EHX2})also conjectured 
the existence of another sequence of constraints for the free energy functions.
We call them $\widetilde{L}_{n}$ constraints. 
The $\widetilde{L}_{0}$ constraint is the dilaton equation.  
In this section we will prove
the $\widetilde{L}_{1}$ and $\widetilde{L}_{2}$ constraints for the genus-$0$
free energy function $F_{0}$.  

\subsection{$\widetilde{L}_{1}$ constraint}
\label{sec:tilde1}
The
 $\widetilde{L}_{1}$ constraint  predicts the vanishing of the
following function:
\[ \widetilde{\Psi}_{0, 1} := 
		- \sum_{m, \alpha} \tilde{t}^{\alpha}_{m}
		 \left<\left< \tau_{m+1}({\cal O}_{\alpha})
			\right>\right>_{0}
		+ \frac{1}{2} \sum_{\alpha}
			\left<\left< {\cal O}_{\alpha}
			\right>\right>_{0}
			\left<\left< {\cal O}^{\alpha}
			\right>\right>_{0}. \]
We now prove that this prediction is true.
\begin{pro}
\[ \widetilde{\Psi}_{0, 1} = 0. \]
\end{pro}
{\bf Proof}: As noted at the end of section~\ref{sec:VirOp},
we only need to show that all second derivatives of
$\widetilde{\Psi}_{0, 1}$ vanish. In fact, for any
$(k, \mu)$ and $(\l, \nu)$,
\begin{eqnarray*}
	\frac{\partial^{2}}{\partial t^{\mu}_{k} \partial t^{\nu}_{l}}
	\widetilde{\Psi}_{0, 1}   
& = &  	- \sum_{m, \alpha} \tilde{t}^{\alpha}_{m}
		 \left<\left< \tau_{m+1}({\cal O}_{\alpha})
			\tau_{k}({\cal O}_{\mu})
			\tau_{l}({\cal O}_{\nu})
			\right>\right>_{0}        \\
	&&	-  \left<\left< 
			\tau_{k+1}({\cal O}_{\mu})
			\tau_{l}({\cal O}_{\nu})
			\right>\right>_{0}
		-  \left<\left< 
			\tau_{k}({\cal O}_{\mu})
			\tau_{l+1}({\cal O}_{\nu})
			\right>\right>_{0}               \\
	&&	+ \sum_{\alpha}
			\left<\left< \tau_{k}({\cal O}_{\mu})
				\tau_{l}({\cal O}_{\nu})
				{\cal O}_{\alpha}
			\right>\right>_{0}
			\left<\left< {\cal O}^{\alpha}
			\right>\right>_{0}  \\
	&&	+ \sum_{\alpha}
			\left<\left< \tau_{k}({\cal O}_{\mu})
				{\cal O}_{\alpha}
			\right>\right>_{0}
			\left<\left< {\cal O}^{\alpha}
				\tau_{l}({\cal O}_{\nu})
			\right>\right>_{0}.
\end{eqnarray*} 
By formula (\ref{eqn:StringRec}), the second, third, and fifth terms
cancelled with each other. Hence, applying the topological recursion
relation to the first term,  we have
\begin{eqnarray*}
	\frac{\partial^{2}}{\partial t^{\mu}_{k} \partial t^{\nu}_{l}}
	\widetilde{\Psi}_{0, 1}   
& = &  	- \sum_{m, \alpha, \beta} \tilde{t}^{\alpha}_{m}
		 \left<\left< \tau_{m}({\cal O}_{\alpha})
			{\cal O}_{\beta}
			\right>\right>_{0}     
		 \left<\left< {\cal O}^{\beta}
			\tau_{k}({\cal O}_{\mu})
			\tau_{l}({\cal O}_{\nu})
			\right>\right>_{0}        \\
	&&	+ \sum_{\alpha}
			\left<\left< \tau_{k}({\cal O}_{\mu})
				\tau_{l}({\cal O}_{\nu})
				{\cal O}_{\alpha}
			\right>\right>_{0}
			\left<\left< {\cal O}^{\alpha}
			\right>\right>_{0}    \\
&=&   	\sum_{\beta} 
		 \left<\left< {\cal D}
			{\cal O}_{\beta}
			\right>\right>_{0}     
		 \left<\left< {\cal O}^{\beta}
			\tau_{k}({\cal O}_{\mu})
			\tau_{l}({\cal O}_{\nu})
			\right>\right>_{0}        \\
	&&	+ \sum_{\alpha}
			\left<\left< \tau_{k}({\cal O}_{\mu})
				\tau_{l}({\cal O}_{\nu})
				{\cal O}_{\alpha}
			\right>\right>_{0}
			\left<\left< {\cal O}^{\alpha}
			\right>\right>_{0}.    \\
\end{eqnarray*}
By Lemma~\ref{lem:DilatonCorr} (2), we have
\[	\frac{\partial^{2}}{\partial t^{\mu}_{k} \partial t^{\nu}_{l}}
	\widetilde{\Psi}_{0, 1}    = 0. \]
$\Box$

\subsection{$\widetilde{L}_{2}$ constraint}
\label{sec:tilde2}

The $\widetilde{L}_{2}$ constraint predicts the vanishing of the function
\begin{eqnarray*}
	\widetilde{\Psi}_{0, 2} 
& := & 	\sum_{m, \alpha} (m+b_{\alpha}+1) \tilde{t}^{\alpha}_{m}
		\left<\left< \tau_{m+2}({\cal O}_{\alpha})
		 \right>\right>_{0}
	+ \sum_{m, \alpha, \beta} {\cal C}_{\alpha}^{\beta}
		\tilde{t}^{\alpha}_{m}
		\left<\left< \tau_{m+1}({\cal O}_{\beta})
		 \right>\right>_{0}  \\
&&	- \sum_{\alpha} b_{\alpha} 
		\left<\left< {\cal O}^{\alpha}
		 \right>\right>_{0}
		\left<\left< \tau_{1}({\cal O}_{\alpha})
		 \right>\right>_{0}
	- \frac{1}{2} \sum_{\alpha, \beta} {\cal C}_{\alpha}^{\beta}
		\left<\left< {\cal O}^{\alpha}
		 \right>\right>_{0}
		\left<\left< {\cal O}_{\beta}
		 \right>\right>_{0}.
\end{eqnarray*}
To prove this constraint, we need to study correlation functions
involving the following vector field,
\[ \widetilde{{\cal L}}_{1} 
	:= \sum_{m, \alpha} \tilde{t}^{\alpha}_{m}
				\frac{\partial}{\partial t^{\alpha}_{m+1}}.
			\]
The $\widetilde{L}_{1}$ constraint can be reformulated as
\begin{equation} \label{eqn:1pttilde1}
	\left<\left< \widetilde{{\cal L}}_{1} 
		 \right>\right>_{0} =
	\frac{1}{2} \sum_{\alpha}
                        \left<\left< {\cal O}_{\alpha}
                        \right>\right>_{0}
                        \left<\left< {\cal O}^{\alpha}
                        \right>\right>_{0}.
\end{equation}
Using the fact that
\[ [ \widetilde{{\cal L}}_{1} , \, \frac{\partial}{\partial t^{\alpha}_{m}}]
	= - \frac{\partial}{\partial t^{\alpha}_{m+ 1}}, 
  \]
we can show the following
\begin{lem} \label{lem:tilde1Corr}
\begin{eqnarray*}
&(1)& \left<\left< \widetilde{{\cal L}}_{1} \tau_{m}({\cal O}_{\alpha}) 
		 \right>\right>_{0} =
		- \left<\left< \tau_{m+1}({\cal O}_{\alpha}) 
                        \right>\right>_{0} +
		\sum_{\sigma}
                        \left<\left< \tau_{m}({\cal O}_{\alpha}) 
				{\cal O}_{\sigma}
                        \right>\right>_{0}
                        \left<\left< {\cal O}^{\sigma}
                        \right>\right>_{0}, \\
&(2)& \left<\left< \widetilde{{\cal L}}_{1} 
		\tau_{m}({\cal O}_{\alpha}) 
		\tau_{n}({\cal O}_{\beta}) 
		 \right>\right>_{0} =
		\sum_{\sigma}
                        \left<\left< \tau_{m}({\cal O}_{\alpha}) 
				\tau_{n}({\cal O}_{\beta}) 
				{\cal O}_{\sigma}
                        \right>\right>_{0}
                        \left<\left< {\cal O}^{\sigma}
                        \right>\right>_{0}.
\end{eqnarray*}
\end{lem}
{\bf Proof}:
	The first equation follows directly from (\ref{eqn:1pttilde1})
and the fact that
\begin{eqnarray*}
 \left<\left< \widetilde{{\cal L}}_{1}
	 \tau_{m}({\cal O}_{\alpha}) 
		 \right>\right>_{0}
& = & \widetilde{{\cal L}}_{1} \frac{\partial}{\partial t^{\alpha}_{m}} F_{0}
\, \,  = \, \, \frac{\partial}{\partial t^{\alpha}_{m}} 
	(\widetilde{{\cal L}}_{1} F_{0})
	- \frac{\partial}{\partial t^{\alpha}_{m+ 1}} F_{0}.
\end{eqnarray*}

Now we prove the second equation.
\begin{eqnarray*}
\left<\left< \widetilde{{\cal L}}_{1} 
		\tau_{m}({\cal O}_{\alpha}) 
		\tau_{n}({\cal O}_{\beta}) 
		 \right>\right>_{0} 
&=& \widetilde{{\cal L}}_{1} 
	\frac{\partial}{\partial t^{\alpha}_{m}}
	\frac{\partial}{\partial t^{\beta}_{n}} F_{0}  \\
&= & \left\{ \frac{\partial}{\partial t^{\alpha}_{m}}
	\frac{\partial}{\partial t^{\beta}_{n}} 
 	\widetilde{{\cal L}}_{1} 
	- \frac{\partial}{\partial t^{\alpha}_{m+1}}
	\frac{\partial}{\partial t^{\beta}_{n}} 
	- \frac{\partial}{\partial t^{\alpha}_{m}}
	\frac{\partial}{\partial t^{\beta}_{n+1}} \right\} F_{0}.
\end{eqnarray*}
By (\ref{eqn:1pttilde1}), we have
\begin{eqnarray*}
\left<\left< \widetilde{{\cal L}}_{1} 
		\tau_{m}({\cal O}_{\alpha}) 
		\tau_{n}({\cal O}_{\beta}) 
		 \right>\right>_{0} 
&=& \sum_{\sigma}
                \left<\left< \tau_{m}({\cal O}_{\alpha}) 
				\tau_{n}({\cal O}_{\beta}) 
				{\cal O}_{\sigma}
                        \right>\right>_{0}
                        \left<\left< {\cal O}^{\sigma}
                        \right>\right>_{0}  \\
&& + 	\sum_{\sigma}
                        \left<\left< \tau_{m}({\cal O}_{\alpha}) 
				{\cal O}_{\sigma}
                        \right>\right>_{0}
                        \left<\left< {\cal O}^{\sigma}
				\tau_{n}({\cal O}_{\beta}) 
                        \right>\right>_{0} \\
&& -  \left<\left< \tau_{m+1}({\cal O}_{\alpha}) 
				\tau_{n}({\cal O}_{\beta}) 
                        \right>\right>_{0}
	-  \left<\left< \tau_{m}({\cal O}_{\alpha}) 
				\tau_{n+1}({\cal O}_{\beta}) 
                        \right>\right>_{0}.
\end{eqnarray*}
By formula (\ref{eqn:StringRec}), the last three terms are canceled
with each other. This proves the second equation.
$\Box$

We also need the following
\begin{lem} \label{lem:QuadForm:2}
\begin{eqnarray*}
&&	\sum_{\sigma} b_{\sigma}
                      \left<\left< \tau_{m}({\cal O}_{\alpha}) 
				{\cal O}^{\sigma}
                        \right>\right>_{0}
                        \left<\left< \tau_{1}({\cal O}_{\sigma})
				\tau_{n}({\cal O}_{\beta}) 
                        \right>\right>_{0} \\
&& + \sum_{\sigma} b_{\sigma}
                      \left<\left< \tau_{m}({\cal O}_{\alpha}) 
				\tau_{1}({\cal O}_{\sigma})
                        \right>\right>_{0}
                        \left<\left< {\cal O}^{\sigma}
				\tau_{n}({\cal O}_{\beta}) 
                        \right>\right>_{0} \\
& = & (m+b_{\alpha}+1) \left<\left< \tau_{m+2}({\cal O}_{\alpha}) 
				\tau_{n}({\cal O}_{\beta}) 
                        \right>\right>_{0}
	+ \sum_{\sigma}{\cal C}_{\alpha}^{\sigma}
		\left<\left< \tau_{m+1}({\cal O}_{\sigma}) 
				\tau_{n}({\cal O}_{\beta}) 
                        \right>\right>_{0}   \\
&& + (n+b_{\beta}+1) \left<\left< \tau_{m}({\cal O}_{\alpha}) 
				\tau_{n+2}({\cal O}_{\beta}) 
                        \right>\right>_{0}
	+ \sum_{\sigma}{\cal C}_{\beta}^{\sigma}
		\left<\left< \tau_{m}({\cal O}_{\alpha}) 
				\tau_{n+1}({\cal O}_{\sigma}) 
                        \right>\right>_{0}   \\
&& - \sum_{\sigma, \rho} {\cal C}_{\sigma}^{\rho}
                      \left<\left< \tau_{m}({\cal O}_{\alpha}) 
				{\cal O}^{\sigma}
                        \right>\right>_{0}
                        \left<\left< {\cal O}_{\rho}
				\tau_{n}({\cal O}_{\beta}) 
                        \right>\right>_{0}. 
\end{eqnarray*}
\end{lem}
{\bf Proof}: By Lemma~\ref{lem:EulerCorr} (3), we have
\begin{eqnarray} 
&&	\sum_{\sigma} b_{\sigma}
                      \left<\left< \tau_{m}({\cal O}_{\alpha}) 
				{\cal O}^{\sigma}
                        \right>\right>_{0}
                        \left<\left< \tau_{1}({\cal O}_{\sigma})
				\tau_{n}({\cal O}_{\beta}) 
                        \right>\right>_{0} \nonumber \\
& = & 	\sum_{\sigma} 
                      \left<\left< \tau_{m}({\cal O}_{\alpha}) 
				{\cal O}^{\sigma}
                        \right>\right>_{0} (n+b_{\beta}+b_{\sigma}+1)
                        \left<\left< \tau_{1}({\cal O}_{\sigma})
				\tau_{n}({\cal O}_{\beta}) 
                        \right>\right>_{0} \nonumber \\
&&	-(n+b_{\beta}+1)\sum_{\sigma}
                      \left<\left< \tau_{m}({\cal O}_{\alpha}) 
				{\cal O}^{\sigma}
                        \right>\right>_{0}
                        \left<\left< \tau_{1}({\cal O}_{\sigma})
				\tau_{n}({\cal O}_{\beta}) 
                        \right>\right>_{0}  \nonumber \\
&=& 	\sum_{\sigma} 
                      \left<\left< \tau_{m}({\cal O}_{\alpha}) 
				{\cal O}^{\sigma}
                        \right>\right>_{0} 
	\left\{ \left<\left< {\cal X}\tau_{1}({\cal O}_{\sigma})
				\tau_{n}({\cal O}_{\beta}) 
                        \right>\right>_{0} 
		- \sum_{\rho}{\cal C}_{\sigma}^{\rho}
		\left<\left< {\cal O}_{\rho}
				\tau_{n}({\cal O}_{\beta}) 
                        \right>\right>_{0}      \right.      \nonumber \\
&& \textrm{ \hspace{150pt}} \left.
		- \sum_{\rho}{\cal C}_{\beta}^{\rho}
		\left<\left< \tau_{1}({\cal O}_{\sigma})
				\tau_{n-1}({\cal O}_{\rho}) 
                        \right>\right>_{0} 
		\right\}
	\nonumber \\
&&	-(n+b_{\beta}+1)\sum_{\sigma}
                      \left<\left< \tau_{m}({\cal O}_{\alpha}) 
				{\cal O}^{\sigma}
                        \right>\right>_{0}
                        \left<\left< \tau_{1}({\cal O}_{\sigma})
				\tau_{n}({\cal O}_{\beta}) 
                        \right>\right>_{0}  .  \label{eqn:QuadForm2:1}
\end{eqnarray}
On the other hand,
\begin{eqnarray} 
&&	\sum_{\sigma} b_{\sigma}
                      \left<\left< \tau_{m}({\cal O}_{\alpha}) 
				\tau_{1}({\cal O}_{\sigma})
                        \right>\right>_{0}
                        \left<\left< {\cal O}^{\sigma}
				\tau_{n}({\cal O}_{\beta}) 
                        \right>\right>_{0} 
		\nonumber \\
& = & 	- \sum_{\sigma} 
                       \left<\left< 
				\tau_{m}({\cal O}_{\alpha})
				\tau_{1}({\cal O}_{\sigma}) 
                        \right>\right>_{0} 
		(n+b_{\beta}+1-b_{\sigma})
			 \left<\left< {\cal O}^{\sigma}
				\tau_{n}({\cal O}_{\beta}) 
                        \right>\right>_{0} 
	\nonumber \\
&&	+(n+b_{\beta}+1)\sum_{\sigma}
                       \left<\left< 
				\tau_{m}({\cal O}_{\alpha})
				\tau_{1}({\cal O}_{\sigma}) 
                        \right>\right>_{0} 
                        \left<\left< {\cal O}^{\sigma}
				\tau_{n}({\cal O}_{\beta}) 
                        \right>\right>_{0}  \nonumber \\
&=& 	-\sum_{\sigma} \left<\left< 
				\tau_{m}({\cal O}_{\alpha})
				\tau_{1}({\cal O}_{\sigma}) 
                        \right>\right>_{0} 
	\left\{ \left<\left< {\cal X}{\cal O}^{\sigma}
				\tau_{n}({\cal O}_{\beta}) 
                        \right>\right>_{0} 
		- \delta_{n, 0}{\cal C}_{\beta}^{\sigma}    
		- \sum_{\rho}{\cal C}_{\beta}^{\rho}
		\left<\left< {\cal O}^{\sigma}
				\tau_{n-1}({\cal O}_{\rho}) 
                        \right>\right>_{0} 
		\right\}
	\nonumber \\
&&	+(n+b_{\beta}+1)\sum_{\sigma}
                       \left<\left< 
				\tau_{m}({\cal O}_{\alpha})
				\tau_{1}({\cal O}_{\sigma}) 
                        \right>\right>_{0} 
                        \left<\left< {\cal O}^{\sigma}
				\tau_{n}({\cal O}_{\beta}) 
                        \right>\right>_{0}. \label{eqn:QuadForm2:2}
\end{eqnarray}
Summing up (\ref{eqn:QuadForm2:1}) and (\ref{eqn:QuadForm2:2}) together,
then using Lemma~\ref{lem:QuadRel} (ii) and (iii) to simplify it, we obtain
\begin{eqnarray*}
&&	\sum_{\sigma} b_{\sigma}
                      \left<\left< \tau_{m}({\cal O}_{\alpha}) 
				{\cal O}^{\sigma}
                        \right>\right>_{0}
                        \left<\left< \tau_{1}({\cal O}_{\sigma})
				\tau_{n}({\cal O}_{\beta}) 
                        \right>\right>_{0} \\
&& + \sum_{\sigma} b_{\sigma}
                      \left<\left< \tau_{m}({\cal O}_{\alpha}) 
				\tau_{1}({\cal O}_{\sigma})
                        \right>\right>_{0}
                        \left<\left< {\cal O}^{\sigma}
				\tau_{n}({\cal O}_{\beta}) 
                        \right>\right>_{0} \\
&=&     \left<\left< \tau_{m+2}({\cal O}_{\alpha}) 
			{\cal X}
			\tau_{n}({\cal O}_{\beta}) 
                        \right>\right>_{0} \\
&& - \sum_{\sigma, \rho} {\cal C}_{\sigma}^{\rho}
                      \left<\left< \tau_{m}({\cal O}_{\alpha}) 
				{\cal O}^{\sigma}
                        \right>\right>_{0}
                        \left<\left< {\cal O}_{\rho}
				\tau_{n}({\cal O}_{\beta}) 
                        \right>\right>_{0} \\
&& - \sum_{\sigma} {\cal C}_{\beta}^{\sigma}
	\left\{ \left<\left< \tau_{m+2}({\cal O}_{\alpha}) 
			\tau_{n-1}({\cal O}_{\sigma}) 
                        \right>\right>_{0} 
		-\left<\left< \tau_{m}({\cal O}_{\alpha}) 
			\tau_{n+1}({\cal O}_{\sigma}) 
                        \right>\right>_{0} 
		+\delta_{n, 0} \left<\left< \tau_{m}({\cal O}_{\alpha}) 
			\tau_{1}({\cal O}_{\sigma}) 
                        \right>\right>_{0} 
	\right\} \\
&& + \delta_{n, o} \sum_{\sigma} {\cal C}_{\beta}^{\sigma}
		\left<\left< \tau_{m}({\cal O}_{\alpha}) 
			\tau_{1}({\cal O}_{\sigma}) 
                        \right>\right>_{0} \\
&& - (n + b_{\beta}+1) \left\{ \left<\left< \tau_{m+2}({\cal O}_{\alpha}) 
			\tau_{n}({\cal O}_{\beta}) 
                        \right>\right>_{0} 
			- \left<\left< \tau_{m}({\cal O}_{\alpha}) 
			\tau_{n+2}({\cal O}_{\beta}) 
                        \right>\right>_{0} \right\}.
\end{eqnarray*}
Using Lemma~\ref{lem:EulerCorr} (3) to remove ${\cal X}$ in the first
term and simplifying, we obtain the desired equation.
$\Box$

Now we are ready to prove the $\widetilde{L}_{2}$ constraint.
\begin{pro}
\[ \widetilde{\Psi}_{0, 2} = 0. \] 
\end{pro}
{\bf Proof}:
As in the proof of $\widetilde{L}_{1}$ constraint, we only need to show
that all second derivatives of $\widetilde{\Psi}_{0, 2}$ are equal to
zero. For arbitrary $(\mu, k)$ and $(\nu, l)$, we have
\begin{eqnarray*}
\frac{\partial^{2}}{\partial t^{\mu}_{k} \partial t^{\nu}_{l}}
        \widetilde{\Psi}_{0, 2} 
&=& \sum_{m, \alpha} (m+b_{\alpha}+1) \tilde{t}^{\alpha}_{m}
                \left<\left< \tau_{m+2}({\cal O}_{\alpha})
		        \tau_{k}({\cal O}_{\mu}) 
			\tau_{l}({\cal O}_{\nu}) 
                 \right>\right>_{0}   \\
 &&       + \sum_{m, \alpha, \beta} {\cal C}_{\alpha}^{\beta}
                \tilde{t}^{\alpha}_{m}
                \left<\left< \tau_{m+1}({\cal O}_{\beta})
			 \tau_{k}({\cal O}_{\mu}) 
			\tau_{l}({\cal O}_{\nu}) 
                 \right>\right>_{0}  \\
&& + (k+b_{\mu}+1) 
                \left<\left< 
		        \tau_{k+2}({\cal O}_{\mu}) 
			\tau_{l}({\cal O}_{\nu}) 
                 \right>\right>_{0}   \\
&& + \sum_{\sigma} {\cal C}_{\mu}^{\sigma}
                \left<\left< 
			 \tau_{k+1}({\cal O}_{\sigma}) 
			\tau_{l}({\cal O}_{\nu}) 
                 \right>\right>_{0}  \\
&& + (l+b_{\nu}+1) 
                \left<\left< 
		        \tau_{k}({\cal O}_{\mu}) 
			\tau_{l+2}({\cal O}_{\nu}) 
                 \right>\right>_{0}   \\
&& + \sum_{\sigma} {\cal C}_{\nu}^{\sigma}
                \left<\left< 
			 \tau_{k}({\cal O}_{\mu}) 
			\tau_{l+1}({\cal O}_{\sigma}) 
                 \right>\right>_{0}  \\
&& - \sum_{\sigma} b_{\sigma}
                      \left<\left< \tau_{k}({\cal O}_{\mu}) 
			\tau_{l}({\cal O}_{\nu}) 
                                {\cal O}^{\sigma}
                        \right>\right>_{0}
                        \left<\left< \tau_{1}({\cal O}_{\sigma})
                        \right>\right>_{0} \\
&& - \sum_{\sigma} b_{\sigma}
                      \left<\left< \tau_{k}({\cal O}_{\mu}) 
                                {\cal O}^{\sigma}
                        \right>\right>_{0}
                        \left<\left< \tau_{1}({\cal O}_{\sigma})
                                \tau_{l}({\cal O}_{\nu}) 
                        \right>\right>_{0} \\
&& - \sum_{\sigma} b_{\sigma}
                      \left<\left< \tau_{k}({\cal O}_{\mu}) 
                                \tau_{1}({\cal O}_{\sigma})
                        \right>\right>_{0}
                        \left<\left< {\cal O}^{\sigma}
                                \tau_{l}({\cal O}_{\nu}) 
                        \right>\right>_{0} \\
&& - \sum_{\sigma} b_{\sigma}
                      \left<\left< \tau_{k}({\cal O}_{\mu}) 
			\tau_{l}({\cal O}_{\nu}) 
                        \tau_{1}({\cal O}_{\sigma})
                        \right>\right>_{0}
                        \left<\left< {\cal O}^{\sigma}
                        \right>\right>_{0} \\
&& - \sum_{\sigma, \rho} {\cal C}_{\sigma}^{\rho}
                      \left<\left< \tau_{k}({\cal O}_{\mu}) 
			\tau_{l}({\cal O}_{\nu}) 
                                {\cal O}^{\sigma}
                        \right>\right>_{0}
                        \left<\left< {\cal O}_{\rho}
                        \right>\right>_{0} \\
&& - \sum_{\sigma, \rho} {\cal C}_{\sigma}^{\rho}
                      \left<\left< \tau_{k}({\cal O}_{\mu}) 
			     {\cal O}^{\sigma}
                        \right>\right>_{0}
                        \left<\left< {\cal O}_{\rho}
			\tau_{l}({\cal O}_{\nu}) 
                        \right>\right>_{0} .
\end{eqnarray*}
Applying the topological recursion relation to the first two terms, and
Lemma~\ref{lem:QuadForm:2} to the eighth and nineth terms, we have
\begin{eqnarray} 
\frac{\partial^{2}}{\partial t^{\mu}_{k} \partial t^{\nu}_{l}}
        \widetilde{\Psi}_{0, 2} 
&=& \sum_{m, \alpha, \sigma} (m+b_{\alpha}+1) \tilde{t}^{\alpha}_{m}
                \left<\left< \tau_{m+1}({\cal O}_{\alpha})
		        {\cal O}_{\sigma}
                 \right>\right>_{0}  
		\left<\left< {\cal O}^{\sigma}
		        \tau_{k}({\cal O}_{\mu}) 
			\tau_{l}({\cal O}_{\nu}) 
                 \right>\right>_{0}   \nonumber \\
 &&       + \sum_{m, \alpha, \beta, \sigma} {\cal C}_{\alpha}^{\beta}
                \tilde{t}^{\alpha}_{m}
                \left<\left< \tau_{m}({\cal O}_{\beta})
			{\cal O}_{\sigma} 
                 \right>\right>_{0}  
		\left<\left< {\cal O}^{\sigma}
			 \tau_{k}({\cal O}_{\mu}) 
			\tau_{l}({\cal O}_{\nu}) 
                 \right>\right>_{0}  \nonumber \\
&& - \sum_{\sigma} b_{\sigma}
                      \left<\left< \tau_{k}({\cal O}_{\mu}) 
			\tau_{l}({\cal O}_{\nu}) 
                                {\cal O}^{\sigma}
                        \right>\right>_{0}
                        \left<\left< \tau_{1}({\cal O}_{\sigma})
                        \right>\right>_{0} \nonumber \\
&& - \sum_{\sigma} b_{\sigma}
                      \left<\left< \tau_{k}({\cal O}_{\mu}) 
			\tau_{l}({\cal O}_{\nu}) 
                        \tau_{1}({\cal O}_{\sigma})
                        \right>\right>_{0}
                        \left<\left< {\cal O}^{\sigma}
                        \right>\right>_{0} \nonumber \\
&& - \sum_{\sigma, \rho} {\cal C}_{\sigma}^{\rho}
                      \left<\left< \tau_{k}({\cal O}_{\mu}) 
			\tau_{l}({\cal O}_{\nu}) 
                                {\cal O}^{\sigma}
                        \right>\right>_{0}
                        \left<\left< {\cal O}_{\rho}
                        \right>\right>_{0}.   \label{eqn:tildederiv}
\end{eqnarray}
We now use Lemma~\ref{lem:EulerCorr} (3) to compute the first two terms.
Let
\begin{eqnarray*}
f & := & \sum_{m, \alpha, \sigma} (m+b_{\alpha}+1) \tilde{t}^{\alpha}_{m}
                \left<\left< \tau_{m+1}({\cal O}_{\alpha})
		        {\cal O}_{\sigma}
                 \right>\right>_{0}  
		\left<\left< {\cal O}^{\sigma}
		        \tau_{k}({\cal O}_{\mu}) 
			\tau_{l}({\cal O}_{\nu}) 
                 \right>\right>_{0}   \nonumber \\
 &&       + \sum_{m, \alpha, \beta, \sigma} {\cal C}_{\alpha}^{\beta}
                \tilde{t}^{\alpha}_{m}
                \left<\left< \tau_{m}({\cal O}_{\beta})
			{\cal O}_{\sigma} 
                 \right>\right>_{0}  
		\left<\left< {\cal O}^{\sigma}
			 \tau_{k}({\cal O}_{\mu}) 
			\tau_{l}({\cal O}_{\nu}) 
                 \right>\right>_{0}  \nonumber \\
& = & \sum_{m, \alpha, \sigma}  \tilde{t}^{\alpha}_{m}
                \left<\left< {\cal X} \tau_{m+1}({\cal O}_{\alpha})
		        {\cal O}_{\sigma}
                 \right>\right>_{0}  
		\left<\left< {\cal O}^{\sigma}
		        \tau_{k}({\cal O}_{\mu}) 
			\tau_{l}({\cal O}_{\nu}) 
                 \right>\right>_{0}   \nonumber \\
&& - \sum_{m, \alpha, \sigma} b_{\sigma} \tilde{t}^{\alpha}_{m}
                \left<\left< \tau_{m+1}({\cal O}_{\alpha})
		        {\cal O}_{\sigma}
                 \right>\right>_{0}  
		\left<\left< {\cal O}^{\sigma}
		        \tau_{k}({\cal O}_{\mu}) 
			\tau_{l}({\cal O}_{\nu}) 
                 \right>\right>_{0}   \nonumber \\
& = & \sum_{\sigma} 
                \left<\left< {\cal X} \widetilde{{\cal L}}_{1}
		        {\cal O}_{\sigma}
                 \right>\right>_{0}  
		\left<\left< {\cal O}^{\sigma}
		        \tau_{k}({\cal O}_{\mu}) 
			\tau_{l}({\cal O}_{\nu}) 
                 \right>\right>_{0}   \nonumber \\
&& - \sum_{\sigma} b_{\sigma} 
                \left<\left< \widetilde{{\cal L}}_{1}
		        {\cal O}_{\sigma}
                 \right>\right>_{0}  
		\left<\left< {\cal O}^{\sigma}
		        \tau_{k}({\cal O}_{\mu}) 
			\tau_{l}({\cal O}_{\nu}) 
                 \right>\right>_{0}.  
\end{eqnarray*}
Using the generalized WDVV equation and Lemma~\ref{lem:tilde1Corr}, we have
\begin{eqnarray*}
f &=& \sum_{\sigma} 
                \left<\left< {\cal X} \tau_{k}({\cal O}_{\mu}) 
		        {\cal O}^{\sigma}
                 \right>\right>_{0}  
		\left<\left< {\cal O}_{\sigma}
			\widetilde{{\cal L}}_{1}
			\tau_{l}({\cal O}_{\nu}) 
                 \right>\right>_{0}   \nonumber \\
&& - \sum_{\sigma} b_{\sigma} 
                \left<\left< \widetilde{{\cal L}}_{1}
		        {\cal O}_{\sigma}
                 \right>\right>_{0}  
		\left<\left< {\cal O}^{\sigma}
		        \tau_{k}({\cal O}_{\mu}) 
			\tau_{l}({\cal O}_{\nu}) 
                 \right>\right>_{0}  \\
& = & \sum_{\sigma, \rho} 
                \left<\left< {\cal X} \tau_{k}({\cal O}_{\mu}) 
		        {\cal O}^{\sigma}
                 \right>\right>_{0} 
		 \left<\left< {\cal O}_{\sigma}
                                \tau_{l}({\cal O}_{\nu}) 
                                {\cal O}_{\rho}
                        \right>\right>_{0}
                        \left<\left< {\cal O}^{\rho}
                        \right>\right>_{0} \\
&& - \sum_{\sigma} b_{\sigma} 
		\left\{- \left<\left< \tau_{1}({\cal O}_{\sigma}) 
                        \right>\right>_{0} +
                \sum_{\rho}
                        \left<\left< {\cal O}_{\sigma} 
                                {\cal O}_{\rho}
                        \right>\right>_{0}
                        \left<\left< {\cal O}^{\rho}
                        \right>\right>_{0}
		\right\}
       		\left<\left< {\cal O}^{\sigma}
		        \tau_{k}({\cal O}_{\mu}) 
			\tau_{l}({\cal O}_{\nu}) 
                 \right>\right>_{0}.
\end{eqnarray*}
Using the generalized WDVV first and then applying Lemma~\ref{lem:EulerCorr}
(3), we have
\begin{eqnarray*}
f &=& \sum_{\sigma, \rho} 
                \left<\left< {\cal X} {\cal O}_{\rho}
		        {\cal O}_{\sigma}
                 \right>\right>_{0} 
		 \left<\left< {\cal O}^{\sigma}
                                \tau_{l}({\cal O}_{\nu}) 
				\tau_{k}({\cal O}_{\mu})                       
                        \right>\right>_{0}
                        \left<\left< {\cal O}^{\rho}
                        \right>\right>_{0} \\
&& + \sum_{\sigma} b_{\sigma} 
	\left<\left< \tau_{1}({\cal O}_{\sigma}) 
                        \right>\right>_{0}
	       		\left<\left< {\cal O}^{\sigma}
		        \tau_{k}({\cal O}_{\mu}) 
			\tau_{l}({\cal O}_{\nu}) 
                 \right>\right>_{0}  \\
&&  - \sum_{\sigma, \rho} b_{\sigma} 
                        \left<\left< {\cal O}_{\sigma} 
                                {\cal O}_{\rho}
                        \right>\right>_{0}
                        \left<\left< {\cal O}^{\rho}
                        \right>\right>_{0}
       		\left<\left< {\cal O}^{\sigma}
		        \tau_{k}({\cal O}_{\mu}) 
			\tau_{l}({\cal O}_{\nu}) 
                 \right>\right>_{0}  \\
& = & \sum_{\sigma, \rho} 
                \left\{ (b_{\rho} + b_{\sigma})
			\left<\left< {\cal O}_{\rho}
		       	 {\cal O}_{\sigma}
                	 \right>\right>_{0} 
			+ {\cal C}_{\sigma \rho}
		\right\}
		 \left<\left< {\cal O}^{\sigma}
                                \tau_{l}({\cal O}_{\nu}) 
				\tau_{k}({\cal O}_{\mu})                       
                        \right>\right>_{0}
                        \left<\left< {\cal O}^{\rho}
                        \right>\right>_{0} \\
&& + \sum_{\sigma} b_{\sigma} 
	\left<\left< \tau_{1}({\cal O}_{\sigma}) 
                        \right>\right>_{0}
	       		\left<\left< {\cal O}^{\sigma}
		        \tau_{k}({\cal O}_{\mu}) 
			\tau_{l}({\cal O}_{\nu}) 
                 \right>\right>_{0}  \\
&&  - \sum_{\sigma} b_{\sigma} 
                \sum_{\rho}
                        \left<\left< {\cal O}_{\sigma} 
                                {\cal O}_{\rho}
                        \right>\right>_{0}
                        \left<\left< {\cal O}^{\rho}
                        \right>\right>_{0}
       		\left<\left< {\cal O}^{\sigma}
		        \tau_{k}({\cal O}_{\mu}) 
			\tau_{l}({\cal O}_{\nu}) 
                 \right>\right>_{0}  \\
& = & \sum_{\sigma, \rho} 
                 {\cal C}_{\sigma \rho}
		\left<\left< {\cal O}^{\sigma}
                                \tau_{l}({\cal O}_{\nu}) 
				\tau_{k}({\cal O}_{\mu})                       
                        \right>\right>_{0}
                        \left<\left< {\cal O}^{\rho}
                        \right>\right>_{0} \\
&& + \sum_{\sigma} b_{\sigma} 
	\left<\left< \tau_{1}({\cal O}_{\sigma}) 
                        \right>\right>_{0}
	       		\left<\left< {\cal O}^{\sigma}
		        \tau_{k}({\cal O}_{\mu}) 
			\tau_{l}({\cal O}_{\nu}) 
                 \right>\right>_{0}  \\
&& + \sum_{\sigma} b_{\rho} 
                \sum_{\rho}
                        \left<\left< {\cal O}_{\sigma} 
                                {\cal O}_{\rho}
                        \right>\right>_{0}
                        \left<\left< {\cal O}^{\rho}
                        \right>\right>_{0}
       		\left<\left< {\cal O}^{\sigma}
		        \tau_{k}({\cal O}_{\mu}) 
			\tau_{l}({\cal O}_{\nu}) 
                 \right>\right>_{0}.
\end{eqnarray*}
Plugging this formula into (\ref{eqn:tildederiv}) and using topological
recursion relation again, we obtain
\[ \frac{\partial^{2}}{\partial t^{\mu}_{k} \partial t^{\nu}_{l}}
        \widetilde{\Psi}_{0, 2}  = 0. \]
This proves the proposition.
$\Box$


\vspace{30pt}

\noindent
Department of Mathematics \\ 
Massachusetts Institute of Technology \\
Cambridge, MA 02139 \\
USA \\

\noindent
e-mail address: 

    {\it xbliu@math.mit.edu} 

    {\it tian@math.mit.edu}

\end{document}